\documentclass[11pt,reqno]{amsart}
\usepackage{amsmath,amssymb,amsthm,mathrsfs,enumitem}
\usepackage[margin=1in]{geometry}
\usepackage{needspace}
\usepackage[colorlinks=true,linkcolor=blue,citecolor=blue]{hyperref}
\usepackage{enumitem}

\newtheorem{Theorem}{Theorem}[section]
\newtheorem{Proposition}[Theorem]{Proposition}
\newtheorem{Lemma}[Theorem]{Lemma}

\newtheorem{Corollary}[Theorem]{Corollary}
\newtheorem{Conjecture}[Theorem]{Conjecture}

\theoremstyle{definition}
\newtheorem{Definition}[Theorem]{Definition}
\theoremstyle{remark}
\newtheorem{Remark}[Theorem]{Remark}

\newtheorem{Example}[Theorem]{Example}

\def\bbC{\mathbb C}
\def\bbZ{\mathbb Z}
\def\dN{\underline N}
\def\dG{\underline G}
\def\dV{\underline V}
\def\dW{\underline W}
\def\dQ{\underline Q}
\def\dg{\underline{\mathfrak g}}
\def\scrX{\mathscr X}
\def\scrZ{\mathscr Z}
\def\scrO{\mathscr O}
\def\B{\mathrm B}
\def\v{\mathrm v}
\def\w{\mathrm w}
\def\BM{{\mathrm BM}}
\def\calA{\mathcal A}

\def\calV{\mathcal V}
\providecommand{\calO}{\mathcal O}

\providecommand{\End}{\operatorname{End}}
\def\calN{\mathcal N}
\def\calE{\mathcal E}
\def\bfU{\mathbf U}
\def\Fl{\mathrm{Fl}}
\def\DEC{\mathrm{DEC}}
\def\Spec{\operatorname{Spec}}
\def\Ker{\operatorname{Ker}}
\def\IC{\operatorname{IC}}
\def\Hom{\operatorname{Hom}}
\def\Ext{\operatorname{Ext}}
\def\Up{\operatorname{Up}}
\def\Cyc{\operatorname{Cyc}}
\def\ev{\mathrm{ev}}
\def\muhom{\mu hom}
\def\frakp{\mathfrak p}
\def\bR{\mathbf R}
\def\bT{\mathbf T}
\def\M{\mathfrak M}

\providecommand{\bfM}{\mathbf M}

\newcommand{\llb}{[\![}
\newcommand{\rrb}{]\!]}
\newcommand{\llp}{(\!(}	
\newcommand{\rrp}{)\!)}

\title{The non symmetric Hikita conjecture}
\date{}
\author{Michela Varagnolo} 
\address{\scriptsize{M.V. : ~CY Cergy Paris Universit\'e,  95302 Cergy-Pontoise, France,
UMR8088 (CNRS), ANR-24-CE40-3389 (GRAW).}}
\author{Eric Vasserot} 
\address{\scriptsize{E.V. :~Universit\'e Paris Cit\'e, 75013 Paris, France, UMR7586 (CNRS), 
ANR-24-CE40-3389 (GRAW).
}}

\begin{document}

\begin{abstract}
We construct local and global comparison maps for the Hikita
conjecture in quiver gauge theories and study their analogues
for finite non symmetric Cartan pairs. 
The construction uses quiver Hecke description
of Coulomb branch and microlocal geometry.
It yields comparison maps from  B-algebras to the cohomology of
some Lagrangians of quiver varieties. 

For arbitrary symmetric types,  the local and global comparison maps are surjective and
are isomorphisms under a B-algebra dimension bound.
In finite non symmetric types, the local image is the image of the Kirwan map, i.e., the tautological subalgebra. 
It can be proper. Under some assumptions we also
construct a non symmetric global comparison map.

Crystal theoretic methods identify the set of  nonempty Lagrangians with the set of isomorphism classes of
simple modules of the category $\scrO$ of the Coulomb branch, and 
proves that the Lagrangians are nonempty for 
$\ell$-dominant $\ell$-highest weights. 

The proof uses microlocal modules of Coulomb branches that we define.
We prove they form integral bases of the Grothendieck group of the category $\scrO$, 
with Jordan-H\"older multiplicities given by characteristic cycles. 

In symmetric
type we prove an isomorphism between based quasimap
modules at homogeneous polarized points of quiver varieties and co-specialized microlocal
modules. The proof uses hyperbolic restriction and Braden's theorem. We conjecture an extension to all
homogeneous points and a relative version.
\end{abstract}

\maketitle
\setcounter{tocdepth}{1}
\tableofcontents

\section{Introduction}

\subsection{The conjecture}
Let $(C,D)$ be a Cartan pair of finite type, with Cartan matrix
$C=(c_{ij})_{i,j\in I}$ and normalized symmetrizer
$D=\operatorname{diag}(d_i)$. We call the datum symmetric when
$D=1$, and non symmetric otherwise. Fix an orientation $Q$ of its
Dynkin graph. In \cite{VVsh} we study the category
$\scrO_\mu^\rho$ of a truncated shifted quantum affine algebra
$\calA_\mu^\rho$, realized as the quantized K-theoretic Coulomb branch with
symmetrizers of a quiver gauge theory, with integral flavour parameter $\rho$,
shift $\mu$, which is related with the rank of the gauge group,
and quantum parameter $\zeta$ which is not a root of unity, see also \cite{NW23}. 
The quantized Coulomb branch $\calA_\mu^\rho$ is graded by charges, 
which label the connected components of the affine Grassmanian of the gauge group $G$.
The B-algebra
$\B(\calA_\mu^\rho)$ is a quotient of the charge zero subalgebra of $\calA_\mu^\rho$ by
products of negative and positive charge. Its spectrum parametrizes
the simple objects of $\scrO_\mu^\rho$. 
This simple objects are also by labelled by $\ell$-highest weight of the affine quantum groups, or, equivalently, by integral points 
$\zeta^{2\gamma}$ of the gauge group,
where $\gamma$ is a cocharacter of a maximal torus $T$ of  $G$.
A combinatorial notation of red and black types to describe $\rho$
and $\gamma$ is introduced in
\S\ref{subsec:data}.

For symmetric Cartan pairs, Hikita's conjecture compares the B-algebra of the quantized homological Coulomb
branch with the equivariant cohomology of a Nakajima quiver variety
\cite{HikNotes}. Its integral flavour specialization at $\rho$ and completion
at $\gamma$ involve two maps from the completed Cartan algebra:
\[
 \widehat R_\gamma\xrightarrow{\ r_\gamma\ }
 \widehat{\B(\calA_\mu^\rho)}_\gamma,
 \quad
 \widehat R_\gamma\xrightarrow{\ \kappa_\gamma\ }
 H^*\bigl(\M_{\dQ}^\bullet(v_\gamma,w_\rho),\bbC\bigr).
\]
Here $\dQ$ is the unfolded quiver associated with the Cartan pair $(C,D)$ and the orientation $Q$,
and $\M_{\dQ}^\bullet(v_\gamma,w_\rho)$ is a graded Nakajima quiver
variety \cite{Nak01} of type $\dQ$. 
The completed Cartan algebra $\widehat R_\gamma$ is the completion of the complexified representation ring
of $G$ at the point $\zeta^{2\gamma}$.  
It is identified with a completion of the charge 0 subalgebra of $\calA_\mu^\rho$
via the dressed monopole operators.
The map $\kappa_\gamma$ is the local Kirwan map. It sends completed symmetric
polynomials in Chern roots to tautological Chern classes. The
{local Hikita conjecture} is the equality of kernels
$$\Ker( r_\gamma)=\Ker(\kappa_\gamma)$$ in \eqref{locHik}.
Constructing a local comparison map 
\begin{align}\label{Phigamma}
 \Phi_\gamma:\ \widehat{\B(\calA_\mu^\rho)}_\gamma
 \longrightarrow
 H^*\bigl(\M_{\dQ}^\bullet(v_\gamma,w_\rho),\bbC\bigr),
 \quad \Phi_\gamma\circ r_\gamma=\kappa_\gamma.
\end{align}
first requires
proving the inclusion $\Ker (r_\gamma)\subseteq\Ker(\kappa_\gamma)$.
The conjecture says that the resulting map $\Phi_\gamma$ is an isomorphism.
The global conjecture compares the B-algebra
$\B^\theta_{F,\hbar}$ of the quantized Coulomb branch, with full
flavour torus $F=T_W$ and a generic charge cocharacter $\theta$, to
$H^*_{F\times\bbC^\times}(\M_Q(v,w),\bbC)$, see \eqref{GHC}.
Here $\theta$ is identified with  the stability parameter of the quiver variety, 
which is assumed to be the framing generated one.
For non symmetric Cartan pairs there is no ordinary
Higgs branch for the Coulomb theory with symmetrizers
\cite[\S1]{NW23}. We therefore formulate an analogue of the Hikita conjecture
using stable Lagrangian loci in the unfolded representation space.

\subsection{ADE type}
Our first result is a construction of the local comparison map that
does not use quasimaps. For $Q$ of type ADE,
Theorem~\ref{thm:main}(a),(b) gives a surjective homomorphism
$ \Phi_\gamma$.
By
Theorem~\ref{thm:main}(c), the map is an isomorphism if the
classical B-algebra dimension bound \eqref{inequality} holds. This bound
is available for types $A,D,E_6$ and, with the framing restrictions
\eqref{FR}, for types $E_7,E_8$ \cite{KTWWY19a,DK25}.
Independently of this bound,
Corollary~\ref{cor:points} proves
\[
 \M_{\dQ}^\bullet(v_\gamma,w_\rho)\ne\emptyset
 \quad\Longleftrightarrow\quad
 \zeta^{2\gamma}\in\Spec\B(\calA_\mu^\rho).
\]
The same local maps, together with specialization at a Zariski dense
set of integral flavour parameters, yield a surjective equivariant
global comparison. Under \eqref{inequality},
Corollary~\ref{cor:global} identifies it with the global Hikita
isomorphism.

The existence of the map $\Phi_\gamma$ in finite simply laced
type is accessible through affine Grassmannian slices
\cite{KTWWY19a}. The point here is a different, microlocal construction. 
It does not use quasimaps.
Indeed, the localization theorem of \cite{VV25},
Theorem~\ref{thm:V} here, identifies integral blocks of the Coulomb
algebra with a completion of the equivariant Yoneda algebra
\[
 \bR=\bigoplus_{\eta,\nu}
 \Ext^*_{\dG}(\mathscr L_\eta,\mathscr L_\nu)
\]
of flag-incidence complexes. For each type
$\gamma$, a smooth stable conormal locus $Z_\gamma^{\rm st}$ is introduced,
and the sheaf $\mathscr L_\gamma$ is microlocally simple along it.
The restriction of microlocal Hom therefore defines an algebra map
\[
 \Psi_\gamma:e_\gamma\bR e_\gamma
 \longrightarrow H^*_{\dG}(Z_\gamma^{\rm st},\bbC).
\]
The microsupport estimate
$\operatorname{SS}(\mathscr L_\eta)\subseteq Z_\eta$ shows that
$\Psi_\gamma$ kills the types higher than $\gamma$ and the cyclotomic one, 
defining a morphism from the local B-algebra. Its restriction
to the Cartan algebra coincides with the local Kirwan map.
We do not need to compute the full microsupport of the flag-incidence complex $\mathscr L_\eta$.
When $D=1$, an algebraic contraction identifies the
cohomology of the stable conormal locus with that of the graded quiver
variety. The Kirwan surjectivity for Nakajima varieties
\cite{MN18} then proves the surjectivity in \eqref{Phigamma}.

\subsection{Symmetric type}
The construction extends beyond Dynkin type. In
\S\ref{sec:general} we allow an arbitrary finite quiver $Q$ with
$D=1$, including loops and multiple arrows. In place of the
finite type shifted quantum affine algebra, we use the quantized
Coulomb branch of $Q$. An homological localization and the
quiver Hecke identification of its integral blocks in
Lemma~\ref{VQ} provide the same completed B-algebras and
microlocal construction. Theorem~\ref{thm:gen-surj} gives the
surjective local comparison maps, and
Theorem~\ref{thm:gen-global-map} gives a surjection
\[
 \Phi_{F,\hbar}:\ \B^\theta_{F,\hbar}
 \twoheadrightarrow H^*_{F\times\bbC^\times}
 (\M_Q(v,w),\bbC).
\]
Another approach to the symmetric case using quasimaps  appears in \cite{BT26}.
For general $Q$, the injectivity is a separate question: The sufficient
condition used in Theorem~\ref{thm:gen-inj} is
\begin{align}\label{inequality}
 \dim\B^\theta(A_0)
 \leqslant \dim H^*(\M_Q(v,w),\bbC),
\end{align}
where $A_0$ is the classical Coulomb algebra at zero flavour and
quantization parameters. Under this bound all local maps and the
global map are isomorphisms. Besides the finite Dynkin cases above,
the bound is known for the Jordan quiver, including general framing
\cite{Hik17,KS24}. 
A general proof of this bound using  IA was recently proposed, see \cite{HikNotes} for details.

\subsection{Non symmetric types}
For finite non symmetric Cartan pairs, the graded retraction used in
the symmetric proof is not generally available. The relevant
replacement for the graded quiver variety is the quotient
\begin{align}\label{LRG1}
 \mathfrak L_{\rho,\gamma}
 =Z_\gamma^{\rm st}/\dG.
\end{align}
Its points are stable representations of the unfolded data whose
minimal $\bbZ$-flag has type $\gamma$. If nonempty, this is a
smooth irreducible quasi-projective Lagrangian subvariety of the
Nakajima variety for $\dQ$, see Lemmas~\ref{lem:stable-irreducible},
\ref{F}. Its cohomology is finite dimensional, but need not be
concentrated in even degrees. The local analogue of the Hikita map
is
\begin{align}\label{Phigammans}
 \Phi^{\rm ns}_\gamma:\ 
 \widehat{\B(\calA_\mu^\rho)}_\gamma
 \longrightarrow H^*_{\dG}(Z_\gamma^{\rm st},\bbC)
 \cong H^*(\mathfrak L_{\rho,\gamma},\bbC).
\end{align}
By Theorem~\ref{thm:ns-map}, it is defined for every integral pair
$(\rho,\gamma)$, and its image equals the image of the local Kirwan
map $\kappa^{\rm ns}_\gamma$. Hence, we have
\[
 \operatorname{Im}\Phi^{\rm ns}_\gamma
 =\operatorname{Im}\kappa^{\rm ns}_\gamma
 =H^*(\mathfrak L_{\rho,\gamma},\bbC)_{\rm taut}.
\]
It is the subalgebra generated by tautological Chern classes. 
The same theorem gives
\[
 \Phi^{\rm ns}_\gamma\ne0
 \quad\Longleftrightarrow\quad Z_\gamma^{\rm st}\ne\emptyset
 \quad\Longleftrightarrow\quad
 \zeta^{2\gamma}\in\Spec\B(\calA_\mu^\rho).
\]

The tautological image may be proper.
Proposition~\ref{prop:B2-local-kirwan} exhibits a type $B_2$
example with
\[
 \mathfrak L_{\rho,\gamma}\cong\bbC^\times\times\mathbb A^2,
 \quad
 \operatorname{Im}\Phi^{\rm ns}_\gamma=\bbC
 \subsetneq H^*(\mathfrak L_{\rho,\gamma},\bbC).
\]
The missing class has cohomological degree one. 
The non symmetric Kirwan map is not always
surjective. We also describe $\mathfrak L_{\rho,\gamma}$ as a
critical locus, in Proposition~\ref{prop:ns-crit}, and compare this
model with the graded triple quiver critical locus considered in other works, e.g., in  \cite{VV23}.

A global comparison construction is given in
Appendix~\ref{app:ns-global-technical}. We use locally free
moduli spaces $\mathfrak M_{N_0}$ associated with the symmetrizer.
Proposition~\ref{prop:ns-global-torsion} constructs the global comparison map
\[
 \Phi^{\rm ns}_{\rm tf}:\B^\theta_{F,\hbar}
 \longrightarrow
 H^*_{F\times\bbC^\times}(\mathfrak M_{N_0},\bbC)
 /\operatorname{tor}.
\]
Its image is tautological. If the target is torsion free over $H^*_{F\times\bbC^\times}$ then the map lifts to the equivariant
cohomology, and it is surjective when the global Kirwan
map is. This condition cannot be suppressed:
Proposition~\ref{prop:G2-global-kirwan} gives a type $G_2$
example with nonzero odd cohomology for which the global
Kirwan surjectivity fails.

\subsection{Application to representation theory}
The first application is an existence and classification result for highest types, i.e., 
highest $\ell$-weights of simple modules of the category $\scrO_\mu^\rho$. 
For every finite Cartan pair and every weight block
$\Gamma$ of $\scrO_\mu^\rho$, the Kashiwara-Saito geometry of
preprojective algebra crystals \cite{KS97} implies in Proposition~\ref{prop:microlocal-basis}(a) that
\[
Z_\gamma^{\rm st}\ne\emptyset,\quad\forall \gamma\in \Gamma.
\]
This identifies the set of types $\gamma$ with nonempty stable locus $Z_\gamma^{\rm st}$
with the set of highest types or, equivalently, the set of highest $\ell$-weights of the 
simple objects in $\scrO_\mu^\rho$.
Independently, Theorem~\ref{thm:ns-existence} proves that every
$\ell$-dominant type $\gamma$ has
$Z_\gamma^{\rm st}\ne\emptyset$.
Together with Theorem~\ref{thm:ns-map}(d) and the character criterion
of Lemma~\ref{lem:Bcrit}, this proves
Theorem~\ref{Theorem:431}: for each $\ell$-dominant $\ell$-weight
$Z_{\rho,\gamma}$, there is a unique algebra character
$\varepsilon:\calA_\mu^\rho\to\bbC$ restricting to evaluation at
$\zeta^{2\gamma}$ on the Cartan algebra, and its one dimensional
representation is the simple module $L(\rho,\gamma)$. 

Our second application is the construction of Morse and microlocal
modules over the completed Coulomb algebra, in
\S\S\ref{sec:morse}-\ref{sec:micro}.  Theorems~\ref{thm:morse-JH} and
\ref{thm:microlocal-JH} show that the relevant specialized modules
have finite length in $\scrO_\mu^\rho$, with a unique copy of
$L(\rho,\gamma)$ at the highest type. 
Here, highest refers to the charge type. The other weight pieces
are controlled by the support order.
Their Jordan-H\"older multiplicities are expressed by vanishing cycles or microlocal Hom
and, at generic points of $Z_\gamma^{\rm st}$, by
characteristic cycle multiplicities of the simple summands of the
flag-incidence complexes. Corollary~\ref{cor:microlocal-morse-K0}
identifies the generic ungraded Grothendieck classes of the Morse
and specialized microlocal modules. By
Proposition~\ref{prop:microlocal-basis}(b), these classes
form an integral unitriangular basis of the Grothendieck group of
each finite weight block. 

\subsection{Quasimap modules}
The microlocal construction suggests a relation with the
Coulomb branch modules defined by Botta and Tamagni via the
critical cohomology of based quasimaps \cite{BT26}. In
\S\ref{sec:QMM} we conjecture a comparison with the
\emph{co-specialized} microlocal module
\[
 M^!(\rho,\gamma)^\xi
 \cong\bigoplus_\eta H^*\bigl((i_\xi)^!\mathcal K_\eta[2n]\bigr),
 \quad
 \mathcal K_\eta
 =\muhom(\mathscr L_\gamma,\mathscr L_\eta)
      |_{Z_\gamma^{\rm st}}.
\]
The exceptional restriction is essential.
The conjecture is formulated for homogeneous stable points in
symmetric types. Such a point lies in a linear Lagrangian
$L_{F_\gamma}$, but this does not by itself make its image $\mathfrak L_{\rho,\gamma}$ in the
Nakajima variety polarized in the sense of \cite[\S1.3.5]{BT26}, i.e., BT-polarized.
For BT-polarized homogeneous stable points, Theorem~\ref{thm:BT-polarized}
proves the  isomorphism of modules over the quantized Coulomb branch. The
proof uses Braden's
hyperbolic restriction theorem \cite{Bra03}.
A partial Fourier transform gives the change of polarization.
 Without the BT-polarization the isomorphism
remains conjectural. 

\subsection{Plan}
In \S\ref{sec:recall} we recall the Coulomb, B-algebra and
quiver-Hecke ingredients. In \S\ref{sec:cotangent} we construct
the conormal and stable Lagrangian loci, and in
\S\ref{sec:sym} the contraction for $D=1$. The type ADE
comparison is stated in \S\ref{sec:main}, while the Morse and
microlocal modules are developed in
\S\S\ref{sec:morse}-\ref{sec:micro}. The proof of the
ADE comparison is completed in \S\ref{sec:surj}, while 
\S\ref{sec:general} treats arbitrary symmetric quivers,
and \S\ref{sec:ns} finite non symmetric types.
In \S\ref{sec:QMM} we prove the quasimap comparison at polarized
homogeneous stable points and formulate conjectural extensions. The appendices develop
the non symmetric global construction, the critical locus models and many examples. 
Sections~\ref{sec:sym}, \ref{sec:main}
and \ref{sec:surj} assume the type is ADE unless stated
otherwise. The Morse and microlocal results apply to all types.

\subsection{Conventions}
All varieties are over $\bbC$ and cohomology has coefficients in
$\bbC$ unless otherwise specified. We write
$\bullet=\Spec\bbC$ and
$H_G^*=H_G^*(\bullet,\bbC)$ for an algebraic group $G$.

\section{Coulomb branches, B-algebras and quiver Hecke algebras}\label{sec:recall}

\subsection{Cartan pairs, reds and blacks, $\ell$-weights}\label{subsec:data}
A Cartan pair is a pair $(C,D)$ where $C=(c_{ij})_{i,j\in I}$ is the incidence matrix
of a finite quiver $Q$ with vertex set $I$, and $D=(d_i)_{i\in I}$ is a symmetrizer with
$d_ic_{ij}=d_jc_{ji}$ and $\min \{d_i\}_{i\in I}=1$. 
Let $H$ denote the graph underlying the quiver $Q$.
The pair is of symmetric type if $D=1$.
It is of finite type if $C$ is a Cartan matrix,
equivalently, if $DC$ is positive definite. 
Fix $\zeta\in\bbC^\times$, not a root of unity.

Let $V=\bigoplus_{i\in I}V_i$ and $W=\bigoplus_{i\in I}W_i$ be finite dimensional $I$-graded vector spaces, of dimension vectors
$v=(v_i)$ and $w=(w_i)$, and put
\[
G_V=\prod_iGL(V_i),\quad G_W=\prod_iGL(W_i),
\]
with maximal tori $T_V\subset G_V$ and $T_W\subset G_W$ and Weyl groups $\mathfrak S_V=\prod_i\mathfrak S_{v_i}$
and  $\mathfrak S_W$. The character lattice
of $T_V$ is $X^*(T_V)=\bigoplus_i\bbZ^{v_i}$, and the complexified representation ring of $G_V$ is 
$$R(G_V)=\bbC[X^*(T_V)]^{\mathfrak S_V},$$ the ring of
symmetric Laurent polynomials in the Chern roots 
$$x_{i,1},\dots,x_{i,v_i},\quad i\in I.$$ Its complex points are
\[
\Spec\big(R(G_V)\big)=T_V/\mathfrak S_V ,
\]
a point being a $\mathfrak S_V$-orbit $[t]$ of elements $t\in T_V$, with character 
$$\ev_t:f\mapsto f(t).$$
When no confusion is possible, we abbreviate
$$G=G_V,\quad T=T_V$$

An element $s=(s_V,s_W)$ in $T_V\times T_W$ or its image in $(T_V\times T_W)/(\mathfrak S_V\times\mathfrak S_W)$ 
is {integral} if all its eigenvalues are integral powers of $\zeta$. Such an $s$ defines the {$s$-grading}
\begin{align}\label{VW}
V_i=\bigoplus_{a\in\bbZ}V_i[a],\quad V_i[a]=\Ker(s-\zeta^a)|_{V_i},\quad W_i=\bigoplus_{l\in\bbZ}W_i[l],\quad W_i[l]=\Ker(s-\zeta^l)|_{W_i}.
\end{align}
A vector is {$s$-homogeneous} of degree $a$ if it lies in some $V_i[a]$ or $W_i[a]$.

\begin{Definition}\label{def:blacks}
Let $s=(s_V,s_W)$ be integral.
Fix bases of $V_i$, $W_i$  of $s$-homogeneous vectors.
\begin{enumerate}[label=$\mathrm{(\alph*)}$,leftmargin=8mm]
\item 
A \emph{black} is a vector of the basis of $V$. A black in $V_i[a]$ is recorded as the pair $(i,a)$.
Its {vertex} is $i$ and its {uncentred position} is $a$. 
\item 
A \emph{red} is a vector of the basis of $W$. A red in $W_i[l]$ is recorded as the pair $(i,l)$. 
Its vertex is $i$ and its uncentred position is $l$.
\item 
The multiset of black pairs $(i,a)$ is denoted $\gamma$ and the multiset of red pairs $(i,l)$ is denoted  $\rho$. 
The black $(i,a)$ occurs in $\gamma$ with multiplicity $\dim V_i[a]$, and the red $(i,l)$ occurs in $\rho$ with multiplicity $\dim W_i[l]$.
\item
A \emph{type} of dimension $v\in\mathbb N I$ 
is a multiset $\eta$ of pairs $(i,a)$ with $v_i$ pairs of vertex $i$ for each $i\in I$.
We also call its elements blacks.
A type $\gamma$ determines the integral point $\zeta^{2\gamma}$ in $T_V/\mathfrak S_V$ 
with one eigenvalue $\zeta^a$ for each black $(i,a)$, and the evaluation character
\begin{align}\label{EVG}\ev_\gamma:R(G_V)\to\mathbb C,\quad f\mapsto f\big((\zeta^a)_{(i,a)\in\gamma}\big).\end{align}
\end{enumerate}
\end{Definition}

Since $\zeta$ is not a root of unity, types, integral Weyl-group
orbits, and their evaluation characters are equivalent data.
In the conventions of \cite{VVsh} we have
$a=2d_i\gamma_{i,r}$ and $l=2d_i\rho_{i,s}$.
The localization parameter is
$$s_{\rho,\gamma}=(\zeta^{2\gamma},\zeta^{2\rho})\in T_V\times T_W.$$

\begin{Definition}\hfill
\begin{enumerate}[label=$\mathrm{(\alph*)}$,leftmargin=8mm]
\item
We define the \emph{centred positions} of a  red $(i,l)$ or black $(i,a)$ to be the integers
\[
L=l+d_i,\quad P=a+d_i.
\]
\item
The \emph{colour} of a red or black of vertex $i$ is the pair $(i,r)$, where $r$ is the class of its centred position modulo $2d_i$. 
We may call $r$ the residue.
Two blacks of the same vertex have the same colour iff their $s$-eigenvalues differ by an integral power of $\zeta^{2d_i}$.
In the symmetric case ($D=1$) the colour is the vertex together with the parity of the position. 
\end{enumerate}
\end{Definition}

The $\ell$-weight attached to $(\rho,\gamma)$ is the Laurent monomial in the variables $Z_{i,k}$ given by
\begin{align}\label{ZRG}
Z_{\rho,\gamma}=\prod_{(i,l)\in\rho}Z_{i,l+2d_i}\ \prod_{(i,a)\in\gamma}A_{i,a+d_i}^{-1},\quad
A_{i,a+d_i}=Z_{i,a}Z_{i,a+2d_i}\prod_{j\neq i}\prod_{n=1}^{|c_{ij}|}Z_{j,a+d_ic_{ij}+2d_in}^{-1}
\end{align}
It is {$\ell$-dominant} if all its $Z_{i,k}$-exponents are $\geqslant0$.
We have $Z_{\rho,\gamma}=Z_\rho A_\gamma^{-1}$ in the notation of \cite{VVsh}.

\subsection{B-algebras}\label{subsec:algebras}
The B-algebra of a $\bbZ$-graded algebra $A=\bigoplus_{n\in\bbZ}A_{[n]}$ is the algebra
\[
\B(A)=A_{[0]}\Big/\sum_{n>0}A_{[-n]}A_{[n]}.
\]
 In the rest of this section we apply this construction to the shifted quantum affine groups,
 assuming that the Cartan pair
$(C,D)$ is of finite type. 
For any coweight $\mu$, let $\bfU_\mu$ be the corresponding shifted affine quantum group, with
the Drinfeld generators $x^\pm_{i,n}$, the Cartan generators $\psi^\pm_{i,\pm n}$, and triangular decomposition
$$\bfU_\mu=\bfU^-_\mu\otimes\bfU^0_\mu\otimes\bfU^+_\mu.$$
 The shifted affine quantum group is graded by the root lattice. 
 Let $\bfU_{\mu,\beta}$ be the
piece of weight $\beta$ and $\bfU_{\mu,>0}$, $\bfU_{\mu,<0}$ the sums over positive and negative degrees. 

The truncated shifted affine quantum group $\calA^\rho_\mu$ is the quotient 
$$\Psi:\bfU_\mu\twoheadrightarrow\calA^\rho_\mu.$$ 
given by the quantized K-theoretic Coulomb branch 
with symmetrizers used in \cite{VV25}, \cite{VVsh} specialized at the flavour parameter $\zeta^{2\rho}$ and the quantum parameter 
$q=\zeta$,  with $\zeta$ not a root of unity. 
It allows us to identify the category $\scrO$ of $\calA_\mu^\rho$ with a subcategory $\scrO_\mu^\rho$ of the category
$\scrO$ of the shifted affine quantum group $\bfU_\mu$.  
We only consider here the geometric truncation given by K-theoretic Coulomb branch, 
not the analogous algebraic truncation considered in \cite{H23}.
These truncations are known to be different.
Let 
$$R_G=\Psi(\bfU^0_\mu)=\calA^{\rho,0}_\mu.$$
We call $R_G$ the \emph{Cartan algebra}. It is the ring of symmetric Laurent polynomials in
the Chern roots of the gauge group $G$.
The set of integral points of $R_G$ is $\{\zeta^{2\gamma}\,:\,\gamma\text{ is a type}\}$.

We consider now the B-algebra $\B(A)$ of the truncated shifted affine quantum group $\calA^\rho_\mu$ with
$$A=\calA^\rho_\mu,\quad A_\beta=\Psi(\bfU_{\mu,\beta}).$$
The $\mathbb Z$-grading comes from a cocharacter $\theta$ of
 the torus $(\mathbb C^\times)^I$, see \S\ref{subsec:global-conjecture}.
For a cocharacter positive on the simple roots, i.e., a dominant cocharacter, the triangular
decomposition gives the presentation
\[
\B(A)=A_0\big/(A_{<0}A_{>0})\cap A_0.
\]
Let $p_0:\bfU_{\mu,0}\to\bfU^0_\mu$ be the projection along $\bigoplus_{\beta>0}\bfU^-_{-\beta}\bfU^0_\mu\bfU^+_\beta$. 

\begin{Lemma}\label{lem:Bpres}\hfill
\begin{enumerate}[label=$\mathrm{(\alph*)}$,leftmargin=8mm]
\item $\B(A)\cong R_G/J$ and $J=\Psi\big(p_0(\Ker\Psi\cap\bfU_{\mu,0})\big).$ 
\item The map $q:R_G\to\B(A)$ is surjective. 
\end{enumerate}
\end{Lemma}

\begin{proof}
Put $K=\Ker\Psi$. Then $A_0=\bfU_{\mu,0}/K_0$. 
The map $p_0$ is the projection along the space $N=(\bfU_{<0}\bfU_{>0})_0$. 
Hence 
$$\B(A)=\bfU_{\mu,0}/(N+K_0)\cong\bfU^0_\mu/p_0(K_0),$$ and
$\Psi$ maps $\bfU^0_\mu$ onto $R_G$ with kernel $K\cap\bfU^0_\mu\subset p_0(K_0)$.
\end{proof}

We call \emph{Cartan quotient} the surjective map in Lemma \ref{lem:Bpres}
\begin{align}\label{CQ} q:R_G\twoheadrightarrow \B(A)
\end{align}
The points of $\Spec\B(A)$ classify simple objects of $\scrO_\mu^\rho$ via the simple heads of Verma modules. 
Let $L(\rho,\gamma)\in\scrO_\mu^\rho$ 
be the simple $\bfU_\mu$-module of highest $\ell$-weight given by the Laurent monomial $Z_{\rho,\gamma}$
defined in \eqref{ZRG}.

\begin{Lemma}\label{lem:onedim}
A one dimensional $\bfU_\mu$-module
with $\ell$-weight $Z_{\rho,\gamma}$ exists if and only if $Z_{\rho,\gamma}$ is $\ell$-dominant.
Then it is $L(\rho,\gamma)$.
\end{Lemma}

\begin{proof}
A character $\phi$ of $\bfU_\mu$ kills the $x^\pm_{i,n}$: the relation 
$$\psi^+_ix^\pm_{i,n}=\zeta_i^{\pm1}x^\pm_{i,n}\psi^+_i$$ with
$\phi(\psi^+_i)\neq0$ and $\zeta_i\neq1$ forces $\phi(x^\pm_{i,n})=0$. So $\phi$ is determined by an $\ell$-weight, and the only
remaining relation says that the expansions of the $\ell$-weight at $u=\infty$ and $u=0$ agree, i.e., that it is a
polynomial in $u^{-1}$. For $Z_{\rho,\gamma}$ this is $\ell$-dominance, see also \cite[Thm.~5.5]{H23}.
\end{proof}

Let $\bfM^\eta(f)\in\calA^\rho_\mu$ be the dressed monopole operator of charge the minuscule coweight $\eta$ of the gauge group $G$.
Recall that the assignment $f\mapsto \bfM^0(f)$ yields an algebra isomorphism
$$R(G)\cong R_G.$$
One of the results of this paper is the following theorem.

\begin{Theorem}\label{Theorem:431}
Assume that $Z_{\rho,\gamma}$ is $\ell$-dominant.
\begin{enumerate}[label=$\mathrm{(\alph*)}$,leftmargin=8mm]
\item
 There is a unique character $\varepsilon:\calA^\rho_\mu\to\bbC$ with
$\varepsilon|_{R_G}=\ev_\gamma$. 
\item
The character $\varepsilon$ vanishes on the operators $\bfM^\eta(f)$ with $\eta\neq0$ minuscule.
\item
We have 
$\varepsilon^*\bbC\cong L(\rho,\gamma).$
\qed
\end{enumerate}
\end{Theorem}

The theorem is proved in \S\ref{sec:ns}.
Indeed, by Theorem~\ref{thm:ns-existence}, the $\ell$-dominance of $Z_{\rho,\gamma}$ implies that
the set $Z_\gamma^{\rm st}$ (defined below) is non empty.
Theorem~\ref{thm:ns-map}(d) then gives
$\zeta^{2\gamma}\in\Spec\B(\calA^\rho_\mu)$ and 
the assertions of the theorem follow from the following lemma that we quote here for a future use.

\begin{Lemma}\label{lem:Bcrit}
Assume that $Z_{\rho,\gamma}$ is $\ell$-dominant. Then
$$\text{Theorem $\ref{Theorem:431}$ holds for $(\rho,\gamma)$}\iff \zeta^{2\gamma}\in\Spec\B(\calA^\rho_\mu).$$
\end{Lemma}

\begin{proof}
Theorem $\ref{Theorem:431}$ says that the one dimensional module 
$L(\rho,\gamma)$ of $\bfU_\mu$ in Lemma \ref{lem:onedim} descends to $\calA^\rho_\mu$.
A simple module of $\scrO^\rho_\mu$ with highest weight $\zeta^{2\gamma}$ exists iff $\zeta^{2\gamma}\in\Spec\B$, 
and it is then $L(\rho,\gamma)$. 
\end{proof}

\subsection{The unfolded data}\label{subsec:unfolded}
Fix a pair of types $(\rho,\gamma)$. 
Following \cite{VV25}, \cite{VVsh}, we consider
the fixed point set $\dN$ of the localization parameter $s_{\rho,\gamma}=(\zeta^{2\gamma},\zeta^{2\rho})$
in the regular representation space (defined below).
It is described in the following way.
For each colour $c=(i,r)$ let $\dV_c$ be the vector space with basis
the blacks of colour $c$, and $\dW_c$ the one with basis the reds of colour $c$. Put 
$$\dV=\bigoplus_c\dV_c,\quad
\dW=\bigoplus_c\dW_c,\quad
\dG=\prod_cGL(\dV_c),\quad
\dg=\operatorname{Lie}\dG$$
and
$$v_c=v_c(\gamma)=\dim \dV_c,\quad w_c=w_c(\rho)=\dim\dW_c.$$
We have
\[
\dN=\bigoplus_c\Hom(\dW_c,\dV_c)\ \oplus\bigoplus_{h:i\to j}\ \bigoplus_{(c,c')}\Hom(\dV_c,\dV_{c'}).
\]
The second sum runs over the arrows $h$ of $Q$ and the pairs of colours 
\begin{align}\label{colour1}c=(i,r),\quad c'=(j,s)\end{align} 
where $d_{ij}=\gcd(d_i,d_j)$ and
\begin{align}\label{parity}
r-s\equiv d_{ij}\pmod{2d_{ij}}
\end{align}
A point of $\dN$ is a tuple
$$x=(A,B),\quad
A=(A_c),\quad
B=(B_h).$$ 
In the symmetric case $d_{ij}=1$, hence $c'$ is the colour of opposite parity.

Let $\dQ$ be the quiver whose vertices are the colours with one arrow
$c\to c'$ for each arrow $h:i\to j$ of $Q$ and each pair of colours as in \eqref{colour1}.
We call it the \emph{unfolded quiver}. 
Then $\dN$ is the representation space of $\dQ$ framed by $\dW$, and
$\dG$ is its gauge group.

Above, the expression \emph{fixed points} refers to the combined $s_{\rho,\gamma}$ and twisted loop
rotation action of \cite{VV25}, \cite{VVsh} on the \emph{Laurent representation space} 
$$N_K=\bigoplus_i\Hom(W_i,V_i)\llp z^{d_i}\rrp
 \ \oplus\!\bigoplus_{h:i\to j}
 \Hom(V_i,V_j)\llp z^{d_{ij}}\rrp$$
The \emph{regular representation space} is
\[
 N_\calO=\bigoplus_i\Hom(W_i,V_i)\llb z^{d_i}\rrb
 \ \oplus\!\bigoplus_{h:i\to j}
 \Hom(V_i,V_j)\llb z^{d_{ij}}\rrb z^{d_i+b_{ij}},
 \quad b_{ij}=d_ic_{ij}.
\]

Now, we explain formula \eqref{parity}. 
An element of $N_\calO$ is a sum of elements of the form
\begin{align*}
 &Bz^m,\quad B:V_i[a]\longrightarrow V_j[a'],\quad
 m=d_i+b_{ij}+kd_{ij},\quad k\geqslant0,\\
 &Az^m,\quad  A:W_i[l]\to V_i[a],\quad m=kd_i,\quad k\geqslant 0
 \end{align*}
In the convention of \cite{VV25}, \cite{VVsh}
loop rotation is $z\mapsto t^2z$, together with multiplication by
$t^{2\sigma_i}$ on the spaces at vertex $i$, where
\[
 \sigma_i-\sigma_j-d_i+d_j=\frac12 b_{ij}.
\]
The gauge action on $B$ contributes $t^{a'-a}$, the monomial $z^m$
contributes $t^{2m}$, and the vertex twists contribute
$t^{2(\sigma_j-\sigma_i)}$.
Thus $Bz^m$ is fixed precisely when
\[
 a'-a+2m+2(\sigma_j-\sigma_i)=0.
\]
Passing to the centred positions $P=a+d_i$ and $Q=a'+d_j$, we obtain
\[
 \begin{aligned}
 P-Q=d_i+d_j+b_{ij}+2kd_{ij}.
 \end{aligned}
\]
For a finite type edge we have
$-b_{ij}=\max(d_i,d_j)$ and $d_{ij}=\min(d_i,d_j)$. Hence
\begin{align}\label{CP}P-Q=(2k+1)d_{ij}.
\end{align}
Since $P\equiv r\pmod{2d_i}$ and $Q\equiv s\pmod{2d_j}$, reduction
modulo $2d_{ij}$ gives  \eqref{parity}.
For a framing coefficient
$Az^{kd_i}$, fixedness gives
$L-P=2kd_i$,
Regularity means $k\geqslant0$, i.e.,  $Q<P$ for $B$ and $P\leqslant L$ for $A$.
For $B$ the exponent of $z$ itself may be negative because of
the shift $d_i+b_{ij}$. 

Evaluation $N_K\to N$ identifies the fixed Laurent space with $\dN$, and its
fixed regular subspace with $\dN_{F_\gamma}$, defined below
\cite[(3.1.60), (3.1.66)-(3.1.70)]{VV25}.

\begin{Example}\hfill
\begin{enumerate}[label=$\mathrm{(\alph*)}$,leftmargin=8mm]
\item
If $D=1$, then an arrow $i\to j$ gives
$(i,r)\to(j,r+1)$, $r\in\bbZ/2$. A loop gives two arrows between
the two parities, so $\dQ$ is loopless. 
\item
If $Q$ is loopless and
bipartite with bipartition $\epsilon$, the invariant $r+\epsilon(i)$
splits $\dQ$ into two copies of $Q$. Bipartiteness is essential:
the parity cover of a triangle is a connected six-cycle.
\end{enumerate}
\end{Example}

A {$\bbZ$-flag} $F$ of $\dV$ is a family of subspaces $F_c^{\leqslant P}\subset\dV_c$, for each colour $c$ and each
position $P$ in the class of $c$, increasing in $P$, equal to $0$ for $P\ll0$ and to $\dV_c$ for $P\gg0$. 
A  $\bbZ$-flag is not a flag in the usual sense: it is indexed by integers, 
and infinitely many consecutive terms coincide. 
Forgetting the indexing, a $\bbZ$-flag  maps to a partial flag of each colour space $\dV_c$. We write
$$F^{\leqslant P}=\bigoplus_c F_c^{\leqslant P}
,\quad
F_c^{<P}=F_c^{\leqslant P-2d_i}$$ 
for a colour $c$ of vertex $i$. 
Hence $F$ may be viewed as the sequence $(F^{\leqslant P})$ where $P$ runs in $\mathbb Z$.
The \emph{type} of the $\mathbb Z$-flag $F$ is the type
 having $\dim F_c^{\leqslant P}/F_c^{<P}$ blacks of colour $c$ at centred position $P$. 
For each type $\gamma$, the
\emph{standard $\bbZ$-flag} $F_\gamma$ of type $\gamma$ is 
\begin{align}\label{SZF}F_\gamma^{\leqslant P}= \text{the\ span\ in $\dV$\ of\ the\ blacks\ of $\gamma$\ of\ positions\ }\leqslant P\end{align}

Extend the grading on $V_i=\bigoplus_p V_i[p]$ to $V_i\llp z_i\rrp$ via the convention
\begin{align}\label{SGC}\deg(vz_i^n)=p+2d_i n,\quad\forall v\in V_i[p].\end{align}
Hence, the multiplication by $z_i$ has degree $2d_i$.
We have the following lattice/flag correspondence.

\begin{Lemma}[\cite{VV25}, \cite{VVsh}]\label{lem:semilinear-flags}\hfill
\begin{enumerate}[label=$\mathrm{(\alph*)}$,leftmargin=8mm]
\item
The $k$th homogeneous subspace of $V_i\llp z_i\rrp$ for the grading \eqref{SGC} is identified with the
colour space $\dV_{(i,k\bmod 2d_i)}$ under the linear map such that
$v\otimes z_i^n\mapsto v$. 
\item
The $s_{\rho,\gamma}$-fixed lattices correspond
bijectively to the $\bbZ$-flags of $\dV$ under $L\mapsto F$,
where $F^{\leqslant P}$ is the degree $P$ part $L_{i,P}$ of $L_i$.
\item
The standard lattice $V\otimes\calO$ corresponds to $F_\gamma$,
and the type of $F$ is the multiset of degrees of a homogeneous
$\calO$-basis of $L$.
\end{enumerate}
\qed
\end{Lemma}

Fixing the numbers of blacks of each colour fixes the unfolded dimension vector $v=(v_c)_c$.
Let 
$$C_v=\{\eta\,:\,\eta\text{ has $v_c$ blacks of colour }c\}
,\quad
\Gamma_v=\{\eta\in C_v\,:\,L(\rho,\eta)\in\scrO_\mu^\rho\}$$
Then, we consider the full subcategory $\mathscr C_v\subset\scrO_\mu^\rho$ given by
\begin{align*}
\mathscr C_v&=\{M\in\scrO_\mu^\rho\,:\,\text{the composition factor of }M\text{ are }L(\rho,\eta)\,,\,\eta\in\Gamma_v\},\\
&=\{M\in\scrO_\mu^\rho\,:\,e_\eta M=0\,,\,\forall \eta\notin C_v\}.
\end{align*}
In finite types, under the decategorification isomorphism of \cite{VV25}, the complexified Grothendieck group $K(\mathscr C_v)$
identifies with a weight subspace of the representation in $\bigoplus_\mu K(\scrO_\mu^\rho)$ of the lie algebra of type $\dQ$.
We call $C_v$, $\Gamma_v$ and $\mathscr C_v$ \emph{weight blocks}.
If $v_c=v_c(\gamma)$ we may write 
\begin{align}\label{CCG}
C_\gamma=C_v,\quad
\Gamma_\gamma=\Gamma_v,\quad
\mathscr C_\gamma=\mathscr C_v.
\end{align}

For a type $\eta$ let $\Fl_\eta$ be the variety of
$\bbZ$-flags of type $\eta$.
It is isomorphic to the homogeneous space $\dG/P_\eta$, where $P_\eta$ is the stabilizer of $F_\eta$ in $\dG$. 
A point $x=(A,B)$ of $\dN$ {preserves} a $\mathbb Z$-flag $F$ if
\begin{equation}\label{eq:NF}
A(\dW^{\leqslant L})\subset F^{\leqslant L},\quad B_h(F_i^{\leqslant P})\subset F_j^{<P},\quad \forall h:i\to j,
\end{equation}
Here $\dW^{\leqslant L}$ is spanned by the reds of positions $\leqslant L$. 
Let 
\begin{align}\label{NF}\dN_F\subset\dN\end{align}
 be the subspace of all $x$ preserving $F$. 
The incidence spaces of \S\ref{subsec:unfolded}
are smooth vector bundles
\begin{align}\label{X}
 \scrX_\eta=\{(F,x):F\in\Fl_\eta,\ x\in\dN_F\}
 \cong\dG\times_{P_\eta}\dN_{F_\eta}.
\end{align}
Their dimension is $d_\eta=\dim\Fl_\eta+\dim\dN_{F_\eta}$, and the map
\begin{align}\label{pi}
 \pi_\eta:\scrX_\eta\to\dN,\quad(F,x)\mapsto x
\end{align}
is projective with image 
\begin{align}\label{Y}Y_\eta=\pi_\eta(\scrX_\eta)=\dG\cdot\dN_{F_\eta}.\end{align}
Let $\mathscr L=(\mathscr L_\eta)_\eta $ be the family of constructible complexes
\begin{align}\label{L}\mathscr L_\eta=(\pi_\eta)_*\omega_{\scrX_\eta}=(\pi_\eta)_*\bbC_{\scrX_\eta}[2d_\eta].\end{align}
We call them the \emph{flag-incidence complexes}.
Note that $\mathscr L$ is not an infinite direct sum which does not make sense in $D^b_c(\dN)$.

For types $\eta,$ $\nu$ in the same weight block, 
we write $\eta\leqslant\nu$ if a $\bbZ$-flag of type
$\eta$ contains a $\bbZ$-flag of type $\nu$ termwise, i.e., 
if there are $\bbZ$-flags $F,$ $F'$ of types $\eta,$ $\nu$ such that
\begin{align}\label{PO}\eta\leqslant\nu\iff\dim F'^{\leqslant P}\leqslant\dim F^{\leqslant P}\,,\,\forall P\end{align}
If $\eta\leqslant\nu$ we say that $\eta$ is obtained from $\nu$ by \emph{moving blacks down} within their colours.
We write
$$\eta<\nu\iff \eta\leqslant\nu\,,\,\eta\neq\nu.$$

For a type $\eta$ and a vertex $i$, let $|\eta|_i$ be the sum of the centred positions of the blacks of $\eta$ of vertex $i$. 
If $\eta$ and $\gamma$ are in the same weight block, then  $|\eta|_i-|\gamma|_i\in2d_i\bbZ.$ 
The  \emph{charge} of $\eta$ relative to $\gamma$ is the tuple 
\begin{align}\label{CEG}\mathrm{ch}(\eta,\gamma)=\big((|\eta|_i-|\gamma|_i)/2d_i\big)_{i\in I}\end{align}
It is the degree, in the grading of $\calA^\rho_\mu$ by $\bbZ I$, 
of the elements mapping the $\zeta^{2\gamma}$-weight space to the
$\zeta^{2\eta}$-weight space. 
The set of types obtained from $\gamma$ by \emph{moving blacks up} is
\[
\Up(\gamma)=\big\{\eta:\ |\eta|_i\geqslant|\gamma|_i\ \text{for all }i\in I,\ \text{and}\ |\eta|_i>|\gamma|_i\ \text{for some }i\big\}
\]
A type is called \emph{cyclotomic} if some black lies above all reds. Let 
$$\Cyc(\rho)=\{\text{cyclotomic types}\}.$$

\begin{Remark}\label{rem:Up}\hfill
\begin{enumerate}[label=$\mathrm{(\alph*)}$,leftmargin=8mm]
\item 
We have $\{\eta:\eta<\gamma\}\cap\Up(\gamma)=\emptyset$. 
\item  A type $\eta\neq\gamma$ 
with $|\eta|_i=|\gamma|_i$ for all $i$ has charge $0$,
and is not in $\Up(\gamma)$.
\item $\Up(\gamma)$ is defined by the charge, not by the order on $\bbZ$-flags. The inclusion below may be strict
 $$\{\eta:\gamma<\eta\}\subseteq\Up(\gamma)$$
\end{enumerate}
\end{Remark}

\subsection{Quiver Hecke algebras}\label{subsec:VV25}
The {quiver Hecke algebra} is the Yoneda algebra
\begin{align}\label{KLR}
\bR=\bigoplus_{\eta,\gamma}\Ext^*_{\dG}(\mathscr L_\eta,\mathscr L_\gamma)
\end{align}
where $\mathscr L_\eta$ is defined below.
Multiplication is Yoneda composition, with orthogonal idempotents
$e_\eta=\mathrm{id}_{\mathscr L_\eta}$ and
$e_\nu\bR e_\eta=\Ext^*_{\dG}(\mathscr L_\eta,\mathscr L_\nu)$.
For a finite set $J$ of types, put
\[
 e_J=\sum_{\eta\in J}e_\eta,\quad
 \mathscr L_J=\bigoplus_{\eta\in J}\mathscr L_\eta,\quad
 \bR_J=e_J\bR e_J=\Ext^*_{\dG}(\mathscr L_J,\mathscr L_J).
\]
These are \emph{finite corners}, meaning that $J$, not necessarily
$\bR_J$, is finite. Every finite collection of algebra elements has
a two sided identity $e_J$. All constructions with $\mathscr L$ use such
finite sums.

The algebra $\bR$ is the convolution algebra of the Steinberg varieties
$\scrX_\eta\times_{\dN}\scrX_{\nu}$ by \cite{CG97}. At type
$\eta$, we call the diagonal polynomial sub-algebra
\[
 H^*_{\dG}(\scrX_\eta,\bbC)=H^*_{P_\eta}
 \subset e_\eta\bR e_\eta
\]
the \emph{local Cartan algebra}. Its action is induced by cup product on
$\omega_{\scrX_\eta}$ followed by $\pi_{\eta*}$, see
\cite[\S2.1.2]{VV25} and \cite[\S8.6]{CG97}. It is described as follows. 
Let $\eta$ have $v_c(\eta,P)$ blacks of colour $c$ at centred position $P$, so that a grading of $\dV$ of type $\eta$ has
graded pieces of these dimensions. The stabilizer $P_\eta$ of $F_\eta$ has Levi subgroup 
$$\dG_\eta=\prod_{c,P}GL_{v_c(\eta,P)},$$ whose Weyl group is
\[
\mathfrak S_\eta=\prod_{c,P}\mathfrak S_{v_c(\eta,P)}\ \subset\ \mathfrak S_{\dV}=\prod_c\mathfrak S_{\dim\dV_c}.
\]
Here $\mathfrak S_{\dV}$ is the Weyl group of $\dG$. Equivalently, $\mathfrak S_\eta$ is the stabilizer in $\mathfrak S_{\dV}$ of the point $\zeta^{2\eta}$ of the maximal torus of $\dG$, whose coordinates on $\dV_c$
are the $\zeta^{P-d_i}$ for the blacks of $\eta$ of colour $c$.
Since $\zeta$ is not a root of unity, two coordinates are equal iff the positions are. Then
\[
H^*_{P_\eta}=H^*_{\dG_\eta}=\bbC[\mathfrak t]^{\mathfrak S_\eta}
\]
is the ring of polynomials in the Chern roots of the graded pieces,
symmetric under $\mathfrak S_\eta$. 

Let $\mathfrak m=H^{>0}_{\dG}\subset H^*_{\dG}$ be the augmentation
ideal. For each pair of types define
\[
 e_\nu\bR^\wedge e_\eta
 =\varprojlim_{n\geqslant1}
 \frac{e_\nu\bR e_\eta}
      {\mathfrak m^n(e_\nu\bR e_\eta)}.
\]
Let $\bT$ be the quotient of $\bR$
by the two sided ideal generated by the idempotents $e_\eta$ with $\eta\in\Cyc(\rho)$.
Define $e_\nu\bT^\wedge e_\eta$ similarly. Yoneda composition
extends continuously to these completed bimodules. 
Note that the completion concerns equivariant parameters, not idempotents.
For a module $M$, we call $e_\eta M$ its \emph{$\eta$-piece}.

\begin{Theorem}[\cite{VV25}]\label{thm:V}\label{V}
Fix a finite type datum $(C,D)$ and an integral flavour
parameter $\rho$.
\hfill
\begin{enumerate}[label=$\mathrm{(\alph*)}$,leftmargin=8mm]
\item
There is a topological algebra isomorphism
\begin{align}\label{loc}\calA_{\mu}^{\rho,\wedge}\cong\bR^\wedge.\end{align}
\item
The isomorphism \eqref{loc} identifies the generalized weight projectors on the left with the 
$e_\eta$'s on the right. It preserves the charge
grading.
\item
The isomorphism \eqref{loc} induces an equivalence of the shifted truncated category $\scrO_\mu^\rho$
with the
category of finitely generated smooth locally unital
$\bT^\wedge$-modules with finite dimensional $\eta$-pieces.
The functor identifies $\eta$-pieces with 
generalized weight spaces of $\ell$-weight $Z_{\rho,\eta}$.
\qed
\end{enumerate}
\end{Theorem}

Parts (a) and (b) are \cite[Thm.~4.3.1, Lem.~4.3.5]{VV25}.
The geometric localization and completed Chern-character argument
extends to finite quivers with loops and multiple arrows, see
Lemma~\ref{VQ} for the cyclotomic passage. Outside finite type we
use the Coulomb algebra itself, without a shifted quantum group.

We next define two versions of the completed  B-algebras.
They are identified in Corollary \ref{prop:localB}. 
Recall the map $\ev_\gamma$ in \eqref{EVG}.
Let 
$$\mathfrak m_\gamma=\Ker(\ev_\gamma:R_G\to\bbC).$$

\begin{Definition}\label{def:localB}
Fix a pair of types $(\rho,\gamma)$. Let $A=\calA^\rho_\mu$.
\begin{enumerate}[label=$\mathrm{(\alph*)}$,leftmargin=8mm]
\item
On the Coulomb side put
\[\widehat R_\gamma=\widehat{(R_G)}_{\mathfrak m_\gamma}
,\quad
\widehat{\B(A)}_\gamma=\varprojlim_{n\geqslant1}\B(A)/\mathfrak m_\gamma^n\B(A).\]
The Cartan quotient $q:R_G\twoheadrightarrow\B(A)$ in \eqref{CQ} extends to a continuous surjection
\begin{align}\label{CCQ}
q_\gamma:\widehat R_\gamma\twoheadrightarrow \widehat{\B(A)}_\gamma.\end{align}
The completion on the right hand side is zero if $\ev_\gamma$ does not descend to $\B(A)$. 
Otherwise it is the completed local ring at that character.
\item
On the quiver Hecke side put
\begin{align}\label{BD}
\widehat\B_\gamma=e_\gamma\Delta(\rho,\gamma)
,\quad
 \Delta(\rho,\gamma)=
 \bT^\wedge e_\gamma\Big/
 \overline{\sum_{\eta\in\Up(\gamma)}
     \bT^\wedge e_\eta\bT^\wedge e_\gamma}.
 \end{align}
The sums are finite before taking
closure. Closure is taken separately in each
$e_{\eta}\bT^\wedge e_\gamma$ in the definition of
$\Delta(\rho,\gamma)$.
\end{enumerate}
\end{Definition}

The completed Chern character in Theorem~\ref{thm:V} gives the topological algebra isomorphism
\begin{align}\label{CH}\widehat R_\gamma\cong H^{*,\wedge}_{P_\gamma}.\end{align}
Composing it  with homomorphism to 
$e_\gamma\bT^\wedge e_\gamma$
 and  the quotient gives the map
\begin{align}\label{RG}
r_\gamma:\widehat R_\gamma
   \to H^{*,\wedge}_{P_\gamma}
   \to e_\gamma\bT^\wedge e_\gamma
   \to \widehat\B_\gamma.
\end{align}
The map to $e_\gamma\bT^\wedge e_\gamma$ need not be injective after the cyclotomic quotient. 

\begin{Corollary}
\label{prop:localB} 
For each $(\rho,\gamma)$ the following hold.
\begin{enumerate}[label=$\mathrm{(\alph*)}$,leftmargin=8mm]
\item There is a unique continuous algebra isomorphism
\begin{equation}\label{eq:localB-comparison}
 u_\gamma:\widehat{\B(A)}_\gamma
          \xrightarrow{\cong}\widehat\B_\gamma,
 \quad u_\gamma \circ q_\gamma=r_\gamma.
\end{equation}
\item The map $r_\gamma$
is surjective and $\Ker (r_\gamma)=\Ker (q_\gamma)$.
\end{enumerate}
\end{Corollary}

\begin{proof}
For $n\geqslant1$, put
\[
 \B_{\gamma,n}=\B(A)/\mathfrak m_\gamma^n\B(A),\quad
 \mathcal U_{\gamma,n}
   =(A/AA_{>0})\otimes_{\B(A)}\B_{\gamma,n}.
\]
The vector $v=1\otimes1$ is universal subject to
$A_{>0}v=0$ and $\mathfrak m_\gamma^n v=0$.
The triangular decomposition identifies the degree zero space with
$\B_{\gamma,n}$ and places all degrees in
$-\sum_i\mathbb N\alpha_i$.
Theorem~\ref{thm:V}(a),(b) identifies the generator with $e_\gamma$.
A positive charge element sends it to a type in $\Up(\gamma)$.
Thus $AA_{>0}$ corresponds to the closed submodule
\[
 \overline{\sum_{\eta\in\Up(\gamma)}
       \bR^\wedge e_\eta\bR^\wedge e_\gamma}.
\]
Since $\mathcal U_{\gamma,n}$ is finitely generated, bounded above,
and has finite dimensional generalized weight spaces, it belongs to
$\scrO_\mu^\rho$.
Theorem~\ref{thm:V}(c) therefore shows that its localized action
annihilates the cyclotomic ideal. Replacing $\bR^\wedge$ by
$\bT^\wedge$ adds no relation. The localized module is consequently
\[
 \Delta_{\gamma,n}=\Delta(\rho,\gamma)\big/
       \overline{\Delta(\rho,\gamma)\mathfrak m_\gamma^n},
\]
where the right Cartan action is induced by \eqref{RG}.
Taking the $\gamma$-pieces gives compatible the algebra isomorphism
$\B_{\gamma,n}\cong e_\gamma\Delta_{\gamma,n}$.
Passing to inverse limits  yields
\[
 u_\gamma:
 \widehat{\B(A)}_\gamma
 \xrightarrow{\cong}e_\gamma\Delta(\rho,\gamma)=\widehat\B_\gamma.
\]
The Cartan compatibility gives $u_\gamma q_\gamma=r_\gamma$.
Since $q_\gamma$ is surjective, $r_\gamma$ is surjective, their kernels
coincide, and $u_\gamma$ is unique. Continuity follows from the
inverse limit construction.
\end{proof}

\begin{Proposition}\label{prop:crit}Let $M$ be a locally unital $\bR$-module 
on which  $H^{>0}_{\dG}$ acts locally nilpotently.
Assume $e_\gamma M\neq0$ and $e_\eta M=0$ for each $\eta\in\Up(\gamma)\cup\Cyc(\rho)$. Then:  \hfill
\begin{enumerate}[label=$\mathrm{(\alph*)}$,leftmargin=8mm]
\item
$\zeta^{2\gamma}\in\Spec\B(\calA^\rho_\mu)$. 
\item
$Z_{\rho,\gamma}$ is $\ell$-dominant $\Rightarrow$ Theorem $\ref{Theorem:431}$ holds for $(\rho,\gamma)$.
\end{enumerate}
\end{Proposition}

\begin{proof}
The cyclotomic idempotents act by zero, so $M$ is a $\bT$-module.
Local nilpotence of the augmentation ideal extends each action
$e_\nu\bT e_\eta\times e_\eta M\to e_\nu M$ continuously to
$e_\nu\bT^\wedge e_\eta$: a sufficiently high power of
$\mathfrak m$ kills each vector. Thus closed factorization ideals act by zero whenever
their finite sums do. No finite dimensionality of the total module
is required. For
$0\neq m\in e_\gamma M$, the map 
$$\bT e_\gamma\to M,\quad t\mapsto tm$$ 
kills the subspace $\sum_{\eta\in\Up(\gamma)}\bT e_\eta \bT e_\gamma$.
Thus it factors through
$\Delta(\rho,\gamma)$ and is nonzero. Hence $\widehat\B_\gamma\neq0$, and Corollary \ref{prop:localB} and Lemma \ref{lem:Bcrit}
conclude.
\end{proof}

\section{The cotangent side}\label{sec:cotangent}

The results of this section hold for all Cartan pair $(C,D)$.

\subsection{The conormal correspondence}
Let first us list some properties of the flag-incidence complexes introduced in \eqref{L}.
\begin{enumerate}[label=$\mathrm{(\alph*)}$,leftmargin=8mm]
\item  By proper base change, for $x\in\dN$ we have
\[
H^k(\mathscr L_\eta)_x=H^{k+2d_\eta}\big(\pi_\eta^{-1}(x),\bbC\big),\quad 
H^k\big((i_x)^!\mathscr L_\eta\big)=H^{\BM}_{-k}\big(\pi_\eta^{-1}(x),\bbC\big),
\]
The fibre $\pi_\eta^{-1}(x)$ is the  variety of $\bbZ$-flags of type $\eta$ preserved by $x$.
\item  Since $\pi_\eta$ is proper, we have
$\mathbb D\mathscr L_\eta=\mathscr L_\eta[-2d_\eta]$ and the decomposition theorem \cite{BBD82}
implies that the sheaf $\mathscr L_\eta$ is a direct sum of shifted simple $\dG$-equivariant perverse
sheaves on $\dN$.
\item 
Up to a grading shift, the Ext-algebra is 
$$\bR\cong 
\bigoplus_{\eta,\gamma}H^{\dG,\BM}_*\big(\scrX_\eta\times_{\dN}\scrX_\gamma,\bbC\big).$$ 
The composition of Ext's is identified with the convolution \cite[\S8.6]{CG97}. 
\item $H^*_{\dG}(\scrX_\eta,\bbC)\cong H^*_{P_\eta}$ acts on $\omega_{\scrX_\eta}$ by cup product, 
hence on $\mathscr L_\eta$ through $(\pi_\eta)_*$, yielding
the embedding $H^*_{P_\eta}\subset e_\eta \bR e_\eta$ used in \S\ref{subsec:VV25}.
\end{enumerate}

Next, recall that points of $\dN$ are pairs 
$$x=(A,B),\quad A_c:\dW_c\to\dV_c.$$
Let $\dN^*$ be the dual of $\dN$, with points 
$$y=(A^*,B^*),\quad A^*_c:\dV_c\to\dW_c$$ 
where $B^*_h$ in the
reverse direction of $B_h$. The group $\dG$ acts on 
$T^*\dN=\dN\times\dN^*$ with moment map 
\begin{align}\label{mu}\mu(x,y)(\xi)=\langle\xi\cdot x,y\rangle,\quad \mu(x,y)=\sum_h[B_h,B_h^*]+AA^*\end{align}
in the trace identification $\dg^*\cong\dg$.
Let $\dN_F^\perp$ be the annihilator of the subspace $\dN_F$ in \eqref{NF}. Set
\begin{align}\label{LF}L_F=\dN_F\times\dN_F^\perp
,\quad
\Lambda_F=\{\xi\in L_F\,:\,\mu(\xi)=0\}.\end{align}

\begin{Lemma}\label{lem:NFperp}
For each $y\in\dN^*$ we have
\begin{equation}\label{eq:NFperp}
y\in\dN_F^\perp\iff
B_h^*(F_j^{\leqslant Q})\subset F_i^{\leqslant Q},\quad A^*(F^{\leqslant P})\subset\dW^{<P}.
\end{equation}
\end{Lemma}

\begin{proof}
In a basis adapted to $F$, an elementary map from a basis vector at position $P$ to one at position $Q$ lies in $\dN_F$ iff
$Q<P$ for edges and $Q\leqslant P$ for framings. A covector annihilates all of them iff its matrix coefficients vanish in these
positions, i.e., its own elementary components go from $Q$ to $P$ with $P\leqslant Q$ for edges  $P<Q$ for framings.
\end{proof}

For each type $\eta$ we define
\begin{align}\label{Lambda}
\Lambda_\eta=\{(F,\xi)\,:\,F\in\Fl_\eta,\ \xi\in\Lambda_F\}\subset\Fl_\eta\times T^*\dN
\end{align}
We have the obvious projection
\begin{align}\label{p}p_\eta:\Lambda_\eta\to T^*\dN,\quad(F,\xi)\mapsto\xi\end{align}
whose image is
\begin{align}\label{OE}Z_\eta=p_\eta(\Lambda_\eta)\subset T^*\dN\end{align}
Let $\frakp_\eta=\operatorname{Lie}P_\eta$, where $P_\eta\subset\dG$ is as in \S\ref{subsec:unfolded}.

\begin{Lemma}\label{lem:LZL}
\hfill
\begin{enumerate}[label=$\mathrm{(\alph*)}$,leftmargin=8mm]
\item
$\Lambda_\eta$ is the conormal correspondence of the map $\pi_\eta$.
\item
$Z_\eta$ is a closed conic isotropic subvariety of $T^*\dN$. 
The micro support of $\mathscr L_\eta$ is $\operatorname{SS}(\mathscr L_\eta)\subset Z_\eta$. 
\item
We have
$\Lambda_\eta=\dG\times_{P_\eta}\Lambda_{F_\eta}$ and
$\mu(L_{F_\eta})\subset\frakp_\eta^\perp$.
\end{enumerate}
\end{Lemma}

\begin{proof}\hfill

(a) 
For any $x\in\dN_F$ the image of $d\pi_\eta$ in \eqref{pi} at $(F,x)$ is
$\mathrm{Im}(d\pi_\eta)=\dN_F+\dg\cdot x$.
Hence, for any  covector $y$ we have
$$y(\mathrm{Im}(d_{(F,x)}\pi_\eta))=0\iff (F,x,y)\in\Lambda_\eta.$$

(b) \cite[Ch.~V, Ch.~VIII]{KS90}.

(c) The induced-space description is immediate. Since
$\dN_{F_\eta}$ is $P_\eta$-stable, we have
$$\mu(x,y)(z)=\langle z\cdot x,y\rangle=0
,\quad
\forall z\in\frakp_\eta\,,\, x\in\dN_{F_\eta}\,,\,y\in\dN_{F_\eta}^\perp.$$ 
Hence the moment map takes values in
$\frakp_\eta^\perp$.
\end{proof}

\subsection{Minimal $\bbZ$-flags}
For $\xi=(x,y)$ in $T^*\dN$ let $F_\xi$ be the smallest family of subspaces 
$F_\xi^{\leqslant P}\subset\dV,$ increasing in $P$, with
\begin{equation}\label{eq:Gdef}
A(\dW^{\leqslant L})\subset F_\xi^{\leqslant L},\quad 
B_h(F_{\xi,i}^{\leqslant P})\subset F_{\xi,j}^{<P}
,\quad 
B^*_h(F_{\xi,j}^{\leqslant Q})\subset F_{\xi,i}^{\leqslant Q}.
\end{equation}
There is no condition on $A^*$. 
For a red $u$ at $L$ and a word $w=w_m\cdots w_1$ in the $B_h$ and $B_h^*$, 
we define the \emph{bound} to be the integer
\begin{align}\label{bound}b(w,u)\in\mathbb Z\end{align}
 defined inductively by: 
 \begin{itemize}[leftmargin=8mm]
 \item $b(w,u)=L$ for $m=0$.
\item Adding a letter $B_h$ replaces $b(w,u)$ by the largest position $<b(w,u)$ of the target colour. 
\item Adding a letter $B^*_h$ by the largest position $\leqslant b(w,u)$ of the target colour.
\end{itemize}

\begin{Lemma}\label{lem:words}\hfill
\begin{enumerate}[label=$\mathrm{(\alph*)}$,leftmargin=8mm]
\item
$F_\xi^{\leqslant P}$ is spanned by the vectors $w(\xi)A(u)$ with $b(w,u)\leqslant P$. 
\item
$F_{\xi}^{\leqslant P}$ is spanned by finitely many vectors depending polynomially on $\xi$.
\item
$(F,\xi)\in\Lambda_\eta\Rightarrow F_\xi^{\leqslant P}\subset F^{\leqslant P},\ \forall P.$
\end{enumerate}
\end{Lemma}

\begin{proof}
Let $F_\xi'$ be the family of vector spaces spanned by these vectors.
It is increasing and satisfies \eqref{eq:Gdef}.
Thus, we have $F_\xi\subseteq F_\xi'.$
Conversely, by induction on the length, we have $w(\xi)A(u)\in F_\xi^{\leqslant b(w,u)}$.
Thus we also have
$F_\xi'\subseteq F_\xi.$
This proves part (a). Part (b) is obvious.
For (c), the $\bbZ$-flag 
$F$ satisfies the three conditions \eqref{eq:Gdef}, by \eqref{eq:NF} and \eqref{eq:NFperp}, and $F_\xi$ is the smallest such family.
\end{proof}

We define
\begin{align}\label{LST}
\Lambda^{\rm st}_{F_\gamma}=\{\xi\in\Lambda_{F_\gamma}:\ F_\xi=F_\gamma\}.
\end{align}
We say that the pair $(\rho,\gamma)$ satisfies the condition \eqref{bfE} for the orientation $Q$ if 
\begin{align}\label{bfE}
\Lambda^{\rm st}_{F_\gamma}\neq\emptyset.
\end{align}

\begin{Lemma}\label{lem:noncyc}\hfill
\begin{enumerate}[label=$\mathrm{(\alph*)}$,leftmargin=8mm]
\item
$Z_{\rho,\gamma}$ is $\ell$-dominant $\Rightarrow\gamma$ is not cyclotomic. 
\item
$\Lambda^{\rm st}_{F_\gamma}\neq\emptyset\Rightarrow$ no $\eta\leqslant\gamma$ is
cyclotomic.
\end{enumerate}
\end{Lemma}

\begin{proof}\hfill

(a) Let $(i,a)$ be a black in uncentred notation for which the
centred position $a+d_i$ is maximal among the centred positions of all blacks. 
The exponent of $Z_{i,a+2d_i}$ in $Z_{\rho,\gamma}$ receives $-1$ from the black factor $A_{i,a+d_i}^{-1}$ at $(i,a)$, 
and positive contributions only from a red factor $Z_{i,l+2d_i}$ at $(i,l)$ with $l=a$, or from a black 
factor $A_{j,b+d_j}^{-1}$ at $(j,b)$ with 
$$b+d_jc_{ji}+2d_jn=a+2d_i,\quad
1\leqslant n\leqslant|c_{ji}|$$
In the latter case we have
$$b+d_j\geqslant a+2d_i-d_i|c_{ij}|+d_j>a+d_i,$$ using $d_i|c_{ij}|=d_j|c_{ji}|$ and
$d_j=d_i|c_{ij}|$ when $|c_{ij}|\geqslant2$.
This contradicts maximality. So there is a red at the centred position of $a$.

(b) Choose \(\xi\in\Lambda_{F_\gamma}^{\mathrm{st}}\). 
For a red vector \(u\) at position \(L\), the recursive definition gives \(b(w,u)\leqslant L\leqslant L_{\max}\), 
since adjoining any letter cannot increase the bound. 
Hence Lemma \(\ref{lem:words}\), together with \(F_\xi=F_\gamma\), gives \(F_\gamma^{\leqslant L_{\max}}=\underline V\).
If $\eta\leqslant\gamma$ then a $\bbZ$-flag of type $\eta$ contains a $\bbZ$-flag of type $\gamma$.
Hence
$F^{\leqslant P}=\dV$ for $P\geqslant L_{\max}$, and no black of $\eta$ lies above all reds.
\end{proof}

\begin{Lemma}\label{lem:IU}
Fix $\xi=(x,y)$ in $L_{F_\gamma}$ with $F_\xi=F_\gamma$. We have
\begin{equation}\label{eq:IU}
\frakp_\gamma=\{z\in\dg:\ z\cdot x\in\dN_{F_\gamma},\ z\cdot y\in\dN_{F_\gamma}^\perp\}.
\end{equation}
\end{Lemma}

\begin{proof}
Let $z$ be in the left-hand side, and work over $\bbC[\varepsilon]=\bbC[t]/(t^2)$. 
The $\bbZ$-flag $F_\varepsilon=(1+\varepsilon z)F_\gamma$ of free
$\bbC[\varepsilon]$-modules is preserved by $x$ and satisfies \eqref{eq:NFperp} for $y$, 
since 
$$(1+\varepsilon z)^{-1}\cdot x=x-\varepsilon\,z\cdot x$$
and similarly for $y$. 
The proof of Lemma \ref{lem:words} works over $\bbC[\varepsilon]$, so 
$$F_\gamma^{\leqslant P}\otimes\bbC[\varepsilon]=F_\xi^{\leqslant P}\otimes\bbC[\varepsilon]\subseteq
F_\varepsilon^{\leqslant P}.$$ 
Both are free direct summands of the same rank, hence equal.
Hence, we have 
$$z(F_\gamma^{\leqslant P})\subseteq F_\gamma^{\leqslant P},$$ 
i.e., we have
$z\in\frakp_\gamma$.
\end{proof}

Recall the moment map $\mu$ in \eqref{mu}. 

\begin{Lemma}\label{lem:smooth}
Let $\xi\in\Lambda_{F_\gamma}^{\rm st}$. 
\begin{enumerate}[label=$\mathrm{(\alph*)}$,leftmargin=8mm]
\item
The moment map restricts to a map
$\mu:L_{F_\gamma}\to\frakp_\gamma^\perp$ 
which is submersive at $\xi$.
\item
$\Lambda_\gamma$ is smooth of dimension $\dim\dN$ at $(F_\gamma,\xi)$.
\item
The tangent map of $p_\gamma$ at $(F_\gamma,\xi)$ is injective.
\end{enumerate}
\end{Lemma}

\begin{proof}
The moment map $\mu$ restricts to a map
$\mu:L_{F_\gamma}\to\frakp_\gamma^\perp$ 
by Lemma \ref{lem:LZL}.
Set $\xi=(x,y)$.
By \eqref{mu}, an element $z\in\dg$ annihilating the image of 
$$d\mu_{\xi}:(\delta x,\delta y)\mapsto\mu(\delta x,y)+\mu(x,\delta y)$$ satisfies
\begin{align*}
\langle z\cdot x,\delta y\rangle=
\langle z\cdot\delta x,y\rangle=-\langle\delta x,z\cdot y\rangle=0
,\quad\forall\,\delta x\in\dN_{F_\gamma}\,,\,\delta y\in\dN_{F_\gamma}^\perp,\\
\end{align*} 
i.e., we have
$z\cdot x\in\dN_{F_\gamma}$ and $z\cdot y\in\dN_{F_\gamma}^\perp.$
Thus $z\in\frakp_\gamma$ by Lemma \ref{lem:IU}. 
Hence $d\mu$ is onto
$\frakp_\gamma^\perp=(\dg/\frakp_\gamma)^*$, and $\Lambda_\gamma$ is smooth at $(F_\gamma,\xi)$ of dimension
$$\dim\dg/\frakp_\gamma+\dim\dN_{F_\gamma}+\dim\dN_{F_\gamma}^\perp-\dim\frakp_\gamma^\perp=\dim\dN.$$
For $z\in\dg/\frakp_\gamma$, the tangent map of $p_\gamma$ sends $(z,\delta x,\delta y)$
to $$(z\cdot x+\delta x,z\cdot y+\delta y).$$
If this vanishes then $z\cdot x\in\dN_{F_\gamma}$ and $z\cdot y\in\dN_{F_\gamma}^\perp$, 
hence $z\in\frakp_\gamma$ by Lemma \ref{lem:IU}, and then $$\delta x=\delta y=0.$$
\end{proof}

\subsection{The Lagrangian $Z^{\rm st}_\gamma$}\label{subsec:Zst}
Recall the closed subset $Z_\gamma\subset T^*\dN$ defined in \eqref{OE}. We define
\begin{align}\label{OLZ}
\begin {split}
\Omega_\gamma&=\{\xi\in T^*\dN:\ \dim F_\gamma^{\leqslant P}\leqslant \dim F_\xi^{\leqslant P},\ \forall P\},\\
\Lambda^{\rm st}_\gamma&=\{(F,\xi)\in\Lambda_\gamma:\ F_\xi=F\}=\dG\times_{P_\gamma}\Lambda^{\rm st}_{F_\gamma},\\
Z^{\rm st}_\gamma&=Z_\gamma\cap\Omega_\gamma.
\end{split}
\end{align}
Recall the map $p_\gamma:\Lambda_\gamma\to T^*\dN$ in \eqref{p}.

\begin{Lemma}\label{lem:Z}\hfill
\begin{enumerate}[label=$\mathrm{(\alph*)}$,leftmargin=8mm,itemsep=1mm]
\item $\Omega_\gamma$ is open, conic and $\dG$-stable in $T^*\dN$.
\item $Z^{\rm st}_\gamma=p_\gamma(\Lambda^{\rm st}_\gamma)=\{\xi\in Z_\gamma\,:\,F_\xi\text{ has type }\gamma\}$.
It is closed in $\Omega_\gamma$.
\item 
$\Lambda^{\rm st}_\gamma$ is smooth of dimension $\dim\dN$.
\item
$Z^{\rm st}_\gamma$ is smooth, closed, conic and Lagrangian in $\Omega_\gamma$.
\item
$p_\gamma$ restricts to
a $\dG$-equivariant isomorphism $\Lambda^{\rm st}_\gamma\to Z^{\rm st}_\gamma$.
\item $Z^{\rm st}_\gamma\cap Z_\eta\neq\emptyset\Rightarrow\eta\leqslant\gamma$. 
\item
$\eta\in\Up(\gamma)$ or $\eta\in\Cyc(\rho)\Rightarrow Z^{\rm st}_\gamma\cap Z_\eta=\emptyset$. 
\end{enumerate}
\end{Lemma}

\begin{proof}\hfill

(a) By Lemma \ref{lem:words}, the condition $\dim F_\gamma^{\leqslant P}\leqslant \dim F_\xi^{\leqslant P}$ is a nonvanishing of minors. 
Rescaling $y$ does not change $F_\xi$, and $gF_\xi=F_{g\xi}$. 

(b) If $(F,\xi)\in\Lambda_\gamma$ then $F_\xi\subseteq F$ by Lemma \ref{lem:words}, hence $\dim F_\xi^{\leqslant P}\leqslant d_{\gamma,c,P}$, with
equality for all $P$ iff $F_\xi=F$. 
Therefore, on $Z_\gamma$ the condition defining $\Omega_\gamma$ says that $F_\xi$ has type $\gamma$. 
The set $Z_\gamma$ is closed,
being the image of a proper map. 

(c), (d), (e) If $(F,\xi)\in\Lambda_\gamma$ with $\xi\in Z^{\rm st}_\gamma$, then $F_\xi\subseteq F$ with equal types, so $F=F_\xi$.
Thus $p$ is injective on $\Lambda^{\rm st}_\gamma=p_\gamma^{-1}(Z^{\rm st}_\gamma)$, 
and it is proper as a base change of the proper map $p$. 
A proper continuous bijection of locally compact spaces is a homeomorphism. 
After translation, we have $F=F_\gamma$. Lemma \ref{lem:smooth} gives smoothness,
the dimension and injectivity of the tangent map. An injective immersion which is a homeomorphism onto a closed subset is a
closed embedding. Further $Z^{\rm st}_\gamma$ is conic by (a) and isotropic, hence Lagrangian. 

(f), (g) If $(F,\xi)\in\Lambda_\eta$, then $F\supseteq F_\xi$, of
type $\gamma$, thus we have  $\eta\leqslant\gamma$.
Such $\eta$ is not in $\Up(\gamma)$ by Remark \ref{rem:Up}(c).
Since $Z^{\rm st}_\gamma\neq\emptyset$ forces
$\Lambda^{\rm st}_{F_\gamma}\neq\emptyset$, it is not cyclotomic by Lemma \ref{lem:noncyc}.
\end{proof}

We first record an elementary fact for the preprojective
algebra $\Pi(H)$ of a finite loopless graph $H$. For a representation $B$
on a vertex-graded vector space $E$ and a vertex $i$, let $K_i(B)$
be the kernel of the maximal quotient of $B$ supported at $i$.
It is the subrepresentation given by
\begin{align}\label{K1}
(K_i(B))_j=E_j
,\quad
(K_i(B))_i=\sum_{h:t(h)=i}\operatorname{Im}B_h
,\quad \forall j\neq i.
\end{align}
where the sum runs over arrows of the doubled graph.

\begin{Lemma}\label{lem:prepro}
Fix a finite descending flag $E^\bullet$ of $H$-vertex-graded subspaces
\[
 E=E^{(0)}\supseteq\cdots\supseteq E^{(N)}=0,
\]
whose $k$-th quotient $E^{(k-1)}/E^{(k)}$ is supported at $i_k$, allowing zero quotients. 
The set
\[
 \mathcal T(E^\bullet)=
 \{B:K_{i_k}(B|_{E^{(k-1)}})=E^{(k)}\text{ for all }k\}
\]
is  locally closed and is either empty or irreducible. 
\end{Lemma}

\begin{proof}
We argue by induction on $N$.
Put $i=i_1$ and split $E_i=U\oplus T$, where $U=E_i^{(1)}$.
The restriction of the representations $B$ to $E^{(1)}$ satisfying the remaining
conditions form an irreducible variety if nonempty by induction. Restrict to its
open subset where the signed incoming map
\[
 a:\bigoplus_{h:t(h)=i}E_{s(h)}^{(1)}\longrightarrow U
\]
is onto. All incoming maps at $i$ are then fixed, and the only new
variables are the outgoing maps on $T$, assembled into a map $b$.
The additional preprojective equation is $a b=0$.
The equations at the other vertices are unchanged. Thus the extensions
form the vector bundle $\Hom(T,\ker(a))$ over this open
subset. 
\end{proof}

\begin{Lemma}\label{lem:stable-irreducible}
Let $(C,D)$ be of finite type. When nonempty, the sets 
$\Lambda_{F_\gamma}^{\rm st}$, $\Lambda_\gamma^{\rm st}$ and
$Z_\gamma^{\rm st}$ are irreducible, hence connected.
\end{Lemma}

\begin{proof}
Let $H$ be the underlying graph of $\dQ$.
For a colour \(c=(i,r)\) and an integer \(L\equiv r\mod{2d_i}\), 
let \(W_c(L)\) be the span of the reds of colour \(c\) at centred position \(L\).
We abbreviate $p_c=2d_i$.
To the graph $H$  we attach a vertex $r_{c,L}$ for each
nonzero $W_c(L)$, carrying that space, with its edge oriented towards
$\dV_c$ and maps $A_{c,L}:W_c(L)\to\dV_c$,  $A^*_{c,L}:\dV_c\to W_c(L)$. 
Call the extended graph $\widehat H$. 
We consider the vector space 
$$E=\dV\oplus\bigoplus_{c,L}W_c(L).$$ 
It is graded by the vertices
$c$ and $r_{c,L}$ of $\widehat H$. Let $V_c(P)$ be the span of the
blacks of colour $c$ at $P$. We have
$$F_{\gamma,c}^{\leqslant P}=\bigoplus_{Q\leqslant P}V_c(Q).$$

Choose a finite interval containing all black and red positions,
extended at both ends by more than $\max_c p_c$. 
We define a  sequence \(\mathbf i=(i_1,\ldots,i_N)\) of vertices of $\widehat H$ as follows.
Consider all pairs \((c,P)\), where \(c\) is a colour and \(P\) runs 
through its allowed centred positions in the chosen finite interval, including positions with zero multiplicity, i.e.,
for each colour \(c=(i,r)\) we include all integers within the chosen finite interval such that
\[
P\equiv r\mod{p_c},\quad p_c=2d_i,
\]
even if \(\gamma\) has no black at \(P\).
Add all framing vertices \(r_{c,L}\).
Assign heights
\[
h(c,P)=P,\quad h(r_{c,L})=L+\tfrac12.
\]
Order all pairs \((c,P)\) and all framing vertices \(r_{c,L}\) by decreasing height. 
This gives a sequence $$\mathbf i=(i_1,\ldots,i_N)$$ of vertices of the graph \(\widehat H\), 
where each occurrence \((c,P)\) contributes the vertex \(i_k=c\), 
and each occurrence \(r_{c,L}\) contributes \(i_k=r_{c,L}\). 
The same colour vertex may occur several times in \(\mathbf i\), corresponding to different centred positions.
Note that the sequence \(\mathbf i\) is determined by the pair \((\rho,\gamma)\) and the interval choice.

We construct a flag $E^\bullet$ as follows.
Start with $E^{(0)}=E$. At $(c,P)$ we replace the $c$-component
$F_{\gamma,c}^{\leqslant P}$ by $F_{\gamma,c}^{<P}$, i.e., we remove
$V_c(P)$. At $r_{c,L}$ we replace the component $W_c(L)$ by $0$.
Leave all other components unchanged. This defines a finite descending
flag $E^\bullet=(E^{(k)})_{k=0}^N$ of $\widehat H$-vertex-graded
subspaces. Note that $E^\bullet$ is a flag, not a $\mathbb Z$-flag. 
Its index $k$ records the ordered removal steps and it includes the
framing vertices $r_{c,L}$, whereas
the $\bbZ$-flag $F_\gamma$ is indexed by centred positions $P$ and
lives on $\dV$ alone. 

Equal height occurrences belong to non adjacent vertices, so their order is not important.
Since adjacent colours have no common allowed position by \eqref{parity},
both directions of a colour edge strictly lower position.
The representations of $\Pi(\widehat H)$ preserving $E^\bullet$ therefore satisfy 
\eqref{eq:NF} and \eqref{eq:NFperp}. Moreover, we have
$$A_{c,L}(W_c(L))\subseteq F_c^{\leqslant L}
,\quad
A^*_{c,L}(F_c^{\leqslant L})=0,$$ so the relation
$A^*_{c,L}\circ A_{c,L}=0$ holds. The remaining preprojective
relations are  $\mu=0$.

We claim that, as reduced locally closed varieties, we have
\begin{align}\label{LT}\Lambda_{F_\gamma}^{\rm st}=\mathcal T(E^\bullet).\end{align}
Then the lemma follows from Lemma \ref{lem:prepro}.

To prove the claim,
fix an arbitrary \(\xi\in\Lambda_{F_\gamma}\). 
By Lemma~\(\ref{lem:words}\), we have \(F_\xi\subseteq F_\gamma\). We show that 
\begin{align}\label{IFF}F_\xi=F_\gamma\iff \text{ all removals defining \(E^\bullet\) are maximal}.\end{align}
Let $T_h$  denote a colour component of \(B_h\) or \(B_h^*\), according to the direction of the doubled colour arrow \(h\).
By definition of the bound, we have
\[
b(T_hw,u)
=
\max\{Q\,:\,Q\text{ is an allowed position of }c\,,\, Q<b(w,u)\}.
\]
Adjacent colours have disjoint sets of allowed positions.
So, for \(P\) an allowed position of \(c\), we have
\begin{align}\label{BTB}
b(T_hw,u)<P
\quad\Longleftrightarrow\quad
b(w,u)<P.
\end{align}
Now, by Lemma~\ref{lem:words}, the vector space $F_{\xi,c}^{<P}$ is spanned by the framed words ending at colour \(c\) 
whose bound is \(<P\).
Such a word has one of two forms:
\begin{itemize}[leftmargin=8mm]
\item It is \(A_c(u)\), where the red vector \(u\) has position \(L<P\). These words span \(A_c(W_c^{<P})\).
\item It is \(T_hwA(u)\), where \(h:d\to c\) is the last arrow. 
By \eqref{BTB}, its bound  is \(<P\) if and only if \(wA(u)\) has also bound \(<P\). 
These words span \(T_h(F_{\xi,d}^{<P})\).
\end{itemize}
Hence, we have
\begin{align}\label{FAT}
F_{\xi,c}^{<P}=A_c(W_c^{<P})+
       \sum_{h:d\to c}T_h(F_{\xi,d}^{<P}).
\end{align}

The standard flag has
$F_{\gamma,c}^{<P+p_c}=F_{\gamma,c}^{\leqslant P}$.
Adjacent colours have disjoint allowed positions, so arrows strictly lower.
If $F_\xi=F_\gamma$, then \eqref{FAT} makes every colour removal maximal. The framing steps are maximal too, since
$A^*_{c,L}(F_{\gamma,c}^{\leqslant L})=0$.

Conversely, maximality at $(c,P+p_c)$ gives
\[
 F_{\gamma,c}^{\leqslant P}=A_c(W_c^{\leqslant P})+
       \sum_{h:d\to c}T_h(F_{\gamma,d}^{<P+p_c}).
\]
If $F_\xi\ne F_\gamma$, choose the largest $P$ where some $c$ satisfies
$$F_{\xi,c}^{\leqslant P}+F_{\gamma,c}^{<P}\ne F_{\gamma,c}^{\leqslant P}.$$
By the maximality of $P$, for each neighborhood colour $d$ of $c$, we have
$$F_{\xi,d}^{<P+p_c}+F_{\gamma,d}^{\leqslant P}=F_{\gamma,d}^{<P+p_c}.$$
By the word bounds, the map $T_h$ sends the first summand into $F_{\xi,c}^{\leqslant P}$.
By the strict lowering, it sends the second into $F_{\gamma,c}^{<P}$.
Also, we have  $A_c(W_c^{\leqslant P})\subseteq F_{\xi,c}^{\leqslant P}$.
This contradicts the choice of $P$.
Hence $F_\xi=F_\gamma$.

We have proved \eqref{IFF}. The equality \eqref{LT} follows.

The first paragraph then proves
irreducibility of $\Lambda_{F_\gamma}^{\rm st}$. Finally, the variety
$\dG\times\Lambda_{F_\gamma}^{\rm st}$ surjects onto
$\Lambda_\gamma^{\rm st}\cong Z_\gamma^{\rm st}$ by
Lemma~\ref{lem:Z}. These varieties
are irreducible images of irreducible varieties.
\end{proof}

\section{The contraction}\label{sec:sym}

In this section  $(C,D)$ is of type ADE with $D=1$, unless explicitly stated otherwise.

\subsection{Homogeneous points}\label{subsec:graded}
Given a pair of types $(\rho,\gamma)$, we
use centred gradings $V_i(P)$ and $W_i(L)$, with dimension vectors
$v_\gamma$ and $w_\rho$. For $c=(i,r)$ we have
$$\dV_c=\bigoplus_{P\equiv r\, (2)}V_i(P)
,\quad
F_{\gamma,c}^{\leqslant P}=\bigoplus_{P'\leqslant P,\ P'\equiv r\, (2)}V_i(P').$$
We define
$$\dG_\gamma=\prod_{i,P}GL(V_i(P)).$$
The stabilizer $P_\gamma$ has Levi
$\dG_\gamma$ and unipotent radical $U_\gamma$.
At a vertex, a vector in $V_i(P)$ or $W_i(P)$ has \emph{centred degree} $P$.
For a \emph{matrix unit} $T:E(P)\to F(Q)$ of an arrow or framing map, where $E,F$ is some $V_i$ or $W_i$,
its \emph{centred degree}  is $P-Q$. It is the difference of the centred degrees of the source and target. 
Adjacent colours have positions of opposite parities by \eqref{parity}. By \eqref{eq:NF} and \eqref{eq:NFperp}, the matrix
units allowed in $L_{F_\gamma}$ and their \emph{minimal centred degrees} are defined as follows
\begin{equation}\label{tableau}
\begin{tabular}{llll}
component &  centred degrees & minimal centred degree $\delta$ & homogeneous part\\\hline
$B_h$, $h\in Q$ & $1,3,5,\dots$ & $1$ & $V_i(P)\to V_j(P-1)$\\
$B^*_h$, $h\in Q$ & $1,3,5,\dots$ & $1$ & $V_j(P)\to V_i(P-1)$\\
$A$ & $0,2,4,\dots$ & $0$ & $W_i(L)\to V_i(L)$\\
$A^*$ & $2,4,6,\dots$ & $2$ & $V_i(P)\to W_i(P-2)$
\end{tabular}
\end{equation}

We define the \emph{excess} of a matrix unit to be its centred degree minus $\delta$ (=the minimal one). 
It is a nonnegative even integer. A point 
$\xi=(x,y)$ in $L_{F_\gamma}$ 
is said \emph{homogeneous} if all its matrix units have excess $0$.
Let 
\begin{align}\label{h} 
h:L_{F_\gamma}\to L_{F_\gamma},\quad 
\xi=(x,y)\mapsto\xi_0=(x_0,y_0)
\end{align}
be the map taking a pair to its homogeneous part, which 
is obtained by deleting the matrix units of positive excess.

For $\xi=(x,y)$ an homogeneous point in $L_{F_\gamma}$, we have
\[
x=(A_i,B_h)_{i\in I,h\in Q},\quad y=(A_i^*,\varepsilon(\bar h)B_{\bar h})_{i\in I,h\in Q}.
\]
For all arrows $h$ of the double quiver $\bar Q$ we have
$$B_h:V_i(P)\to V_j(P-1)
,\quad
A_i:W_i(L)\to V_i(L)
,\quad
A_i^*:V_i(P)\to W_i(P-2)$$
Let $\varepsilon(h)=\pm1$ according as $h\in Q$ or not. We have 
$$\mu(\xi)=\mu^\bullet(\xi)=\sum_{h\in\bar Q}\varepsilon(h)B_{\bar h}B_h+AA^*:V_i(P)\to V_i(P-2).$$ 
We define
\begin{align}\label{LB}
\begin{split}
\Lambda_{F_\gamma}^\bullet&=\{\xi\in\Lambda_{F_\gamma}\,:\, \xi\ \text{homogeneous}\},\\
\Lambda_{F_\gamma}^{\rm st,\bullet}&
=\{\xi\in\Lambda_{F_\gamma}^\bullet\,:\,\dV\ \text{is spanned by all}\ w(\xi)A(u)\,,\,u\in\dW\,,\, w\ \text{a word in }x,y\},\\
 \M^\bullet_\gamma&=\Lambda_{F_\gamma}^{\rm st,\bullet}/\dG_\gamma.
\end{split}
\end{align}
For $D=1$, the homogeneous maps, moment equation and gauge action
are those of Nakajima's graded variety \cite[\S4]{Nak01}.
Note that the graded varieties for $Q$ and $\dQ$ are identified, not their ungraded varieties.

\begin{Lemma}\label{lem:free}\hfill
\begin{enumerate}[label=$\mathrm{(\alph*)}$,leftmargin=8mm,itemsep=1mm]
\item If $\dV$ is spanned by the vectors $w(\xi)A(u)$, then the stabilizers of $\xi$ in $\dG$ and $\dg$ are
trivial and $P_\gamma$ acts freely on $\Lambda^{\rm st}_{F_\gamma}$.
\item $\dG_\gamma$ acts freely on $\Lambda_{F_\gamma}^{\rm st,\bullet}$.
\item 
$\M^\bullet_\gamma=\M_{\dQ}^\bullet(v_\gamma,w_\rho)$
is smooth, quasi-projective. The quotient map
$\Lambda_{F_\gamma}^{\rm st,\bullet}\to\M^\bullet_\gamma$ is a $\dG_\gamma$-torsor.
\end{enumerate}
\end{Lemma}

\begin{proof}\hfill

(a) If $g\in\dG$ fixes $\xi$, then $gA=A$ and $g$ commutes with all $B_h,B_h^*$, so $g$ fixes all $w(\xi)A(u)$, which span $\dV$.
Thus, we have $g=1$. 
The proof for $\xi\in\dg$ is the same. 
Points of $\Lambda^{\rm st}_{F_\gamma}$ satisfy the hypothesis by Lemma \ref{lem:words}, 
since $F_\xi=F_\gamma$ is
all of $\dV$ in high degree. 

(b) The freeness is proved as in (a). 

(c) This is the definition and the basic properties of graded quiver varieties \cite[\S\S3-4]{Nak01}.
\end{proof}

Let $\M_{\dQ}(v,w)=\mu^{-1}(0)^{\rm st}/\dG$ be the quiver variety
of the unfolded quiver $\dQ$ with spaces $\dV,$ $\dW$ of dimensions $\v$, $\w$.
Here ${\rm st}$ means generation by the framing. 
We define
\begin{align}\label{LRG2}\mathfrak L_{\rho,\gamma}=\Lambda_{F_\gamma}^{\rm st}/P_\gamma.\end{align}
By Lemma \ref{lem:Z} we have a $\dG$-torsor, see also \eqref{LRG1},
\begin{align}\label{QEG}
q:Z_\gamma^{\rm st}\to\mathfrak L_{\rho,\gamma}
,\quad
\mathfrak L_{\rho,\gamma}\cong Z_\gamma^{\rm st}/\dG
\end{align}
We consider the one parameter subgroup 
\begin{align}\label{TR}
T_\rho(t)=(\rho(t^{-1}),t^{-1}),\quad t\in\bbC^\times.
\end{align}
Its action on $\M_{\dQ}(v,w)$ is
\[
T_\rho(t)(B,A,A^*)=
\bigl(t^{-1}B,\ A\rho_W(t),\ t^{-2}\rho_W(t)^{-1}A^*\bigr),
\quad \rho_W(t)|_{W_c(L)}=t^L.
\]
Here $B$ denotes all doubled arrow maps. 
The subgroup and its fixed points are unchanged by this reparametrization. 
The equivariant parameter $a$  is the restriction
of the ambient rotation class $\hbar$.
It is the negative of the
standard degree two generator for the coordinate $t$ above.

\begin{Lemma}\label{lem:attracting-stratum}
 The action of $T_\rho$ on $\M_{\dQ}(v,w)$
has $\M_{\dQ}^\bullet(v_\gamma,w_\rho)$ as a fixed locus. Forgetting the flag gives
an isomorphism
\[
 \mathfrak L_{\rho,\gamma}\xrightarrow{\cong}\mathfrak M_\gamma^+,
 \quad
 \mathfrak M_\gamma^+=\{z\in \M_{\dQ}(v,w):\lim_{t\to0}T_\rho(t)z
                     \text{ exists in }\M_{\dQ}^\bullet(v_\gamma,w_\rho)\}.
\]
The limit map $\mathfrak M_\gamma^+\to\M_{\dQ}^\bullet(v_\gamma,w_\rho)$ is given by the map $h$ in \eqref{h}.
It  is a Zariski locally trivial fibration in affine spaces. 
\end{Lemma}

\begin{proof}
The action preserves generation and multiplies $\mu$ by $t^{-2}$.
Put $\lambda_\gamma(t)|_{V_c(P)}=t^{-P}$. The compensated action
$\lambda_\gamma(t)T_\rho(t)$ has weights
\[
 P-Q-1,\quad L-P,\quad P-L-2
\]
on the maps
$$B:V_c(P)\to V_d(Q)
,\quad
A:W_c(L)\to V_c(P)
,\quad
A^*:V_c(P)\to W_c(L).$$
By parity, nonnegative
weights are  the flag conditions \eqref{eq:NF}, \eqref{eq:NFperp}.
The limit is $h(\xi)$. 
By Lemma~\ref{lem:words}, a framed word starting at position $L$
and of length $m$ has bound $P=L-m$. Its image in
$F_{\gamma,c}^{\leqslant P}/F_{\gamma,c}^{<P}=V_c(P)$, where $c$
is its final colour, is the corresponding word in the homogeneous
maps $h(\xi)=(B_0,A_0,B_0^*,A_0^*)$.
Since $F_\xi\subseteq F_\gamma$, induction on $P$ gives
$F_\xi=F_\gamma$ if and only if these homogeneous framed words span
every $V_c(P)$. Equivalently, the smallest subrepresentation of
the doubled colour quiver containing $A_0(\dW)$ is all of $\dV$.
The outgoing framing map $A_0^*$ plays no role in this condition.
The moment equation passes to the limit, hence forgetting the flag maps
$\mathfrak L_{\rho,\gamma}$ into $\mathfrak M_\gamma^+$, with the
asserted limit map.

Conversely, near a graded fixed point choose homogeneous framed
words forming bases of the colour spaces. Their nonvanishing
determinants define a $T_\rho$-invariant open chart: a word
of bound $P$ scales by $t^P$. In these word bases the action is
 the compensated action above. Each point attracted to
this chart lies in it. The existence of its stable limit therefore
forces the flag conditions and $F_\xi=F_\gamma$. Two such
representatives differ by an element of $P_\gamma$, since their
minimal flags coincide. The word basis charts give regular local
inverses, proving the isomorphism. 
 
Finally, the variety $\M_{\dQ}(v,w)$ is smooth and
quasi-projective, so the Bia\l ynicki-Birula theorem
\cite[Thm.~4.1]{BB73} gives the affine space fibration.
\end{proof}

\subsection{The local and global Hikita conjectures}\label{subsec:global-conjecture}

The tautological bundle on $\M_{\dQ}^\bullet(v_\gamma,w_\rho)$ is
 $$\calV_i(P)=\Lambda_{F_\gamma}^{\rm st,\bullet}\times_{\dG_\gamma}V_i(P).$$
 The local Kirwan map is the map
\begin{align}\label{LKM1}
\kappa_\gamma:\widehat R_\gamma\xrightarrow{\cong} \widehat H^*_{\dG_\gamma}\to H^*(\M_{\dQ}^\bullet(v_\gamma,w_\rho),\bbC)
\end{align}
sending a $\mathfrak S_\gamma$-symmetric polynomial in the Chern roots of $V_i(P)$ to the same polynomial in the
Chern roots of the $\calV_i(P)$. 
The first isomorphism is \eqref{CH}.
It extends to the completion because $H^{>0}(\M_{\dQ}^\bullet(v_\gamma,w_\rho),\bbC)$ is nilpotent.
Recall the map  in \eqref{RG}
$$r_\gamma:\widehat R_\gamma\to \widehat\B_\gamma$$
The {local Hikita conjecture} for $(\rho,\gamma)$ is the statement
\begin{align}\label{locHik}
\Ker(r_\gamma)=\Ker(\kappa_\gamma),
\end{align}
i.e., the map $\kappa_\gamma$ induces an isomorphism of rings
\begin {align}\label{PHG1}
\Phi_\gamma:\ \widehat\B_\gamma\ \xrightarrow{\cong}\ H^*\big(\M_{\dQ}^\bullet(v_\gamma,w_\rho),\bbC\big).
\end{align}
The inclusion $\subseteq$ says that $\Phi_\gamma$ is a well defined ring homomorphism. 
It is the \emph{local comparison map}.

Let $\M_Q(v,w)$ be the ungraded Nakajima quiver variety of $Q$
for the total spaces
$$V_i=\bigoplus_PV_i(P),\quad W_i=\bigoplus_LW_i(L)
,\quad v=(\dim V_i)
,\quad
w=(\dim W_i).$$
The graded varieties occur as fixed loci in $\M_Q(v,w)$, as recalled in
\S\ref{subsec:ungraded}.
Let $F=T_W$ be the torus of $W$ and $\bbC^\times$ the loop rotation. 
Set 
\begin{align}\label{AFH}
A_{F,\hbar}=\text{quantized
homological Coulomb branch of the same gauge theory with flavour } F.
\end{align}
Let $G=G_V$.
The Cartan subalgebra of $A_{F,\hbar}$ and the global Kirwan map are
\begin{align}\label{GKM}
R_{F,\hbar}=H^*_{G\times F\times\bbC^\times}
,\quad
\kappa_{F,\hbar}:R_{F,\hbar}\to H^*_{F\times\bbC^\times}(\M_Q(v,w),\bbC).\end{align}
Let $s_\rho$ be given by the red positions and $\hbar=1$.
The specialization of $R_{F,\hbar}$ at  $s_\rho$ is $H^*_G$. 
Its completion at $\gamma$ is $\widehat R_\gamma$. 
The
connected components of the affine Grassmannian of $G$ are indexed by $\pi_1(G)$. 
This
gives the {topological charge grading}
\[
 A_{F,\hbar}=\bigoplus_{\mathbf m}A^{\rm ch}_{\mathbf m},\quad
 A^{\rm ch}_{\mathbf m}A^{\rm ch}_{\mathbf n}
 \subseteq A^{\rm ch}_{\mathbf m+\mathbf n},
\]
where $\mathbf m=(m_i)\in\bbZ^I$, with $m_i=0$ if $V_i=0$.
Let 
$$\pi_1(G)^\vee=\Hom(\pi_1(G),\bbC^\times).$$
Choose a generic integral cocharacter $\theta:\bbC^\times\to \pi_1(G)^\vee$. 
The sign of $\theta$ is chosen to match the highest weight convention for $\B$ in \S\ref{subsec:algebras}.
Pairing charges with $\theta=(\theta_i)$ gives a $\bbZ$-grading. 
For $A=A_{F,\hbar}$ in \eqref{AFH}, or for a central specialization of this algebra, write 
\begin{align}\label{[k]}A_{[k]}=\{a\in A:\ \theta(t)\cdot a=t^k a\text{ for all }t\in\bbC^\times\}
 =\bigoplus_{\theta\cdot\mathbf m=k}A^{\rm ch}_{\mathbf m}
\end{align}
where $\theta\cdot\mathbf m=\sum_i\theta_i m_i$.
The brackets $[k]$ always denote this charge degree, not the
cohomological degree. The Cartan algebra, the flavour
and the quantization parameters have charge zero. Thus
$A_{[0]}$ is a degree zero subalgebra. Set
\begin{align}\label{BTFH}
 \B^\theta_{F,\hbar}=\B^\theta(A_{F,\hbar}),\quad
 \B^\theta(A)=A_{[0]}\Big/\sum_{k>0}A_{[-k]}A_{[k]}.
\end{align}
The global Hikita conjecture is the following statement: 
there is a graded $R_{F,\hbar}$-algebra isomorphism compatible with the Kirwan map
\begin{align}\label{GHC}\B^\theta_{F,\hbar}\cong H^*_{F\times\bbC^\times}(\M_Q(v,w),\bbC).\end{align}
We distinguish \eqref{GHC} from the
localization theorem used to pass to its local form.
We abbreviate 
$$A^\rho=A_{s_\rho}=A_{F,\hbar}\otimes_S\bbC_{s_\rho}.$$
The integral completions of
$\calA^\rho_\mu$ and $A^\rho$ are isomorphic, preserving
multiplication, the weight idempotents, the charge grading, and the
completed Cartan maps. Integral completion means the
componentwise completion at the left and right integral Cartan
characters, as in \S\ref{subsec:VV25}. It is not the completion of $\hbar$ at zero. 
There is an isomorphism compatible with the maps from $\widehat R_\gamma$
\begin{align}\label{H2}\widehat{\B^\theta(A^\rho)}_\gamma
 \cong\widehat{\B(\calA^\rho_\mu)}_\gamma
 =\widehat\B_\gamma,
\end{align}

Indeed, K-theoretic localization and the completed Chern character
\cite[\S2.1.3, Thm.~4.3.1, Lem.~4.3.5]{VV25}, and homological
localization \cite[Thm.~4.4, \S4.3]{Web24}, identify both completions
with the same convolution algebra of fixed flag correspondences.
These algebra isomorphisms preserve the charge, so they identify the
positive charge relations as well.

\begin{Proposition}\label{prop:global-local}
If the global Hikita conjecture \eqref{GHC} holds, then, the local Hikita conjecture \eqref{PHG1} holds for each $(\rho,\gamma)$.
\end{Proposition}

\begin{proof}
By the localization result \eqref{H2} and right exactness of the B-algebra construction, 
the algebra $\widehat\B_\gamma$ is the completion at $\gamma$ of
$\B^\theta_{F,\hbar}\otimes\bbC_{s_\rho}$, which by \eqref{GHC} is 
$H^*_{F\times\bbC^\times}(\M_Q(v,w),\bbC)\otimes\bbC_{s_\rho}$. Since $\M_Q(v,w)$ is equivariantly formal,
the localization theorem identifies the latter with 
$$H^*(\M_Q(v,w)^{s_\rho},\bbC)=\bigoplus_{v_\gamma}H^*(\M^\bullet(v_\gamma,w_\rho),\mathbb C)$$ 
by Lemma \ref{lem:fixed}, and
completing at $\gamma$ keeps the summand $v=v_\gamma$, see the proof of Lemma \ref{lem:kirwan}. 
All maps are compatible
with the maps from $\widehat R_\gamma$, which on the right is $\kappa_\gamma$. So we have
$$\Ker(r_\gamma)=\Ker(\kappa_\gamma).$$
\end{proof}

\subsection{Retraction onto the graded quiver variety}\label{subsec:retraction}
Recall the definition of the bound in \eqref{bound}.

\begin{Lemma}\label{lem:leading}
Fix $\xi$, $\xi_0$ in $L_{F_\gamma}$ as in \eqref{h}.
Let $u$ be a red at $L$ and $w$ a word of length $m$. We have 
\begin{enumerate}[label=$\mathrm{(\alph*)}$,leftmargin=8mm,itemsep=1mm]
\item 
$b(w,u)=L-m$,
\item
$w(\xi_0)A_0(u)\in V(L-m)$,
\item
$w(\xi)A(u)\equiv w(\xi_0)A_0(u)\mod{F_\gamma^{<L-m}}.$
\end{enumerate}
\end{Lemma}

\begin{proof}
Since adjacent colours have positions of opposite parities, each letter lowers the bound by  $1$. 
The congruence is proved by an induction on $m$. 
For $m=0$ the units of $A-A_0$ have centred degree $\geqslant2$. 
If 
$$v=v_0+v',\quad
v_0\in V(b),\quad
v'\in F_\gamma^{\leqslant b-2},$$ and $C$ is a
letter with homogeneous part $C_0$, then we have
$$Cv=C_0v_0+(C-C_0)v_0+Cv',\quad C_0v_0\in V(b-1),\quad
(C-C_0)v_0+Cv'\,,\, C_0v_0\in F_\gamma^{\leqslant b-3}=F_\gamma^{<b-1}$$ 
since $C-C_0$ has centred degree $\geqslant3$ and $C$ has centred degree $\geqslant1$.
\end{proof}

\begin{Lemma}\label{lem:stab}\hfill
\begin{enumerate}[label=$\mathrm{(\alph*)}$,leftmargin=8mm,itemsep=1mm]
\item 
For $\xi$, $\xi_0$ in $L_{F_\gamma}$ as in \eqref{h},
the following are equivalent: 
\begin{enumerate}[label=$\mathrm{(\arabic*)}$, leftmargin=8mm,itemsep=1mm]
\item 
$F_\xi=F_\gamma$.
\item
 $F_{\xi_0}=F_\gamma$.
\item
For all $i,P$, the vector space $V_i(P)$ is spanned by the $w(\xi_0)A_0(u)$'s with $b(w,u)=P$. 
\end{enumerate}
\item
$\Lambda_{F_\gamma}^{\rm st,\bullet}
=\Lambda_{F_\gamma}^\bullet\cap
\Lambda^{\rm st}_{F_\gamma}.$
\end{enumerate}
\end{Lemma}

\begin{proof}
We have $F_\xi\subseteq F_\gamma$ by Lemma \ref{lem:words} applies to $F_\gamma$. 
By Lemmas \ref{lem:words}, \ref{lem:leading}, the image
of $F_\xi^{\leqslant P}$ in 
$$F_\gamma^{\leqslant P}/F_\gamma^{<P}=V(P)$$ 
is the span of the $w(\xi_0)A_0(u)$ with $b(w,u)=P$. So (3) says
$$F_\xi^{\leqslant P}+F_\gamma^{<P}=F_\gamma^{\leqslant P},\quad\forall P,$$ 
which by increasing induction on $P$ in each colour is equivalent to $F_\xi=F_\gamma$. 
Hence (1)$\Leftrightarrow$(3).
Condition (3) only depends on $\xi_0$, which gives (2)$\Leftrightarrow$(3). 
Claim (3) is the spanning condition defining $\Lambda_{F_\gamma}^{\rm st,\bullet}$ in \eqref{LB}, proving (b).
\end{proof}

\begin{Lemma}\label{lem:mu}\hfill
\begin{enumerate}[label=$\mathrm{(\alph*)}$,leftmargin=8mm,itemsep=1mm]
\item 
Each homogenous matrix coefficient $V_i(P)\to V_i(Q)$ of the moment map endomorphism $\mu(\xi)$ satisfies $P-Q\geqslant 2$.
\item
The matrix-component part of centred degree $2$ of $\mu(\xi)$ is $\mu(\xi_0)$.
\item
The map $h$ in \eqref{h} takes 
$\Lambda_{F_\gamma}$ to $\Lambda_{F_\gamma}^\bullet$ and $\Lambda^{\rm st}_{F_\gamma}$ to $\Lambda_{F_\gamma}^{\rm st,\bullet}$.
\end{enumerate}
\end{Lemma}

\begin{proof}
For $\xi=(x,y)$,  the matrix $\mu(\xi)$ is a sum of products of a unit of $x$ and a unit of $y$.
The centred degrees add, and the minimal centred degrees add up to $2$ in each case
($1+1$ for $B_hB^*_h$ and $B^*_hB_h$, $0+2$ for $AA^*$). Equality holds iff both units have excess $0$. The last claim follows with
Lemma \ref{lem:stab}.
\end{proof}

For $t\in\bbC$, let $c_t$ be the linear operator on $L_{F_\gamma}$ 
which multiply each matrix unit of excess $e$ by $t^e$. 

\begin{Lemma}\label{lem:contr}\hfill
\begin{enumerate}[label=$\mathrm{(\alph*)}$,leftmargin=8mm,itemsep=1mm]
\item 
We have $c_1=\mathrm{id}$ and $c_0=h$ in \eqref{h}.
\item
The map $(t,\xi)\mapsto c_t(\xi)$ is algebraic.
\item
The operator $c_t$ commutes with $\dG_\gamma$ and $c_t(\Lambda^{\rm st}_{F_\gamma})\subseteq\Lambda^{\rm st}_{F_\gamma}$.
\end{enumerate}
\end{Lemma}

\begin{proof}
Only the last claim needs proof.
For $t=0$ it is Lemma \ref{lem:mu}. Let $t\neq0$, let $\lambda(t)\in\dG_\gamma$ act on $V(P)$ by $t^{-P}$ and
$\lambda_W(t)$ on $W(L)$ by $t^{-L}$. 
Conjugation by $(\lambda(t),\lambda_W(t))$ multiplies a unit $V(P)\to V(Q)$ of centred degree $d$ by $t^d$. 
Thus, we have
$$c_t=\sigma_t\circ(\lambda(t),\lambda_W(t)),$$ 
where $\sigma_t$ multiplies $B_h,B_h^*$ by $t^{-1}$, $A$ by $1$ and $A^*$ by $t^{-2}$. These operations
preserve $\dN_{F_\gamma}$ and $\dN_{F_\gamma}^\perp$.
Further, we have  
$$F_{c_t(\xi)}=\lambda(t)F_\xi,$$ 
since $\lambda_W(t)$ preserves each $\dW^{\leqslant L}$ and nonzero rescalings
of $B,$ $B^*$ do not change spans. Finally, 
$$\lambda(t)F_\gamma=F_\gamma.$$
By the proof of Lemma \ref{lem:mu}, an elementary matrix component of the moment map $\mu$ relatively to the decomposition
$\dV_i=\bigoplus_PV_i(P)$ of centred degree $d$ (source degree minus target degree) is multiplied by
$t^{d-2}$. Thus the condition $\mu=0$ is preserved.
\end{proof}

\begin{Proposition}\label{prop:coh}\hfill
\begin{enumerate}[label=$\mathrm{(\alph*)}$,leftmargin=8mm,itemsep=1mm]
\item 
The restrictions to Levi subgroups and $\Lambda_{F_\gamma}^{\rm st,\bullet}$ give isomorphisms of $H^*_{P_\gamma}$-algebras
\[
H^*(\mathfrak L_{\rho,\gamma},\bbC)\cong
H^*_{P_\gamma}(\Lambda^{\rm st}_{F_\gamma},\bbC)\xrightarrow{\cong}
H^*_{\dG_\gamma}(\Lambda^{\rm st}_{F_\gamma},\bbC)\xrightarrow{\cong}
H^*_{\dG_\gamma}(\Lambda_{F_\gamma}^{\rm st,\bullet},\bbC)=H^*(\M_{\dQ}^\bullet(v_\gamma,w_\rho),\bbC).
\]
\item
Before completion, the local Kirwan map $\kappa_\gamma$ in \eqref{LKM1} is identified with a map
\begin{align}\label{LKM2}
\kappa_\gamma:H^*_{P_\gamma}\to H^*(\M_{\dQ}^\bullet(v_\gamma,w_\rho),\bbC)\end{align}
\end{enumerate}
\end{Proposition}

\begin{proof}
The first map is an isomorphism because
$P_\gamma/\dG_\gamma\cong U_\gamma$ is contractible.
For the inclusion $j:\Lambda_{F_\gamma}^{\rm st,\bullet}\hookrightarrow
\Lambda_{F_\gamma}^{\rm st}$, Lemmas~\ref{lem:stab}-\ref{lem:contr}
give $hj=\mathrm{id}$ and the $\dG_\gamma$-equivariant homotopy
$c_t:j h\cong\mathrm{id}$. This proves the second isomorphism.
The quotient identification is Lemma~\ref{lem:free}(c).
All maps preserve characteristic classes, which descend to the
Chern classes of $\calV_i(P)$, proving (b).
\end{proof}

\begin{Remark}\label{rem:A2-not-isomorphic}
The cohomological isomorphism in Proposition \ref{prop:coh} does not, in general,
come from an isomorphism of algebraic varieties. If $Q$ is of type $A_2$, with vertices $1,2$ and graded spaces
\[
 V_1(0)=V_2(1)=\bbC,\quad
 W_1(0)=W_1(2)=W_2(1)=\bbC,
\]
all other graded pieces being zero, then the possible nonzero maps are
$$A_1:W_1(0)\to V_1(0),\quad
A_2:W_1(2)\to V_1(0),\quad
A_3:W_2(1)\to V_2(1),\quad
B^*:V_2(1)\to V_1(0).$$
The moment map equation is automatic, and stability amounts to
$A_1A_3\ne0$. Since $P_\gamma=(\mathbb C^\times)^2$, normalizing
$A_1=A_3=1$ leaves the free coordinates $(A_2,B^*)$. Hence
$\mathfrak L_{\rho,\gamma}\cong\mathbb A^2.$
In the homogeneous locus, the map $A_2$ has positive excess and vanishes,
while $B^*$ has minimal centred degree and remains. Thus
$\M_{\dQ}^\bullet(v_\gamma,w_\rho)\cong\mathbb A^1.$
The two varieties are homotopy equivalent but not isomorphic.
\end{Remark}

\subsection{Equivariant cohomology of the stable Lagrangian}\label{subsec:Zst-coh}
The geometry and irreducibility of $Z_\gamma^{\rm st}$ were established
in \S\ref{subsec:Zst}. We now compare its equivariant cohomology with
that of the graded quiver variety in the symmetric setting.
Recall that in symmetric types the local Kirwan map defined in \eqref{LKM1} is a map
$\kappa_\gamma: \widehat H^*_{\dG_\gamma}\to H^*(\M_{\dQ}^\bullet(v_\gamma,w_\rho),\bbC)$.

\begingroup\emergencystretch=1em
\begin{Proposition}\label{prop:cohZ}
For arbitrary finite type Cartan pairs $(C,D)$, the following hold.
\begin{enumerate}[label=$\mathrm{(\alph*)}$,leftmargin=8mm]
\item
$H^*(\mathfrak L_{\rho,\gamma},\bbC)\cong H^*_{\dG}(Z^{\rm st}_\gamma,\bbC)\cong H^*_{\dG}(\Lambda^{\rm st}_\gamma,\bbC)
\cong H^*_{P_\gamma}(\Lambda^{\rm st}_{F_\gamma},\bbC).$
\item
The  pull-back along the map
$\Lambda^{\rm st}_\gamma\to\scrX_\gamma,$ $(F,x,y)\mapsto(F,x)$ 
is the restriction 
$H^*_{P_\gamma}\to H^*_{P_\gamma}(\Lambda^{\rm st}_{F_\gamma},\bbC).$
If $D=1$ this is the restriction of the local Kirwan map $\kappa_\gamma$  to  
$H^*_{P_\gamma}\cong H_{\dG_\gamma}^*$. 
\end{enumerate}
\end{Proposition}
\endgroup

\begin{proof}
Use Lemma \ref{lem:Z}, the induction isomorphism
$$H^*_{\dG}(\dG\times_{P_\gamma}Y,\bbC)\cong H^*_{P_\gamma}(Y,\bbC),$$ 
and the isomorphism
$H^*_{P_\gamma}(\dN_{F_\gamma},\bbC)\cong H^*_{P_\gamma}.$
The last claim is Proposition $\ref{prop:coh}$.
\end{proof}

\section{Type ADE}\label{sec:main}

Recall the map $r_\gamma:\widehat R_\gamma\to \widehat\B_\gamma$ in \eqref{RG}.
The following theorem is proved in \S\ref{subsec:proof-main-map}, \S\ref{subsec:surj} and \S\ref{subsec:inj}.

\begin{Theorem}\label{thm:main}
Let $Q$ be of type ADE. 
For each $(\rho,\gamma)$ the following hold.
\begin{enumerate}[label=$\mathrm{(\alph*)}$,leftmargin=8mm]
\item $\Ker(r_\gamma)\subseteq\Ker(\kappa_\gamma)$. 
The local comparison map is a well defined ring homomorphism
$$\Phi_\gamma:\widehat\B_\gamma\to H^*(\M_{\dQ}^\bullet(v_\gamma,w_\rho),\bbC)$$ 
\item $\Phi_\gamma$ is surjective.
\item If  \eqref{C2}   is satisfied, then
$\Phi_\gamma$ is an isomorphism, i.e., the local Hikita conjecture holds.
\end{enumerate}
\qed
\end{Theorem}

\begin{Corollary}\label{cor:points}
Let $Q$ be of type ADE. For every integral pair $(\rho,\gamma)$,
\[
 \M_{\dQ}^\bullet(v_\gamma,w_\rho)\ne\emptyset
 \quad\Longleftrightarrow\quad
 \zeta^{2\gamma}\in\Spec\B(\calA^\rho_\mu).
\]
No dimension hypothesis \eqref{C2} is required.
\end{Corollary}

\begin{proof}
By Lemmas~\ref{lem:stab} and \ref{lem:Z}, the left side is equivalent
to $Z_\gamma^{\rm st}\ne\emptyset$. 
Theorem \ref{thm:main} implies the forward implication. The other one follows from
Proposition~\ref{prop:microlocal-basis}(a).
\end{proof}

\begin{Corollary}\label{cor:dim}
Let $Q$ be of type ADE. The following are equivalent: 
\begin{enumerate}[label=$\mathrm{(\alph*)}$,leftmargin=8mm]
\item 
The local Hikita conjecture holds for $(\rho,\gamma)$.
\item
$\Phi_\gamma$ is injective.
\item
$\dim\widehat\B_\gamma\leqslant\dim H^*(\M_{\dQ}^\bullet(v_\gamma,w_\rho),\bbC)$. 
\end{enumerate}
They hold if the global Hikita conjecture \eqref{GHC} or the bound \eqref{C2} hold.
\end{Corollary}

\begin{proof}
For a surjection onto a finite dimensional space, the upper dimension bound in (c) forces finiteness of the source
and is equivalent to injectivity. The
last claim is Proposition \ref{prop:global-local}.
\end{proof}

\begin{Remark}\hfill
\begin{enumerate}[label=$\mathrm{(\alph*)}$,leftmargin=8mm]
\item 
The local maps $\Phi_\gamma$ determine a global equivariant comparison map. 
The argument is  carried out for arbitrary quivers and $D=1$ in \S\ref{sec:general}, in Theorem~\ref{thm:gen-global-map}.
In particular, in types ADE, if \eqref{C2} holds, it implies that the global Kirwan map $\kappa_{F,\hbar}$ in \eqref{GKM} factors through 
the B-algebra $\B^\theta_{F,\hbar}$ in \eqref{BTFH} and gives an isomorphism of graded
$R_{F,\hbar}$-algebras
\[
\Phi:\ \B^\theta_{F,\hbar}\xrightarrow{\cong}H^*_{F\times\bbC^\times}(\M_Q(v,w),\bbC).
\]
\item
For non symmetric finite types, Appendix~\ref{app:ns-global-technical}
constructs smooth moduli spaces $\mathfrak M_N$ and their
truncation tower. 
Proposition~\ref{prop:ns-global-torsion} gives a global comparison modulo equivariant torsion.
Its local specializations are the maps of Theorem~\ref{thm:ns-map}.
The local $B_2$ and global $G_2$ obstructions to surjectivity are given in
Propositions~\ref{prop:B2-local-kirwan} and \ref{prop:G2-global-kirwan}.
\end{enumerate}
\end{Remark}

\section{Morse modules}\label{sec:morse}

\subsection{Construction and existence}\label{subsec:vc}
We recall the notation and the facts on vanishing cycles used below, see \cite[\S4.2]{Dim04} and \cite[\S8.6]{KS90} for more details. 
Let $X$ be a complex
manifold or a smooth complex variety, $f:X\to\bbC$ a holomorphic function.
Let $$i:X_0\hookrightarrow X$$ 
be the inclusion of the zero set $X_0=f^{-1}(0)$.
For a constructible complex $\calE$ on $X$, let $\varphi_f(\calE)$
denote {vanishing cycles} complex of $\calE$ along $f$ on $X_0$.
Recall that 
$$\varphi_f:D^b_c(X)\to D^b_c(X_0)$$ is a triangulated functor, hence
$\calE\mapsto H^*(\varphi_f(\calE)_x)$ is an additive functor for each point $x\in X_0$.
If $\calE$ is locally constant then $\varphi_f(\calE)$ is supported on the critical locus of $f$ in $X_0$,
and if $\pi:Y\to X$ is proper, then the proper base change yields
$$\varphi_f(\pi_*\calE)\cong\pi_*\varphi_{f\circ\pi}(\calE)
,\quad
\varphi_f(\pi_*\calE)_x\cong R\Gamma\big(\pi^{-1}(x),\varphi_{f\circ\pi}(\calE)\big).$$
Further, if $X$ is  of dimension $n$ and $x\in X_0$ is a nondegenerate critical point of $f$, then the Morse lemma
implies that the stalk of $\varphi_f(\bbC_X)$ at $x$ is \(\varphi_f(\mathbb C_X)_x\cong\mathbb C[1-n]\). 
Finally, if $f$ is only defined only on an open neighbourhood $U$ of $x$, we abbreviate
$$\varphi_f(\calE)_x=\varphi_f(\calE|_U)_x.$$

Let $\xi=(x,y)$ in $\Lambda_{F_\gamma}^{\rm st}$. 
For a generic quadratic form $H$ on $\dN$ we define the function
\begin{align}\label{fH}
f:\dN\to\mathbb C,\quad
x'\mapsto\langle y,x'-x\rangle+\frac12H(x'-x,x'-x).\end{align} 
It is a polynomial function with 
$$f(x)=0,\quad df_x=y.$$
We define the Morse fibre to be
\begin{align}\label{MF}MF^\xi_\eta=\varphi_f(\mathscr L_\eta)_x[-d_\gamma-1].\end{align}
Set
\begin{align*}
\mathring M(\rho,\gamma)_\eta=H^*(MF^\xi_\eta),\quad
\mathring M(\rho,\gamma)=\bigoplus_\eta \mathring M(\rho,\gamma)_\eta.
\end{align*}
The shift $[-d_\gamma-1]$  places the one dimensional $\gamma$-Morse fibre in degree zero,
as proved in Lemma~\ref{lem:Morse}(e).

\begin{Lemma}\label{lem:Morse}
Let $\xi=(x,y)$ in $\Lambda_{F_\gamma}^{\rm st}$ and let $f:\dN\to\mathbb C$ be as in \eqref{fH}.
Then:
\begin{enumerate}[label=$\mathrm{(\alph*)}$,leftmargin=8mm,itemsep=1mm]
\item $\mathring M(\rho,\gamma)$ is a locally unital $\bR$-module, with
$e_\eta \mathring M(\rho,\gamma)=\mathring M(\rho,\gamma)_\eta$ finite dimensional for all $\eta$.
\item
The augmentation ideal of $H^*_{\dG}$ acts locally nilpotently.
\item $e_\eta \mathring M(\rho,\gamma)\neq0$ only if $\eta\leqslant\gamma$.
\item The graph of $df$ is transverse at $(x,y)$ to the smooth 
Lagrangian germ $p(\Lambda_\gamma)$.
\item $e_\gamma \mathring M(\rho,\gamma)$ is one dimensional, concentrated in degree zero.
\end{enumerate}
\end{Lemma}

\begin{proof}
First, by proper base change we have
\[
\mathring M(\rho,\gamma)_\eta
 =H^*\bigl(\pi_\eta^{-1}(x),\,
 \varphi_{f\circ\pi_\eta}(\omega_{\scrX_\eta})
 |_{\pi_\eta^{-1}(x)}[-d_\gamma-1]\bigr).
\]

(a), (b) An element in
$\operatorname{Ext}_{\dG}^p(\mathscr L_\eta,\mathscr L_\nu)$ gives a derived
morphism $\mathscr L_\eta\to\mathscr L_\nu[p]$. Applying
$(\varphi_f(-))_x$ and taking cohomology yields a degree $p$ map
$\mathring M(\rho,\gamma)_\eta\to\mathring M(\rho,\gamma)_\nu$.
Functoriality preserves composition, defining a left $\bR$-action.
The idempotent $e_\eta$ projects onto the $\eta$-piece.
Since
$\mathring M(\rho,\gamma)$ is a direct sum, this action is locally unital.
Each $\mathring M(\rho,\gamma)_\eta$ is finite dimensional because $\mathscr L_\eta$ is
constructible. Hence, the classes from $H^{>0}_{\dG}$ act
locally nilpotently. 
The fibre complex is supported at critical points $(F,x)$ of
$f\circ\pi_\eta$. Since $df_x=y$, we have
$d(f\circ\pi_\eta)_{(F,x)}=0$ exactly when $y$ annihilates
$\operatorname{Im}(d_{(F,x)}\pi_\eta)$, equivalently when
$(F,x,y)\in\Lambda_\eta$ by Lemma~\ref{lem:LZL}(a).

(c) For such $F$ we have $F\supseteq F_\xi=F_\gamma$ by Lemma \ref{lem:words}, so $\eta\leqslant\gamma$.

(d) For $\eta=\gamma$, the same argument shows that $F=F_\gamma$.
Thus  $e_\gamma \mathring M(\rho,\gamma)$ is the stalk at the single critical point $(F_\gamma,x)$. By Lemma~\ref{lem:smooth}, $p$ is immersive there with image of
dimension $\dim\dN$; this image is isotropic by Lemma~\ref{lem:LZL}(b),
hence is a smooth Lagrangian germ $L$.
The critical points of $f\circ\pi_\gamma$ near $(F_\gamma,x)$ correspond to $L\cap\Gamma_{df}$.
The tangent space of $\Gamma_{df}$ at $(x,y)$ is $\{(v,Hv)\}$.
For generic symmetric $H$ it is transverse to $T_{(x,y)}L$, since both the transverse Lagrangian subspaces and
the graphs form dense open subsets of the Lagrangian Grassmannian. 

(e) By (d) the point $(F_\gamma,x)$ is a nondegenerate critical point.
Since $\omega_{\scrX_\gamma}=\bbC_{\scrX_\gamma}[2d_\gamma]$,
the Morse lemma shows that the unshifted vanishing-cycle
stalk is $\bbC[d_\gamma+1]$. Hence
$e_\gamma\mathring M(\rho,\gamma)\cong\bbC$ in degree zero. 
\end{proof}

We call the function $f$ a \emph{Morse function} and the $\mathbf R$-module $\mathring M(\rho,\gamma)$ a \emph{Morse module}.
To emphasize the dependence on the point $\xi$ in $Z_\gamma^{\rm st}$ we may write
$$\mathring M(\rho,\gamma)^\xi=\mathring M(\rho,\gamma).$$
The function $f$ is omitted from the notation, since different generic choices give isomorphic modules, 
though not canonically.
Note that an element of
$e_\nu\bR^k e_\eta$ acts with cohomological degree $k$.

\begin{Theorem}\label{thm:point}Assume  $(\rho,\gamma)$ satisfies condition \eqref{bfE}.\hfill
\begin{enumerate}[label=$\mathrm{(\alph*)}$,leftmargin=8mm]
\item
$\zeta^{2\gamma}\in\Spec\B(\calA^\rho_\mu)$. 
\item
$Z_{\rho,\gamma}$ is $\ell$-dominant $\Rightarrow$ Theorem $\ref{Theorem:431}$ holds for $(\rho,\gamma)$.
\end{enumerate}
\end{Theorem}

\begin{proof}
Choose $\xi\in\Lambda_{F_\gamma}^{\rm st}$ and the Morse module
$\mathring M(\rho,\gamma)^\xi$. Lemma~\ref{lem:Morse} gives the local
unitality, the nilpotence of the augmentation ideal, and a nonzero
one dimensional $\gamma$-piece. Every nonzero piece has
$\eta\leqslant\gamma$, so, by Remark~\ref{rem:Up}(c), there is no $\eta$-piece with
$\eta\in\Up(\gamma)$. Similalry, Lemma~\ref{lem:noncyc}(b) excludes all $\eta$-pieces with $\eta\in\Cyc(\rho)$.
Proposition~\ref{prop:crit} therefore applies and proves both assertions.
\end{proof}

\begin{Proposition}\label{prop:topdown}
Let $Z_{\rho,\gamma}$ be $\ell$-dominant and $D=1$. 
\begin{enumerate}[label=$\mathrm{(\alph*)}$,leftmargin=8mm]
\item
We have $\Lambda_{F_\gamma}^{\rm st,\bullet}\neq\emptyset$ and $\M_{\dQ}^\bullet(v_\gamma,w_\rho)\neq\emptyset$.
\item
Theorem $\ref{Theorem:431}$ holds for $(\rho,\gamma)$.
\end{enumerate}
\end{Proposition}

\begin{proof}
Assume $Z_{\rho,\gamma}$ be $\ell$-dominant.
We look for homogeneous points of $\Lambda_{F_\gamma}^{\rm st}$ with $A^*=0$. Write 
$$w_i(L)=\dim W_i(L),\quad v_i(P)=\dim V_i(P).$$ 
For each $i$ and $m$ put $M_i(m+1)=W_i(m)\oplus\textstyle\bigoplus_{j-i}V_j(m+1)$ and
\begin{align*}
\tau_i^{(m)}&=(A_i,B_{j\to i}):\ M_i(m+1)\to V_i(m),\\
\sigma_i^{(m)}&=(0,\varepsilon B_{i\to j}):\ V_i(m+2)\to M_i(m+1).
\end{align*}
each component of $B$ is a component of  one $\tau$ or of  one $\sigma$.
Further, the equation
$\mu^\bullet=0$ is the system of equations $\tau^{(m)}_i\sigma^{(m)}_i=0$. 
We define the maps $\tau^{(m)}$ and $\sigma^{(m)}$
by decreasing induction on $m$. At step $m$ all $\sigma^{(m)}_i$ are known,
since their components have level $\geqslant m+1$.
Choose $\tau^{(m)}_i$ generic among the maps vanishing on $\operatorname{Im}(\sigma^{(m)}_i)$. It is
surjective as soon as 
$$\dim M_i(m+1)-\dim V_i(m+2)\geqslant\dim V_i(m),$$ 
i.e., if we have
\[
w_i(m)+\textstyle\sum_{j-i}v_j(m+1)-v_i(m+2)-v_i(m)\geqslant0.
\]
Note that this is also
 the exponent of $Z_{i,m+1}$ in $Z_{\rho,\gamma}$ in centred positions, which is nonnegative by $\ell$-dominance. 
 The surjectivity of $\tau^{(m)}$ at each level $m$ implies, by decreasing induction on $m$, 
 that $V_i(m)$ is spanned by words in the $B$'s applied to $A(\dW)$.
Each such homogeneous word has bound $m$. Hence Lemma~\ref{lem:words} gives $F_\xi=F_\gamma$.
Thus the point lies in both $\Lambda_{F_\gamma}^{\rm st,\bullet}$ and $\Lambda^{\rm st}_{F_\gamma}$, proving (a).
Part (a) gives condition \eqref{bfE}.
Apply Theorem \ref{thm:point} and Lemma \ref{lem:Bcrit}.
This proves (b).
\end{proof}

\begin{Remark}
Proposition \ref{prop:topdown} can  be deduced from Nakajima's theorem that the regular stratum
$\M_0^{\bullet\rm reg}(v,w)$ is nonempty iff $(v,w)$ is $\ell$-dominant \cite{Nak01}, \cite[Prop.~2.4.9]{KQ14}.
\end{Remark}

\subsection{Jordan-H\"older multiplicities of Morse modules}\label{sec:JHMM}
In this subsection $(C,D)$ is of finite type and $(\rho,\gamma)$
satisfies \eqref{bfE}. Fix $\xi=(x,y)$ as in Lemma~\ref{lem:Morse}.
The preceding results make $\mathring M(\rho,\gamma)^\xi$ a smooth
locally unital $\bT^\wedge$-module with finite dimensional pieces.

Recall that $\Gamma_\gamma$ is the set of types $\eta$ in the fixed weight block $C_\gamma$
for which the simple module $L(\rho,\eta)$ with the highest $\ell$-weight
$Z_{\rho,\eta}$ belongs to $\scrO_\mu^\rho$, see \eqref{CCG}. 
The set $\Gamma_\gamma$ is finite because the simple objects of $\scrO_\mu^\rho$
 are indexed by a subset of
the dual canonical basis of a finite dimensional representation by
\cite[Prop.~4.4.1, Thm.~4.4.2]{VV25}.

For each type $\eta$, we fix a decomposition
\begin{equation}\label{eq:morse-IC-summands}
 \mathscr L_\eta\cong\bigoplus_{\lambda\in\mathcal I_\eta}\ \bigoplus_{s\in\bbZ}
 \IC_\lambda[-s]\otimes E_{\eta,\lambda,s}.
\end{equation}
Here the $\IC_\lambda={\rm IC}(\overline{S_\lambda},\mathcal E_\lambda)$ are
pairwise nonisomorphic simple $\dG$-equivariant perverse sheaves.
The stratum $S_\lambda$ is smooth and locally closed, and $\mathcal E_\lambda$
is an irreducible equivariant local system. We use a common set of labels
$\lambda$, with $E_{\eta,\lambda,s}=0$ when $\lambda\notin\mathcal I_\eta$.
For later use, for any finite set $J\subset C_\gamma$ we put
$\mathcal I_J=\bigcup_{\eta\in J}\mathcal I_\eta.$
Taking direct sums in \eqref{eq:morse-IC-summands} gives
\begin{equation}\label{eq:morse-IC-decomposition}
 \mathscr L_J\cong\bigoplus_{\lambda\in\mathcal I_J}\ \bigoplus_{s\in\bbZ}
 \IC_\lambda[-s]\otimes E_{J,\lambda,s},\quad
 E_{J,\lambda,s}=\bigoplus_{\eta\in J}E_{\eta,\lambda,s}.
\end{equation}
Write $E_{J,\lambda}=\bigoplus_s E_{J,\lambda,s}$ and
$E_{\eta,\lambda}=\bigoplus_sE_{\eta,\lambda,s}$.

\begin{Theorem}\label{thm:morse-JH}
\hfill
\begin{enumerate}[label=$\mathrm{(\alph*)}$,leftmargin=8mm]
\item\label{item:morse-JH-length}
The module $\mathring M(\rho,\gamma)^\xi$ has finite length in
$\scrO_\mu^\rho$ and $[\mathring M(\rho,\gamma)^\xi:L(\rho,\gamma)]=1.$

\item\label{item:morse-JH-global}
For each $\eta\in\Gamma_\gamma$ with
$[\mathring M(\rho,\gamma)^\xi:L(\rho,\eta)]\ne0$ there is a
unique $\lambda(\eta)\in\mathcal I_{\Gamma_\gamma}$ with
$$e_{\Gamma_\gamma} L(\rho,\eta)\cong E_{\Gamma_\gamma,\lambda(\eta)}$$ as an $\bR_{\Gamma_\gamma}$-module.
The type $\eta$ is the unique maximal charge type
$\nu\in\Gamma_\gamma$ with $E_{\nu,\lambda(\eta)}\ne0$, and
\begin{equation}\label{eq:morse-JH-global}
 [\mathring M(\rho,\gamma)^\xi:L(\rho,\eta)]
 =\dim H^*\bigl(\varphi_f(\IC_{\lambda(\eta)})_x\bigr).
\end{equation}

\item\label{item:morse-JH-CC}
Put $\Lambda=\overline{Z_\gamma^{\rm st}}\subset T^*\dN$.
For $\xi\in Z_\gamma^{\rm st}$ generic and $\eta$ as in $\mathrm{(b)}$ we have
\begin{equation}\label{eq:morse-JH-CC}
 [\mathring M(\rho,\gamma)^\xi:L(\rho,\eta)]
 =m_\Lambda(\IC_{\lambda(\eta)}),\quad
 \operatorname{CC}(\IC_{\lambda(\eta)})=\sum_{\Lambda'}
 m_{\Lambda'}(\IC_{\lambda(\eta)})[\Lambda'].
\end{equation}
\end{enumerate}
\end{Theorem}

\begin{proof}
For any finite set $J\subset C_\gamma$,
the projection onto the perverse Ext degree zero component defines a surjective algebra map
\begin{equation}\label{eq:morse-semisimple-quotient}
 q_J:\bR_J\longrightarrow
 \prod_{\lambda\in\mathcal I_J}\End_{\bbC}(E_{J,\lambda}).
\end{equation}
Indeed, negative perverse Ext groups vanish, and degree zero maps
between distinct simple perverse sheaves vanish. Denote by $E_{J,\lambda}$ again
the simple $\bR_J$-module $E_{J,\lambda}$ pulled back through $q_J$.
Then, 
\begin{equation}\label{eq:morse-simple-eta-pieces}
 \dim e_\eta E_{J,\lambda}=\dim E_{\eta,\lambda},\quad\forall\eta\in J.
\end{equation}
We claim that, for every finite $J$, we have
\begin{equation}\label{eq:morse-JH-formula}
 [e_J\mathring M(\rho,\gamma)^\xi:E_{J,\lambda}]
 =\sum_{k\in\bbZ}\dim H^k\bigl(\varphi_f(\IC_\lambda)_x\bigr).
\end{equation}
Hence, at a generic point $\xi$ of an irreducible component
$\Lambda$ of $\operatorname{SS}(\mathscr L_J)$, we have
\begin{equation}\label{eq:morse-JH-CC-corner}
 [e_J\mathring M(\rho,\gamma)^\xi:E_{J,\lambda}]=m_\Lambda(\IC_\lambda)
 ,\quad
 \forall\lambda\in\mathcal I_J.
\end{equation}

To prove the claim, vanishing cycles and \eqref{eq:morse-IC-decomposition} yield
\begin{equation}\label{eq:morse-JH-splitting}
 e_J\mathring M(\rho,\gamma)^\xi
 \cong\bigoplus_{\lambda,s}E_{J,\lambda,s}\otimes
 H^*\bigl(MF^\xi_\lambda[-s]\bigr).
\end{equation}
The kernel of $q_J$ has strictly positive perverse Ext degree, hence
raises the degree in $H^*(\varphi_f(\IC_\lambda)_x)$ after
removing the shifts. Only finitely many degrees occur, so this kernel
acts nilpotently. Every simple factor therefore factors through
$q_J$. For each $\lambda$, projection onto one copy of
$\IC_\lambda[-s]$ defines an idempotent $\epsilon_\lambda$
whose image under $q_J$ has rank one in the $\lambda$-factor.
Since $\dim\epsilon_\lambda E_{J,\mu}=\delta_{\lambda\mu}$,
exactness of applying $\epsilon_\lambda$ to a composition series,
together with \eqref{eq:morse-JH-splitting}, proves
\eqref{eq:morse-JH-formula} and \eqref{eq:morse-JH-CC-corner}.
The common shift in \eqref{eq:morse-JH-splitting} is essential for
this graded decomposition, but does not change the total dimensions
in the ungraded multiplicity formulas.

(a) For any subquotient $N$ of $\mathring M(\rho,\gamma)^\xi$, set
\[
 D(N)=\sum_{\eta\in\Gamma}\dim e_\eta N.
\]
This is finite because $\Gamma$ is finite and each weight piece
is finite dimensional. We claim that $N\ne0$ implies $D(N)>0$.
Indeed, choose $0\ne v\in e_\alpha N$. The cyclic
$\bT^\wedge$-module generated by $v$ is finitely generated,
smooth and locally unital, with finite dimensional pieces.
By Theorem~\ref{thm:V}(c), it belongs to $\scrO_\mu^\rho$ and has
a simple quotient $L(\rho,\nu)$ with $\nu\in\Gamma$.
Exactness of $e_\nu$ gives $e_\nu N\ne0$, proving the claim.
For every strict chain of submodules
$0=N_0\subsetneq\cdots\subsetneq N_a=\mathring M(\rho,\gamma)^\xi$,
\[
 D(\mathring M(\rho,\gamma)^\xi)
 =\sum_{j=1}^a D(N_j/N_{j-1})\geqslant a.
\]
Thus the module has finite length, and its finite composition series
shows it is finitely generated. By Theorem~\ref{thm:V}(c) it
belongs to $\scrO_\mu^\rho$. By Lemma~\ref{lem:Morse}, its nonzero
types are $\leqslant\gamma$ and its $\gamma$-piece has dimension one.
Only $L(\rho,\gamma)$ can contribute to this piece, and its highest
piece is one dimensional. Hence its multiplicity is one.

(b) 
For each $\eta\in\Gamma_\gamma$, the $e_{\Gamma_\gamma}\bR e_{\Gamma_\gamma}$-module 
$e_{\Gamma_\gamma} L(\rho,\eta)$
is simple, since its highest piece is nonzero.
These restrictions are pairwise nonisomorphic: their unique
maximal charge types are the distinct $\eta$'s. Since restriction
by $e_{\Gamma_\gamma}$ is exact, every composition factor occurring in
$\mathring M(\rho,\gamma)^\xi$ restricts to one of the
$E_{\Gamma_\gamma,\lambda}$, with a unique label $\lambda(\eta)$.
Formula~\eqref{eq:morse-JH-formula} for the corner $\Gamma_\gamma$
then proves \eqref{eq:morse-JH-global}. The support identity
\eqref{eq:morse-simple-eta-pieces} and the highest weight property
identify $\eta$ as the unique maximal charge type in that support.

(c) 
By Lemma~\(\ref{lem:Morse}\)(d),(e), for every \(\xi=(x,y)$ in $Z_\gamma^{\rm st}\), there is a function \(f\) with \(df_x=y\) and \(\varphi_f(\mathscr L_\gamma)_x\ne0\). Hence \(Z_\gamma^{\rm st}\subseteq\operatorname{SS}(\mathscr L_\gamma)\).
Thus $Z_\gamma^{\rm st}\subseteq\operatorname{SS}(\mathscr L_{\Gamma_\gamma})$.
Lemmas~\ref{lem:Z} and \ref{lem:stable-irreducible}
show its closure $\Lambda$ is irreducible Lagrangian.
As the microsupport is a union of Lagrangian components,
$\Lambda$ is itself a component. Now apply
\eqref{eq:morse-JH-CC-corner} with $J=\Gamma_\gamma$ and (b).
\end{proof}

\begin{Remark}
Formula \eqref{eq:morse-JH-formula} holds for every $\xi$,
without the genericity assumption in \eqref{eq:morse-JH-CC-corner}.
The Morse module need not be generated by its one dimensional
$\gamma$-piece, even in type $A_1$. Indeed, at the stable
point $\xi=\xi_0$ of Appendix~\ref{ex:coMorse-A1-nonsimple}, with
$\gamma=2$ and $\eta=0$, the  Morse fibres are
$MF^\xi_\gamma\cong\bbC$ and $MF^\xi_\eta\cong\bbC[1]$.
Thus the $\eta$-piece of the Morse module $\mathring M(\rho,\gamma)^{\xi_0}$ 
is one dimensional of degree $-1$, while the $\gamma$-piece is one dimensional of
degree $0$. On the other hand, with the notation in \S\ref{ex:coMorse-A1-nonsimple},
we have $\mathscr L_\gamma=\bbC_Y[2]$ with
$Y=\{a_0=0\}\subset\dN=\bbC^2$, and $\mathscr L_\eta=\bbC_{\dN}[4]$.
So, the purity gives
\[
 e_\eta\bR^{-1}e_\gamma
 =\operatorname{Ext}^{-1}_{\dG}(\mathscr L_\gamma,\mathscr L_\eta)
 \cong H^{-1}_{\dG}(Y,\bbC)=0.
\]
\end{Remark}

\section{Microlocal modules}\label{sec:micro}

\subsection{Microlocal geometry}\label{subsec:micro-inputs}
Let $X$ be a complex manifold, $\pi:T^*X\to X$ its cotangent bundle, and $q_1,q_2:X\times X\to X$ the projections.
The {microlocal Hom} of Kashiwara-Schapira \cite[\S4.4]{KS90} is the bifunctor
\[
\muhom:\ D^b(X)^{\rm op}\times D^b(X)\longrightarrow D^b(T^*X),\quad \muhom(\mathcal F,\mathcal G)=\mu_\Delta\big(R\mathcal Hom(q_2^{-1}\mathcal F,q_1^!\mathcal G)\big),
\]
where $\mu_\Delta:D^b(X\times X)\to D^b(T^*X)$ is Sato's microlocalization 
along the diagonal $\Delta\subset X\times X$, and $T^*_\Delta(X\times X)$ is identified with $T^*X$ by
the first projection. It is a conic complex on $T^*X$ which localizes the sheaf of homomorphisms in the codirections: we have
\[
\pi_*\muhom(\mathcal F,\mathcal G)=R\mathcal Hom(\mathcal F,\mathcal G),
\]
and the restriction of $\muhom(\mathcal F,\mathcal G)$ to the zero section $T^*_XX\cong X$ is $R\mathcal Hom(\mathcal F,\mathcal G)$ 
by \cite[Prop.~4.4.2]{KS90}. 
We use the following facts proved \cite{KS90} for a complex manifold $X$ with an action of an algebraic group $K$ and
constructible complexes. For $K$-equivariant complexes, the functor $\muhom$ is defined on the spaces of the
simplicial resolution as in \cite[\S4.4]{KS90}. It descends, because it commutes with smooth pull-backs.
\begin{enumerate}[label=$\mathrm{(M\arabic*)}$,leftmargin=10mm]

\item\label{M1}  $\operatorname{supp}\muhom(\mathcal F,\mathcal G)\subset \operatorname{SS}(\mathcal F)\cap \operatorname{SS}(\mathcal G)$ by \cite[Cor.~5.4.10]{KS90}.

\item\label{M2}  $R\Gamma(T^*X;\muhom(\mathcal F,\mathcal G))=R\Hom(\mathcal F,\mathcal G)$. There are composition maps
$$\muhom(\mathcal G,\mathcal H)\otimes\muhom(\mathcal F,\mathcal G)\to\muhom(\mathcal F,\mathcal H)$$ 
compatible with the composition of $R\Hom$ and with
restriction to open conic subsets \cite[\S4.4]{KS90}.
The same holds $K$-equivariantly.

\item\label{M3} 
Assume $Z$ is a smooth closed conic Lagrangian submanifold of an open conic $\Omega\subset T^*X$ with
$$\operatorname{SS}(\mathcal F)\cap\Omega\subset Z,$$ such that for each 
$\xi=(x,y)$ in $Z$ there is a holomorphic function $f:X\to\mathbb C$ with $df_x=y$,
the graph of $df$ is transverse to $Z$ at $\xi$, 
and $\varphi_f(\mathcal F)_x$ is of total dimension one, i.e.,  exactly one cohomology group is $\neq 0$ and it is one dimensional.
Then, we say that the sheaf $\mathcal F$ is \emph{simple along $Z$}.
If this holds, by \cite[\S7.5]{KS90}, the assignment $\mathrm{id}\mapsto1$ gives an isomorphism 
$$\muhom(\mathcal F,\mathcal F)|_\Omega\cong\bbC_Z$$ 
This isomorphism is $K$-equivariant if the data are.

\item\label{M4} For a proper map $\pi:Y\to X$ and $\mathcal F$ on $Y$, the functor $\pi_*$ induces a morphism
$$\pi_\mu:Rp_{2!}p_1^{-1}\muhom(\mathcal F,\mathcal F)\to\muhom(\pi_*\mathcal F,\pi_*\mathcal F)$$ 
through 
$$T^*Y\xleftarrow{p_1}Y\times_XT^*X\xrightarrow{p_2}T^*X,$$
compatible with \ref{M2} and with $\pi_*$ on $R\Hom$ \cite[\S4.4, \S5.4]{KS90}. For $\mathcal F=\omega_Y$ we have
$$\muhom(\omega_Y,\omega_Y)=\bbC_{T^*_YY}$$
The map $\pi_\mu$ is compatible with the action of $H^*(Y,\mathbb C)$ by cup product.
\end{enumerate}
\subsection{The map $\Psi_\gamma$}\label{subsec:psi-map}
For $\Omega_\gamma$ and $Z_\gamma^{\rm st}$ as in \eqref{OLZ},
define microlocal restriction by
\begin{equation}\label{PSG}
\begin{split}
\Psi_\gamma: e_\gamma \mathbf{R} e_\gamma =
\Ext^*_{\dG}(\mathscr L_\gamma,\mathscr L_\gamma)
\overset{\text{\ref{M2}}}{=}
\begin{gathered}[t]
H^*R\Gamma_{\dG}\big(T^*\dN\,,\,\muhom(\mathscr L_\gamma,\mathscr L_\gamma)\big)\\
\Big\downarrow\\
H^*_{\dG}\big(\Omega_\gamma\,,\,\muhom(\mathscr L_\gamma,\mathscr L_\gamma)\big).
\end{gathered}
\end{split}
\end{equation}

\begin{Proposition}\label{prop:psi}
Let $(C,D)$ be arbitrary.
\begin{enumerate}[label=$\mathrm{(\alph*)}$,leftmargin=8mm]
\item 
The sheaf $\mathscr L_\gamma$ is simple along $Z^{\rm st}_\gamma$ in the sense of  \ref{M3}.
There is an isomorphism 
$$\muhom(\mathscr L_\gamma,\mathscr L_\gamma)|_{\Omega_\gamma}\cong\bbC_{Z^{\rm st}_\gamma}.$$
The map $\Psi_\gamma$ is a graded algebra homomorphism
$$\Psi_\gamma:
e_\gamma \bR e_\gamma\to H^*_{P_\gamma}(\Lambda^{\rm st}_{F_\gamma},\bbC).$$
\item 
The morphism $\Psi_\gamma$ vanishes on $e_\gamma \bR e_\eta \bR e_\gamma$ for each $\eta\in\Up(\gamma)\cup\Cyc(\rho)$.
\item 
On $H^*_{P_\gamma}\subset e_\gamma \bR e_\gamma$ the map $\Psi_\gamma$ is the restriction $H^*_{P_\gamma}\to H^*_{P_\gamma}(\Lambda^{\rm st}_{F_\gamma},\bbC)$.
\item If $\dim H^*_{P_\gamma}(\Lambda^{\rm st}_{F_\gamma},\bbC)<\infty$  then $\Psi_\gamma$ extends continuously to $e_\gamma \bR^\wedge e_\gamma$ and factors through
$e_\gamma \bT^\wedge e_\gamma$.
\end{enumerate}
\end{Proposition}

\begin{proof}\hfill

(a) By \ref{M1}, we have
$$\operatorname{SS}(\mathscr L_\gamma)\cap\Omega_\gamma\subset Z_\gamma\cap\Omega_\gamma=Z^{\rm st}_\gamma.$$ 
It is a smooth closed conic Lagrangian submanifold of $\Omega_\gamma$
by Lemma \ref{lem:Z}. 
Let $\xi=(x,y)$ in $Z^{\rm st}_\gamma$.
After translation, we have $F_\xi=F_\gamma$. 
Let $f$ be the Morse function in \eqref{fH} with $H$ in the dense open set of quadratic
forms on $\dN$ whose graph  is transverse to $T_\xi Z^{\rm st}_\gamma$. 
This open set depends on $\xi$, and $H$ is chosen separately at each point. 
Recall that $df_x=y$. By Lemma \ref{lem:Morse}(c), the
graph of $df$ is transverse to $Z^{\rm st}_\gamma$ at $\xi$ and the complex
$$\varphi_f(\mathscr L_\gamma)_x$$ has total dimension one. So \ref{M3}
applies, and Proposition \ref{prop:cohZ} yields
$$H^*_{\dG}(\Omega_\gamma\,,\,\bbC_{Z^{\rm st}_\gamma})=
H^*_{P_\gamma}(\Lambda^{\rm st}_{F_\gamma},\bbC)$$ 
The ring
structure given by composition is the cup product, since the isomorphism of \ref{M3} sends $\mathrm{id}$ to $1$ and is compatible with
composition, and restriction commutes with composition by \ref{M2}.

(b) For any elements
$$a\in\Ext_{\dG}(\mathscr L_\eta,\mathscr L_\gamma),\quad b\in\Ext_{\dG}(\mathscr L_\gamma,\mathscr L_\eta),$$
the element $\Psi_\gamma(a\circ b)$ is the composite of the restrictions of $a$
and $b$ to $\Omega_\gamma$, 
viewed as sections of 
$$\muhom(\mathscr L_\eta,\mathscr L_\gamma),\quad\muhom(\mathscr L_\gamma,\mathscr L_\eta).$$
By \ref{M1} these complexes are supported on
$Z^{\rm st}_\gamma\cap Z_\eta$, which is empty by Lemma \ref{lem:Z}.

(c) Let $c\in H^k_{\dG}(\scrX_\gamma,\bbC)$. Cup product by $c$
is a morphism $\omega_{\scrX_\gamma}\to\omega_{\scrX_\gamma}[k]$.
Its pushforward is the corresponding degree $k$ endomorphism
of $\mathscr L_\gamma$. We consider the maps
\[
\begin{aligned}
 q_\gamma^{\rm st}:\Lambda_\gamma^{\rm st}&\longrightarrow\scrX_\gamma,
 & (F,x,y)&\longmapsto(F,x),\\
 p_\gamma^{\rm st}:\Lambda_\gamma^{\rm st}&\xrightarrow{\cong}Z_\gamma^{\rm st},
 & (F,x,y)&\longmapsto(x,y),
\end{aligned}
\]
The second map is the isomorphism of Lemma~\ref{lem:Z}.
In the correspondence of \ref{M4}, the inverse image of the zero
section under $p_1$ is $\Lambda_\gamma$. Over $\Omega_\gamma$, the
microlocal pushforward therefore identifies with
\[
 \bbC_{Z_\gamma^{\rm st}}
 \longrightarrow
 \muhom(\mathscr L_\gamma,\mathscr L_\gamma)|_{\Omega_\gamma}
 \cong\bbC_{Z_\gamma^{\rm st}}.
\]
It sends $1$ to the identity of $\mathscr L_\gamma$, which corresponds
to $1$ under part (a); hence it is the identity. The cup product
compatibility in \ref{M4} now gives
\[
 (p_\gamma^{\rm st})^*\Psi_\gamma(c)
   =(q_\gamma^{\rm st})^*c
 \quad\text{in }H^k_{\dG}(\Lambda_\gamma^{\rm st},\bbC).
\]
By Proposition ~\ref{prop:cohZ}, this is the required restriction map.

(d) By part (c), the positive degree part of $H^*_{\dG}$ acts through $H^{>0}_{P_\gamma}(\Lambda^{\rm st}_{F_\gamma})$, which is nilpotent 
when this ring is finite dimensional.
Thus, the map $\Psi_\gamma$ is continuous and extends. 
By part (b), for $\eta\in\Cyc(\rho)$ it factors through $e_\gamma \bT^\wedge e_\gamma$.
\end{proof}

\subsection{Microlocal modules}
\label{subsec:microlocal-module}
Let $(C,D)$ be arbitrary. Set
\begin{align}\label{microlocal-module}
 M(\rho,\gamma)=
 \bigoplus_\eta
 H_{\dG}^*\bigl(\Omega_\gamma,
             \muhom(\mathscr L_\gamma,\mathscr L_\eta)\bigr)\end{align}

 \begin{Lemma}           
The algebra $\bR$ acts on $M(\rho,\gamma)$
by composition in the second argument of $\muhom$. 
\end{Lemma}

\begin{proof}
By \ref{M2}, the restriction followed by the composition defines an associative map
\[
 \Ext^p_{\dG}(\mathscr L_\eta,\mathscr L_\nu)\otimes
 H_{\dG}^q(\Omega_\gamma,\muhom(\mathscr L_\gamma,\mathscr L_\eta))
 \longrightarrow
 H_{\dG}^{p+q}(\Omega_\gamma,\muhom(\mathscr L_\gamma,\mathscr L_\nu)).
\]
This defines the action. The idempotent $e_\eta$ projects onto the
$\eta$-summand. Finiteness of the support makes the action locally unital.
\end{proof}

We call the locally unital graded left $\bR$-module $M(\rho,\gamma)$ the \emph{microlocal module}.
Throughout this subsection we abbreviate
\begin{align}\label{abv1}
n=\dim\dN,\quad
 i_\gamma:Z^{\rm st}_\gamma\hookrightarrow\Omega_\gamma,\quad
 \dg=\operatorname{Lie}\dG,\quad
 d_\alpha=\dim\scrX_\alpha.
\end{align}
We keep the unshifted convention for $\varphi_f$ of
\S\ref{subsec:vc}, and put
\begin{align}\label{phiphi} \phi_f=\varphi_f[-1].\end{align}
Thus $\phi_f$ is the perverse normalized vanishing cycle functor.
For two types $\eta,$ $\gamma$ we put
\begin{align}\label{YW}
\begin{split}
Y_{\eta,\gamma} =\scrX_\eta\times\scrX_\gamma\times \dN^*
,\quad
 W_{\eta,\gamma}=\{(u,v,y)\in Y_{\eta,\gamma}:\ \pi_\eta(u)=\pi_\gamma(v)\},
   \end{split}
\end{align}
and
\begin{align}\label{FEG}
\begin{split}
  f_{\eta,\gamma}&: \, Y_{\eta,\gamma} \to\mathbb C ,\ 
 (u,v,y)\mapsto\langle y,\pi_\eta(u)-\pi_\gamma(v)\rangle,\\
  r_{\eta,\gamma}&:\,W_{\eta,\gamma}\longrightarrow T^*\dN,
  \quad ((F,x),(F',x),y)\longmapsto(x,y).
 \end{split}
\end{align}
The variety $Y_{\eta,\gamma}$ is smooth.
The map $r_{\eta,\gamma}$ is proper because $W_{\eta,\gamma}$ is closed in
$\Fl_\eta\times\Fl_\gamma\times T^*\dN$.
The critical scheme of $f_{\eta,\gamma}$ is taken in the smooth variety $Y_{\eta,\gamma}$. 
Identifying a point
$$((F,x),(F',x),y)\in W_{\eta,\gamma}$$ with $(F,F',x,y)$, the  equations of the critical set are
\[
\begin{aligned}
\operatorname{Crit}(f_{\eta,\gamma})=\bigl\{(F,F',x,y)\,:\,
 F\in\Fl_\eta,\ F'\in\Fl_\gamma\,,\,
 x\in\dN_F\cap\dN_{F'}\,,\,
 y\in\dN_F^\perp\cap\dN_{F'}^\perp\,,\,
 \mu(x,y)=0\bigr\}.
\end{aligned}
\]
The equations are scheme theoretic, as proved below.
There are natural closed inclusions
\[
\operatorname{Crit}(f_{\eta,\gamma})
\xrightarrow{j_{\eta,\gamma}}W_{\eta,\gamma}
\xrightarrow{k_{\eta,\gamma}} f_{\eta,\gamma}^{-1}(0).
\]
We consider the composed map
\begin{align}\label{QEG}q_{\eta,\gamma}=r_{\eta,\gamma}\circ j_{\eta,\gamma}:\operatorname{Crit}(f_{\eta,\gamma})\longrightarrow T^*\dN,
  \quad (F,F',x,y)\longmapsto(x,y).
\end{align}
The variety $Y_{\eta,\gamma}$ is smooth of dimension $d_\eta+d_\gamma+n$.
We define the shifted perverse sheaf
\begin{align}\label{PEG2}
\phi_{\eta,\gamma}=\left.(\phi_{f_{\eta,\gamma}}
\bbC_{Y_{\eta,\gamma}})\right|_{\operatorname{Crit}(f_{\eta,\gamma})}[2d_\eta].
\end{align}
We regard sheaves on a closed subset as sheaves on the ambient
space by pushforward when this causes no ambiguity. In the next
proposition, the set $\Lambda_\gamma$ and the fibre product carry the
scheme structures defined by the embedding
$$\Lambda_\gamma\subset {\rm Fl}_\gamma\times T^*\dN$$
Constructible sheaves use the underlying analytic spaces.

\begin{Lemma}
\label{lem:offdiag-muhom}
\hfill
\begin{enumerate}[label=$\mathrm{(\alph*)}$,leftmargin=8mm]
\item
There is a scheme isomorphism
$\operatorname{Crit}(f_{\eta,\gamma})\cong\Lambda_\eta\times_{T^*\dN}\Lambda_\gamma.$

\item
The map $q_{\eta,\gamma}$ is proper. 
\item
There is a $\dG$-equivariant isomorphism
 $\muhom(\mathscr L_\gamma,\mathscr L_\eta)
 \cong (q_{\eta,\gamma})_*\phi_{\eta,\gamma}.$
\end{enumerate}
\end{Lemma}

\begin{proof}
Differentiating $f_{\eta,\gamma}$ gives the equations
$\pi_\eta(u)=\pi_\gamma(v)$ and
$d\pi_\eta^*y=d\pi_\gamma^*y=0$.
Since $\operatorname{Im}d\pi_\alpha=\dN_F+\dg\cdot x$,
these are the defining ideals of the fibre product in (a).
The flag varieties are projective, so its projection is proper,
proving (b).

For (c), put $\pi_{\eta,\gamma}=\pi_\eta\times\pi_\gamma$ and
let $p_1,p_2$ be the projections from $\dN\times\dN$.
Write 
$$\calE=R\mathcal Hom(p_2^*\mathscr L_\gamma,p_1^!\mathscr L_\eta)\in D^b(\dN\times\dN).$$
Since $\pi_\gamma$ is proper, we have $\mathbb D\mathscr L_\gamma=(\pi_\gamma)_*\bbC$,
hence
\begin{align}\label{EPI}
 \mathcal E\cong\mathscr L_\eta\boxtimes\mathbb D\mathscr L_\gamma
 \cong(\pi_{\eta,\gamma})_*\bbC[2d_\eta].
\end{align}
Define
\begin{align}\label{DF}
\begin{split}
 \delta&:\dN\times\dN\to\dN,\quad(x_1,x_2)\mapsto x_1-x_2,\\
 h&:\dN\times\dN\times\dN^*\to\mathbb C,\quad (x_1,x_2,y)\mapsto\langle y,x_1-x_2\rangle.
 \end{split}
\end{align}
The diagonal $\Delta=\delta^{-1}(0)$ is smooth, with conormal
identified with $T^*\dN$ by the map 
$$(x,y)\mapsto(x,x;y,-y).$$
Write
\begin{align}\label{iota}
\begin{split}
 \iota&:T^*\dN\hookrightarrow h^{-1}(0),\quad(x,y)\mapsto(x,x,y),\\
 r&:\dN\times\dN\times\dN^*\to\dN\times\dN,\quad(x,_1,x_2,y)\mapsto (x_1,x_2).
 \end{split}
\end{align}
The pairing formula
\cite[Thm.~4.1, Def.~A.28, Rem.~A.30]{Sch25} gives
\begin{align}\label{MP}\mu_\Delta \calE\cong\iota^*\phi_h(r^*\calE)\end{align}
with the normalized functor \eqref{phiphi} and no further shift.
The proper map
\begin{align}\label{OP}
 P_{\eta,\gamma}=\pi_{\eta,\gamma}\times\mathrm{id}_{\dN^*}:
 Y_{\eta,\gamma}\longrightarrow\dN\times\dN\times\dN^*
\end{align}
satisfies $h\circ P_{\eta,\gamma}=f_{\eta,\gamma}$. Put
$$\varphi_{\eta,\gamma}=\phi_{f_{\eta,\gamma}}\bbC_{Y_{\eta,\gamma}}[2d_\eta]\in D^b(f_{\eta,\gamma}^{-1}(0))$$
Let $P_{\eta,\gamma,0}$ be the restriction of $P_{\eta,\gamma}$ to the central fiber $f_{\eta,\gamma}^{-1}(0)$.
By definition, we have
$$\muhom(\mathscr L_\gamma,\mathscr L_\eta)=\mu_\Delta\calE.$$
Consider the Cartesian square
\[
\begin{array}{ccc}
Y_{\eta,\gamma}
=\mathscr X_\eta\times\mathscr X_\gamma\times\underline N^*
&\xrightarrow{\ q\ }&
\mathscr X_\eta\times\mathscr X_\gamma
\\[2mm]
{\scriptstyle P_{\eta,\gamma}}\downarrow
&&\downarrow{\scriptstyle\pi_{\eta,\gamma}}
\\[2mm]
\underline N\times\underline N\times\underline N^*
&\xrightarrow{\ r\ }&
\underline N\times\underline N ,
\end{array}
\]
 The base change along the square above and \eqref{EPI} give the isomorphism
$$r^*\calE\cong (P_{\eta,\gamma})_*\bbC[2d_\eta].$$
Together with the isomorphism \eqref{MP}, this yields the first isomorphism in
\begin{equation}\label{isom0}
 \muhom(\mathscr L_\gamma,\mathscr L_\eta)
 \cong\iota^*\phi_h(P_{\eta,\gamma})_*\bbC[2d_\eta]
 \cong\iota^*(P_{\eta,\gamma,0})_*\varphi_{\eta,\gamma}.
\end{equation}
The second one follows by commutation of vanishing cycles with proper pushforward.
The following square is Cartesian
\[
\begin{array}{ccc}
 W_{\eta,\gamma} & \xrightarrow{k_{\eta,\gamma}} & f_{\eta,\gamma}^{-1}(0)\\[2mm]
 \Big\downarrow\, {\scriptstyle r_{\eta,\gamma}} && {\scriptstyle P_{\eta,\gamma,0}}\,\Big\downarrow\\[2mm]
 T^*\dN & \xrightarrow{\iota} & h^{-1}(0)
\end{array}
\]
Therefore
\begin{align}\label{isom1}
 \iota^*(P_{\eta,\gamma,0})_*\varphi_{\eta,\gamma}
 \cong(r_{\eta,\gamma})_*(k_{\eta,\gamma})^*\varphi_{\eta,\gamma}.
\end{align}
Away from the critical scheme the function $f_{\eta,\gamma}$ is a submersion
on the smooth space $Y_{\eta,\gamma}$, so its vanishing cycles vanish.
Thus
\begin{align}\label{isom2}
 (k_{\eta,\gamma})^*\varphi_{\eta,\gamma}
 \cong(j_{\eta,\gamma})_*\phi_{\eta,\gamma}.
\end{align}
Combining \eqref{isom0}-\eqref{isom2} proves (c). No smoothness of
$W_{\eta,\gamma}$ or of the critical scheme is required.
All constructions commute with smooth base change, hence descend
from the equivariant simplicial presentation.
\end{proof}

Recall the projection $p_\eta$ in \eqref{p}. For each type $\eta$ set
\begin{align}\label{LL}
\begin{split}
\Lambda_{\eta,\gamma}
 =\Lambda_\eta\times_{T^*\dN}Z_\gamma^{\rm st}
 =\{(F,\xi):\ \xi=(x,y)\in Z_\gamma^{\rm st},\ F\in\Fl_\eta,
       \ x\in\dN_F,\ y\in\dN_F^\perp\}.
\end{split}
\end{align}
The restriction of the map $p_\eta$ in \eqref{p} is the proper map
\begin{align}\label{PEG1}
p_{\eta,\gamma}:\Lambda_{\eta,\gamma}\longrightarrow Z_\gamma^{\rm st},
\quad (F,\xi)\longmapsto\xi.
\end{align}
Lemmas \ref{lem:offdiag-muhom}(a) and \ref{lem:Z} give a map
\[
s_{\eta,\gamma}:\Lambda_{\eta,\gamma}\longrightarrow
\operatorname{Crit}(f_{\eta,\gamma}),\quad
(F,\xi)\longmapsto(F,F_\xi,\xi).
\]
We define the constructible complex
\[
\mathscr L_{\eta,\gamma}=(s_{\eta,\gamma})^*\phi_{\eta,\gamma}
\quad\text{on }\Lambda_{\eta,\gamma}.
\]

\begin{Lemma}
\label{lem:microlocal-flag-module}\hfill
\begin{enumerate}[label=$\mathrm{(\alph*)}$,leftmargin=8mm]
\item
There are isomorphisms
\begin{align}\label{eq:offdiag-stable}
\begin{split}
 \muhom(\mathscr L_\gamma,\mathscr L_\eta)|_{\Omega_\gamma}
 \cong (i_\gamma)_*(p_{\eta,\gamma})_*\mathscr L_{\eta,\gamma},\quad
 e_\eta M(\rho,\gamma)\cong H_{\dG}^*(\Lambda_{\eta,\gamma},\mathscr L_{\eta,\gamma}).
\end{split}
\end{align}

\item
For each $\eta$ the $\eta$-piece $e_\eta M(\rho,\gamma)$ is finite dimensional, and
\[
 e_\eta M(\rho,\gamma)=0\quad\text{if }\eta\not\leqslant\gamma
                           \text{ or }\eta\in\Cyc(\rho),
\quad
 e_\gamma M(\rho,\gamma)\cong H^*_{P_\gamma}(\Lambda^{\rm st}_{F_\gamma},\bbC).
\]
The local Cartan algebra $H^*_{P_\gamma}\subset e_\gamma \bR e_\gamma$  acts via the restriction
$H^*_{P_\gamma}\to H^*_{P_\gamma}(\Lambda^{\rm st}_{F_\gamma},\bbC)$.

\item
For $\xi\in Z_\gamma^{\rm st}$ the stalk of the left side of
\eqref{eq:offdiag-stable} is
$R\Gamma\bigl(p_{\eta,\gamma}^{-1}(\xi),
 \mathscr L_{\eta,\gamma}|_{p_{\eta,\gamma}^{-1}(\xi)}\bigr).$
\end{enumerate}
\end{Lemma}

\begin{proof}
(a) Over $\Omega_\gamma$, Lemma~\ref{lem:Z} uniquely identifies
the second flag with $F_\xi$. Hence $s_{\eta,\gamma}$ gives
\begin{align}\label{OLL1}
 q_{\eta,\gamma}^{-1}(\Omega_\gamma)\cong\Lambda_{\eta,\gamma},
\end{align}
under which
\begin{align}\label{OLL2}
 q_{\eta,\gamma}|_{q_{\eta,\gamma}^{-1}(\Omega_\gamma)}
 =i_\gamma\circ p_{\eta,\gamma}.
\end{align}
Restrict Lemma~\ref{lem:offdiag-muhom}(c) and apply proper base
change, then equivariant hypercohomology. This proves (a).
Taking a stalk instead proves (c).

(b) By \ref{M1} and Lemma~\ref{lem:LZL}(b), the
microlocal Hom is supported on $Z_\eta\cap Z_\gamma^{\rm st}$.
The vanishing in (b) follows from Lemma~\ref{lem:Z}.
Propositions~\ref{prop:psi}(a),(c) and \ref{prop:cohZ} give the
diagonal and Cartan identifications. 
Finally, since the complex descends under the free action to the quasi-projective variety
$\mathfrak L_{\rho,\gamma}$, see Lemma~\ref{F}, each piece has finite dimensional hypercohomology. 
\end{proof}

The next proposition describes the $\mathbf R$-action on $M(\rho,\gamma)$ using cohomological correspondences.

\begin{Proposition}\label{prop:microlocal-convolution}\hfill
\begin{enumerate}[label=$\mathrm{(\alph*)}$,leftmargin=8mm]
\item
An element $a\in e_\nu\bR e_\eta$ acts on $M(\rho,\gamma)$
by convolution along the Steinberg correspondence $\scrX_\nu\times_{\dN}\scrX_\eta$
with coefficients in vanishing cycle complexes, yielding a map
$$H_{\dG}^*(\Lambda_{\eta,\gamma},\mathscr L_{\eta,\gamma})
 \rightarrow
 H_{\dG}^{*}(\Lambda_{\nu,\gamma},\mathscr L_{\nu,\gamma}).$$
The $\bR$-action respects Yoneda composition and identities.

\item
The $\bR$-action on $M(\rho,\gamma)$ factors through $\bT$ and extends to $\bT^\wedge$.
\item
There is a  right $H^*_{P_\gamma}(\Lambda^{\rm st}_{F_\gamma},\bbC)$-action on $M(\rho,\gamma)$ 
commuting with the left $\bR$-action.
\end{enumerate}
\end{Proposition}

\begin{proof}
(a) 
Choose a morphism $a:\mathscr L_\eta\to\mathscr L_\nu[k]$. 
By base change and adjunction, it yields a map
\[
 (s_\eta)^*\omega_{\scrX_\eta}
 \longrightarrow (s_\nu)^!\omega_{\scrX_\nu}[k]
 \quad\text{on }\scrZ_{\nu,\eta}=\scrX_\nu\times_{\dN}\scrX_\eta,
\]
where $s_\eta,$ $s_\nu$ are the proper projections \cite[\S8.6]{CG97}. 
Recall $Y_{\alpha,\gamma}$ from \eqref{YW}. Set
$$D=\scrZ_{\nu,\eta}\times\scrX_\gamma\times\dN^*
,\quad
c_\alpha:D\to Y_{\alpha,\gamma}.$$
Taking the external product with
$\bbC_{\scrX_\gamma\times\dN^*}$ gives a map
\begin{align}\label{mor2}
 (c_\eta)^*\calE_\eta\longrightarrow (c_\nu)^!\calE_\nu[k],
 \quad\calE_\alpha=\bbC_{Y_{\alpha,\gamma}}[2d_\alpha].
\end{align}
We have
$$f_{\eta,\gamma}c_\eta=f_{\nu,\gamma}c_\nu:D\longrightarrow\mathbb C.$$
Let $g$ denote both functions.
For \(\alpha=\eta,\nu\) we set
$$D_0=g^{-1}(0)
,\quad
Y_{\alpha,\gamma,0}=f_{\alpha,\gamma}^{-1}(0).$$ 
Since \(g=f_{\alpha,\gamma}\circ c_\alpha\), the map \(c_\alpha\) restricts to
\[
c_{\alpha,0}:D_0\longrightarrow Y_{\alpha,\gamma,0}.
\]
Thus, we have a Cartesian square
\[
\begin{array}{ccc}
D_0&\xrightarrow{c_{\alpha,0}}&
Y_{\alpha,\gamma,0}\\
\downarrow&&\downarrow\\
D&\xrightarrow{c_\eta}&Y_{\alpha,\gamma}.
\end{array}
\]
This yields the exchange morphism
\[
(c_{\alpha,0})^{*}\phi_{f_{\alpha,\gamma}}(\mathcal E_\eta)
\longrightarrow
\phi_g((c_\alpha)^{*}\mathcal E_\alpha).
\]
Consider the sheaves
$$\varphi_{\alpha,\gamma}=\phi_{f_{\alpha,\gamma}}\calE_\alpha.$$
The exchange change morphism, its Verdier dual, 
and the morphism \(\eqref{mor2}\) yield the diagram
\begin{align}\label{CC}
 \begin{split}
 (c_{\eta,0})^*\varphi_{\eta,\gamma}
 &\longrightarrow\phi_g((c_\eta)^*\calE_\eta)
 \longrightarrow\phi_g((c_\nu)^!\calE_\nu)[k]
 \longrightarrow (c_{\nu,0})^!\varphi_{\nu,\gamma}[k].
 \end{split}
\end{align}
Both zero fibres map properly to $h^{-1}(0)$ through $P_{\alpha,0}$,
with common composite $$q=P_{\eta,0}c_{\eta,0}=P_{\nu,0}c_{\nu,0}.$$
The unit for the adjunction $(c_{\eta,0})^*\dashv(c_{\eta,0})_*$,
the functor $q_*$ applied to \eqref{CC}, and the counit for
$(c_{\nu,0})_!\dashv (c_{\nu,0})^!$ yield a map
\[
 (P_{\eta,0})_*\varphi_{\eta,\gamma}
 \longrightarrow(P_{\nu,0})_*\varphi_{\nu,\gamma}[k].
\]
The commutation of proper pushforward with vanishing cycles identifies
this map with the image of $a$ under the pairing construction
\eqref{isom0}. Restricting by $\iota$ and to $\Omega_\gamma$, we recover the microlocal action. 
Cartan classes act by pullback along
the map retaining the $\eta$-flag.
The other generators use the
flag-changing correspondences above.

(c) The right composition with the isomorphism
$$\muhom(\mathscr L_\gamma,\mathscr L_\gamma)|_{\Omega_\gamma}\cong(i_\gamma)_*\bbC_{Z_\gamma^{\rm st}}$$
gives the right cup product action, which commutes with the left action by associativity.

(b) Cyclotomic pieces vanish by Lemma~\ref{lem:microlocal-flag-module},
so the action factors through $\bT$. The left action of
$\mathfrak m=H_{\dG}^{>0}$ equals right multiplication by its image
in $H_{\dG}^{>0}(Z_\gamma^{\rm st},\bbC)$. The latter ideal is
nilpotent, so one power of $\mathfrak m$ kills the whole module.
The action therefore extends to $\bT^\wedge$, componentwise.
\end{proof}

\subsection{Microlocal modules and universal Verma modules}
By Proposition~\ref{prop:psi}(b),(d) and Definition~\ref{def:localB},
the map $\Psi_\gamma$ in \eqref{PSG} descends continuously to
\begin{align}\label{OPS}
 \widehat\Psi_\gamma:\widehat\B_\gamma
 \longrightarrow H^*_{P_\gamma}(\Lambda^{\rm st}_{F_\gamma},\bbC).
\end{align}
The target is a $\widehat\B_\gamma$-module such that
$b\cdot h=\widehat\Psi_\gamma(b)\cup h$.
Recall the Verma module $\Delta(\rho,\gamma)$ from \eqref{BD}.

\begin{Proposition}\label{prop:microlocal-standard-map}
There is a morphism of graded
$(\bT^\wedge,H^*_{P_\gamma}(\Lambda^{\rm st}_{F_\gamma},\bbC))$-bimodules
\begin{equation}\label{eq:microlocal-standard-map}
 \Delta(\rho,\gamma)
       \otimes_{\widehat\B_\gamma}H^*_{P_\gamma}(\Lambda^{\rm st}_{F_\gamma},\bbC)
       \longrightarrow M(\rho,\gamma),
 \quad a\otimes h\longmapsto a\cdot h.
\end{equation}
The map is the identity on the $\gamma$-piece. Its image is
the submodule generated by $e_\gamma M(\rho,\gamma)$.
\end{Proposition}

\begin{proof}
The higher-type pieces vanish by Lemma~\ref{lem:microlocal-flag-module}.
The continuity therefore factors the action of $\bT^\wedge e_\gamma$
on $H^*_{P_\gamma}(\Lambda^{\rm st}_{F_\gamma},\bbC)$ through $\Delta(\rho,\gamma)$.
The associativity gives
$(ab)h=a(\widehat\Psi_\gamma(b)h)$ for $b\in\widehat\B_\gamma$.
Hence, the map is balanced and it is a bimodule morphism. Its $\gamma$-piece is the isomorphism
$$\widehat\B_\gamma\otimes_{\widehat\B_\gamma}H^*_{P_\gamma}(\Lambda^{\rm st}_{F_\gamma},\bbC)\cong 
H^*_{P_\gamma}(\Lambda^{\rm st}_{F_\gamma},\bbC)$$
Its image is the sum of the left translates of this summand.
\end{proof}

\Needspace{5\baselineskip}
\begin{Remark}
\hfill
\begin{enumerate}[label=$\mathrm{(\alph*)}$,leftmargin=8mm]

\item
The specialized microlocal module $M(\rho,\gamma)^\xi$ need not be
generated by its highest vector
$$1_\gamma\in e_\gamma M(\rho,\gamma)^\xi\cong\bbC.$$
This fails in the type $A_5$ configuration in \cite[\S7]{KS97}.
Two nonisomorphic simple perverse IC constituents have
a common conormal component, and a nonzero degree zero microlocal
morphism at a generic covector need not lift to a global morphism.
Indeed, degree zero Hom between the two simple perverse sheaves
vanishes, while positive perverse Ext degrees cannot contribute to
the degree zero microlocal stalk. Thus the distinguished class
$1_\gamma$ need not generate the specialization. 

\item
Let $\bbC_\xi$ be the $\widehat\B_\gamma$-module given  by the evaluation of $\widehat\Psi_\gamma$ at $\xi$.
We don't know if the pointwise analogue 
$\Delta(\rho,\gamma)\otimes_{\widehat\B_\gamma}\bbC_\xi
 \to M(\rho,\gamma)^\xi$
of the map \eqref{eq:microlocal-standard-map}
is surjective, see part (a). 
Note that the restriction map
\[
 M(\rho,\gamma)\otimes_{H^*_{P_\gamma}(\Lambda^{\rm st}_{F_\gamma},\bbC)}\bbC_\xi
 \longrightarrow M(\rho,\gamma)^\xi
\]
may not be an isomorphism because the right $H^*_{P_\gamma}(\Lambda^{\rm st}_{F_\gamma},\bbC)$-module $M(\rho,\gamma)$
may not be free. 

\end{enumerate}

\end{Remark}


\subsection{Specialized microlocal modules and Morse modules}
\label{subsec:microlocal-morse}
We compare restriction of microlocal Hom to a point with the
Morse construction of Lemma \ref{lem:Morse}.
Assume $Z_\gamma^{\rm st}\ne\emptyset$, equivalently assume \eqref{bfE} by Lemma \ref{lem:Z}.
Fix a closed point $\xi=(x,y)$ in $Z_\gamma^{\rm st}$.
For every $\eta$, the restriction to $\xi$ followed by forgetting
equivariance defines a homogeneous map
\[
H_{\dG}^*(\Omega_\gamma,\muhom(\mathscr L_\gamma,\mathscr L_\eta))
\longrightarrow H^*(\muhom(\mathscr L_\gamma,\mathscr L_\eta)_\xi).
\]
We define
\[
 M(\rho,\gamma)^\xi=
 \bigoplus_\eta H^*(\muhom(\mathscr L_\gamma,\mathscr L_\eta)_\xi).
\]
The composition morphisms \ref{M2} give $M(\rho,\gamma)^\xi$ a
locally unital left $\bR$-action called the \emph{specialized microlocal module} (at $\xi$).
By Proposition~\ref{prop:psi}(a), we have $\mathcal K_\gamma\cong
\bbC_{Z_\gamma^{\rm st}}$. Thus its one dimensional $\gamma$-piece
is $H^0((i_\xi)^*\mathcal K_\gamma)=\bbC$, with no additional shift.
An element of $e_\nu\bR^k e_\eta$ sends degree $j$ to degree $j+k$.

\begin{Lemma}
The restrictions for all $\eta$'s assemble to a graded $\bR$-module map
\begin{equation}\label{MM}
 \operatorname{res}_\xi:
 M(\rho,\gamma)\longrightarrow M(\rho,\gamma)^\xi.
\end{equation}
\end{Lemma}

\begin{proof}
Restriction commutes with the composition maps \ref{M2}, proving $\bR$-linearity.
The $\gamma$-piece is the evaluation $H^*_{P_\gamma}(\Lambda^{\rm st}_{F_\gamma},\bbC)\to\bbC$.
For the comparison with Morse modules, choose a Morse function $f$.
Recall the Morse fibre $MF^\xi_\eta$ in \eqref{MF}.
We have  $MF^\xi_\gamma\cong\bbC$.
Microlocal localization \cite[Ch.~VI]{KS90} identifies
\[
 H^k(\muhom(F,G)_\xi)
 =\Hom_{D^b(\dN;\xi)}(F,G[k]).
\]
The Morse functor annihilates objects whose microsupport avoids
$\xi$ and therefore factors through this localization; on our
objects it is $\mathscr L_\eta\mapsto MF^\xi_\eta$
\cite[\S8.6]{KS90}. It induces
\begin{equation}\label{eq:micro-Morse-morphism}
 H^k(\muhom(\mathscr L_\gamma,\mathscr L_\eta)_\xi)
 \longrightarrow\Hom_{D^b(\bbC)}(MF^\xi_\gamma,MF^\xi_\eta[k]),
\end{equation}
compatible with composition, hence with the $\bR$-actions.
This is a morphism, not an asserted isomorphism.
\end{proof}

The map \eqref{eq:micro-Morse-morphism} yields a pairing
$$M(\rho,\gamma)^\xi\otimes H^0(MF^\xi_\gamma)
\to\mathring M(\rho,\gamma)^\xi$$
with $H^0(MF^\xi_\gamma)\cong\mathbb C$. 
Choosing a generator of $H^0(MF^\xi_\gamma)$ we get a graded $\bR$-module homomorphism
\begin{equation}\label{eq:micro-to-Morse}
 \theta_{\xi,f}:M(\rho,\gamma)^\xi
 \longrightarrow\mathring M(\rho,\gamma)^\xi.
\end{equation}
In the type $A_1$ calculation of
\S\ref{sec:micro-to-Morse}
the map $\theta_{\xi,f}$ is neither injective nor surjective.
Combining \eqref{eq:micro-to-Morse} with \eqref{MM} gives a graded $\bR$-module map
\begin{equation}\label{eq:global-micro-to-Morse}
 M(\rho,\gamma)\xrightarrow{\operatorname{res}_\xi}
 M(\rho,\gamma)^\xi\xrightarrow{\theta_{\xi,f}}
 \mathring M(\rho,\gamma)^\xi.
\end{equation}
The Morse module has highest piece of $\ell$-weight $Z_{\rho,\gamma}$
and finite length in finite type by Theorem~\ref{thm:morse-JH}. 
The microlocal module requires neither a point nor a Morse function.

\begin{Remark}\label{rem:morse-evaluation-dual}\hfill
\begin{enumerate}[label=$\mathrm{(\alph*)}$,leftmargin=8mm]

\item
The restriction \eqref{MM} is surjective
in the multiplicity free case, i.e., when
$\dim\dV_c\leqslant1$ for every vertex $c$ of the unfolded quiver
$\dQ$ and $Z_\gamma^{\rm st}\ne\emptyset$.
Then the $\mathbb Z$-flag varieties $\Fl_\eta$ are points. Put
\[
E=\dN,\quad U=\dN_{F_\eta},\quad V=\dN_{F_\gamma}.
\]
Then $\scrX_\eta=U$, $\scrX_\gamma=V$ and the pairing function
of \eqref{FEG} is
\[
f_{\eta,\gamma}:U\times V\times E^*\longrightarrow\bbC,
\quad (u,v,y)\longmapsto\langle y,u-v\rangle.
\]
Its critical scheme is the vector space
\[
\operatorname{Crit}(f_{\eta,\gamma})=(U\cap V)\times(U+V)^\perp,
\]
and its normal Hessian is globally hyperbolic. The normalization
of \eqref{PEG2} therefore gives
\[
\mathscr L_{\eta,\gamma}\cong
\bbC_{\Lambda_{\eta,\gamma}}
 [2(\dim(U\cap V)-\dim V)].
\]
Since the flag of type $\eta$ is unique, the projection
$p_{\eta,\gamma}$ of \eqref{PEG1} is a closed immersion.
If $\xi$ is outside its image, the corresponding microlocal Hom
stalk is zero. If $\xi$ is in its image, restriction to the stalk
is, up to the displayed common shift,
\[
H_{\dG}^*(\Lambda_{\eta,\gamma},\bbC)
\longrightarrow H^*(\{\xi\},\bbC)=\bbC,
\]
which sends $1$ to $1$. Hence \eqref{MM} is surjective on every
pieces.

\item
Fix the point $\xi$ and the Morse function $f$ in \eqref{fH}.
Recall the Morse fibre $MF^\xi_\eta$ considered in \eqref{MF}.
Its cohomology is finite dimensional. We have
\[
\mathring M(\rho,\gamma)^\xi
 =\bigoplus_\eta H^*(MF^\xi_\eta),\quad
\]
 Put
$\mathcal P_\eta=\mathscr L_\eta[-d_\eta]$
and $\phi_f=\varphi_f[-1]$ as in \eqref{phiphi}. Then
\[
MF^\xi_\eta
 \cong i_x^*\phi_f(\mathcal P_\eta)[d_\eta-d_\gamma]
\]
We define the \emph{co-Morse module} to be the left $\bR$-module given by
\[
 \mathring M^{!}(\rho,\gamma)^\xi
 =\bigoplus_\eta
 H^*\bigl(i_x^!\phi_f(\mathcal P_\eta)
                  [d_\eta-d_\gamma]\bigr).
\]
The shift $[d_\eta-d_\gamma]$ is required to keep
this grading equal to that of the Morse fibres $MF^\xi_\eta$.
In particular, the $\gamma$-piece is in degree zero.
The costalk-to-stalk morphisms
\[
 i_x^!\phi_f(\mathcal P_\eta)
 \longrightarrow i_x^*\phi_f(\mathcal P_\eta)
\]
give a graded left $\bR$-module homomorphism
\begin{equation}\label{eq:coMorse-to-Morse}
 \mathring{\mathcal I}_\xi:\mathring M^!(\rho,\gamma)^\xi
       \longrightarrow \mathring M(\rho,\gamma)^\xi.
\end{equation}
The image of
\eqref{eq:coMorse-to-Morse} need not be simple.
The rank one calculation in
\S\ref{ex:coMorse-A1-nonsimple} shows that this
map can be an isomorphism onto a nonsimple Morse module.

\item
Let $i_\xi:\{\xi\}\hookrightarrow Z_\gamma^{\rm st}$ be the embedding of a closed point.  Write
\begin{align}\label{eq:microlocal-K}\mathcal K_\eta=\muhom(\mathscr L_\gamma,\mathscr L_\eta)|_{Z_\gamma^{\rm st}}.\end{align}
The specialized microlocal module defined above has pieces
$$e_\eta M(\rho,\gamma)^\xi=H^*((i_\xi)^*\mathcal K_\eta).$$
The \emph{co-specialized microlocal module} is
\[
 M^!(\rho,\gamma)^\xi
 =\bigoplus_\eta H^*\bigl((i_\xi)^!\mathcal K_\eta[2n]\bigr),
 \quad n=\dim\dN=\dim Z_\gamma^{\rm st}.
\]
The shift $[2n]$ puts the $\gamma$-piece
$H^0((i_\xi)^!\bbC_{Z_\gamma^{\rm st}}[2n])\cong\bbC$
in degree zero, because $Z_\gamma^{\rm st}$ is smooth of
complex dimension $n$. As $[2n]$ is independent of $\eta$,
the microlocal composition gives a graded, locally unital left
$\bR$-module with the original Ext degrees. Shifting the natural
costalk-to-stalk morphism by $[2n]$ gives a graded $\bR$-module map
\begin{equation}\label{eq:co-specialized-to-specialized}
M^!(\rho,\gamma)^\xi
 \longrightarrow
 \bigoplus_\eta H^*\bigl((i_\xi)^*\mathcal K_\eta[2n]\bigr)
 =M(\rho,\gamma)^\xi[2n].
\end{equation}
This map can vanish on the $\gamma$-component, because Proposition \ref{prop:psi}(a) yields 
$\mathcal K_\gamma\cong\bbC_{Z_\gamma^{\rm st}}$ and
\[
 (i_\xi)^!\mathcal K_\gamma[2n]\cong\bbC,
 \quad (i_\xi)^*\mathcal K_\gamma[2n]\cong\bbC[2n].
\]
Therefore, we modify this map in the following way.
Recall the quotient map in \eqref{QEG}
$$q:Z_\gamma^{\rm st}\to\mathfrak L_{\rho,\gamma}$$ 
Put $\bar\xi=q(\xi)$, and let $i_{\bar\xi}:\{\bar\xi\}\hookrightarrow
\mathfrak L_{\rho,\gamma}$ be the inclusion. 
The sheaf $\mathcal K_\eta$ descends to a sheaf $\mathcal F_\eta$ on $\mathfrak L_{\rho,\gamma}$.
A local trivialization of the map $q$, with
$g=\dim\dG$ and $d=n-g$, gives
\[
 (i_\xi)^!\mathcal K_\eta[2n]
 \cong (i_{\bar\xi})^!\mathcal F_\eta[2d].
\]
Thus the co-specialized module has the following realization
$$M^!(\rho,\gamma)^\xi\cong
\bigoplus_\eta H^*(i_{\bar\xi}^!\mathcal F_\eta[2d]).$$
We define the \emph{sliced} costalk-to-stalk morphism to be the graded $\bR$-module map
\begin{equation}\label{eq:co-specialized-slice-map}
 \mathcal I_\xi:M^!(\rho,\gamma)^\xi
 \longrightarrow
 \bigoplus_\eta H^*((i_{\bar\xi})^*\mathcal F_\eta[2d])
 =M(\rho,\gamma)^\xi[2d].
\end{equation}
 In the rank one calculation of
\S\ref{subsubsec:A1-slice-image} the image is  simple at
$\bar\xi=0$. We do not know whether every nonzero image is simple.
\end{enumerate}
\end{Remark}

\subsection{Jordan-H\"older multiplicities of specialized microlocal modules}
\label{subsec:microlocal-JH}
Let $(C,D)$ be of finite type. 
Let $\Gamma$ be a weight block, and $\mathscr C$ the corresponding finite length subcategory of
$\scrO_\mu^\rho$ as in \eqref{CCG}.
Let $K(\mathscr C)$ be the complexified ungraded Grothendieck group of $\mathscr C$.
To the finite set $\Gamma$ correspond the corner $\bR_\Gamma$ and the simple $\bR_\Gamma$-modules
$E_{\Gamma,\lambda}$ of \S\ref{sec:JHMM}.
Put
\begin{align}\label{GST}\Gamma^{\rm st}=\{\gamma\in\Gamma:Z_\gamma^{\rm st}\ne\emptyset\}.\end{align}

\begin{Theorem}\label{thm:microlocal-JH}Let $\gamma\in\Gamma$.\hfill
\begin{enumerate}[label=$\mathrm{(\alph*)}$,leftmargin=8mm]
\item\label{item:microlocal-JH-length}
The module $M(\rho,\gamma)^\xi$ has finite length in $\scrO_\mu^\rho$, and
$[M(\rho,\gamma)^\xi:L(\rho,\gamma)]=1$.

\item\label{item:microlocal-JH-global}
For each $\eta\in\Gamma$ with
$[M(\rho,\gamma)^\xi:L(\rho,\eta)]\ne0$ there is a
unique $\lambda(\eta)\in\mathcal I_\Gamma$ with
$$e_\Gamma L(\rho,\eta)\cong E_{\Gamma,\lambda(\eta)}$$ as an $\bR_\Gamma$-module.
The type $\eta$ is the unique maximal charge type
$\nu\in\Gamma$ with $E_{\nu,\lambda(\eta)}\ne0$, and
\begin{equation}\label{eq:microlocal-JH-global}
 [M(\rho,\gamma)^\xi:L(\rho,\eta)]
 =\sum_k\dim H^k\bigl(\muhom(\mathscr L_\gamma,
             \IC_{\lambda(\eta)})_\xi\bigr).
\end{equation}

\item\label{item:microlocal-JH-CC}
Put $\Lambda=\overline{Z_\gamma^{\rm st}}\subset T^*\dN$.
For $\xi\in Z_\gamma^{\rm st}$  generic and $\eta$ as in (b) we have
\begin{equation}\label{eq:microlocal-JH-CC}
 [M(\rho,\gamma)^\xi:L(\rho,\eta)]=m_\Lambda(\IC_{\lambda(\eta)}),
 \quad \operatorname{CC}(\IC_{\lambda(\eta)})
   =\sum_{\Lambda'}m_{\Lambda'}(\IC_{\lambda(\eta)})[\Lambda'].
\end{equation}
\end{enumerate}
\end{Theorem}

\begin{proof}\hfill

(a) Lemma~\ref{lem:microlocal-flag-module} gives finite dimensional
pieces supported in $\eta\leqslant\gamma$ outside $\Cyc(\rho)$, such that $e_\gamma M(\rho,\gamma)^\xi=\bbC$.
The equivariant parameters act through the augmentation,
so the action extends to $\bT^\wedge$. The proof of
Theorem~\ref{thm:morse-JH}(a) applies with
$D(N)=\sum_{\eta\in\Gamma}\dim e_\eta N$.
Every nonzero subquotient has $D(N)>0$.
Hence the length is $\leqslant D(M(\rho,\gamma)^\xi)$.
The same highest piece argument yields the multiplicity one.

(b) Apply the functor $\muhom(\mathscr L_\gamma,-)_\xi$ to
\eqref{eq:morse-IC-decomposition} with $J=\Gamma$. We get
\[
 e_\Gamma M(\rho,\gamma)^\xi\cong\bigoplus_{\lambda,s}E_{\Gamma,\lambda,s}
 \otimes H^*\bigl(\muhom(\mathscr L_\gamma,\IC_\lambda[-s])_\xi\bigr).
\]
As in the finite corner argument of Theorem~\ref{thm:morse-JH},
the kernel $\ker (q_\Gamma)$ raises the microlocal degree and acts
nilpotently. Projecting onto a copy of $\IC_\lambda[-s]$, we get
\begin{equation}\label{eq:microlocal-JH-corner}
 [e_\Gamma M(\rho,\gamma)^\xi:E_{\Gamma,\lambda}]
 =\sum_k\dim H^k\bigl(\muhom(\mathscr L_\gamma,\IC_\lambda)_\xi\bigr).
\end{equation}
The restrictions $e_\Gamma L(\rho,\eta)$ are nonzero, simple and
pairwise nonisomorphic, since they are distinguished by their highest types.
The exactness of $e_\Gamma$ and the formuals \eqref{eq:microlocal-JH-corner} and
\eqref{eq:morse-simple-eta-pieces} give (b).

(c) 
Recall that $\Lambda$ is irreducible by Lemma~\ref{lem:stable-irreducible}.
The sheaf $\mathscr L_\gamma$ is microlocally simple along
$Z_\gamma^{\rm st}$ by Proposition~\ref{prop:psi}(a).
At a generic point of $Z_\gamma^{\rm st}$ the total dimension of
$H^*(\muhom(\mathscr L_\gamma,\IC_\lambda)_\xi)$ is
$m_\Lambda(\IC_\lambda)$ by the microlocal interpretation
of characteristic cycles \cite[Ch.~IX]{KS90}. Then, apply (b).
\end{proof}

\begin{Corollary}\label{cor:microlocal-morse-K0}
Choose a generic point $\xi\in Z_\gamma^{\rm st}$.
In the {ungraded} Grothendieck group we have
\begin{equation}\label{eq:microlocal-morse-K0-global}
 [M(\rho,\gamma)^\xi]
  =[\mathring M(\rho,\gamma)^\xi]
\end{equation}
\end{Corollary}

\begin{proof}
Apply Theorems \ref{thm:morse-JH} and \ref{thm:microlocal-JH}.
\end{proof}

\begin{Proposition}\label{prop:microlocal-basis}
Let $(C,D)$ be of finite type. 
\begin{enumerate}[label=$\mathrm{(\alph*)}$,leftmargin=8mm]
\item $\Gamma^{\rm st}=\Gamma$.
\item For each $\gamma\in\Gamma$ choose a generic point
$\xi_\gamma\in Z_\gamma^{\rm st}$. The families
$\bigl([\mathring M(\rho,\gamma)^{\xi_\gamma}]\bigr)_{\gamma\in\Gamma}$
and
$\bigl([M(\rho,\gamma)^{\xi_\gamma}]\bigr)_{\gamma\in\Gamma}$
coincide and form a basis of  $K(\mathscr C)$.
It  is independent of the choices of the generic points.
\end{enumerate}
\end{Proposition}

\begin{proof}
(a) 
Use the auxiliary  graph $\widehat H$ introduced in the proof of Lemma \ref{lem:stable-irreducible}.
The oriented representation space on
$$E=\dV\oplus\bigoplus_{c,L}W_c(L)$$ 
 is identified with $\dN$.
After adjoining the framing group $\prod_{c,L}GL(W_c(L))$ to $\dG$,
the flag-incidence complex $\mathscr L_\gamma$ is, up to shift, the Lusztig flag complex  of $\gamma$, 
each framing vertex occuring once with its full dimension.

Choose $\gamma\in\Gamma$.
By Theorem \(\ref{thm:V}\), we have \(e_\gamma L(\rho,\gamma)\neq0\). Choose a decomposition of \(\mathscr L_\gamma\) into shifted simple 
perverse sheaves. Since the corresponding orthogonal projectors sum to \(e_\gamma\), one of them, say \(p\), acts nontrivially on 
\(e_\gamma L(\rho,\gamma)\). Let \(\mathcal P[s]\) be its summand and $b_{\mathcal P}\in B_{\widehat H}(\infty)$ 
be the class of $\mathcal P$ in the crystal.

Fix the sequence \(\mathbf i=(i_1,\ldots,i_N)\) determined by \((\rho,\gamma)\) through the decreasing-height removal order as in the proof of
Lemma~\(\ref{lem:stable-irreducible}\). 
Consider the chain of maximal crystal removals
\begin{align}\label{IB}
 b_0=b_{\mathcal P},\quad
 a_k=\varepsilon_{i_k}(b_{k-1}),\quad
 b_k=\widetilde e_{i_k}^{\,a_k}b_{k-1}
, \quad \forall\,1\leqslant k\leqslant N.
\end{align}
We shall prove below that $b_N$ is the highest weight element $u_\infty$ in $B_{\widehat{H}}(\infty)$.
Let $\mathfrak Z_{b_{\mathcal P}}$ be the corresponding Kashiwara-Saito irreducible component. By
\cite[Thm.~6.2.2]{KS97} and Lemma~\ref{lem:LZL}, we have
\[
 \mathfrak Z_{b_{\mathcal P}}\subseteq\operatorname{SS}(\mathcal P)
 \subseteq\operatorname{SS}(\mathscr L_\gamma)\subseteq Z_\gamma.
\]
Choose an element $\xi\in\mathfrak Z_{b_{\mathcal P}}$ which is generic for the successive
maximal/top correspondences along the sequence $\mathbf i$. 
The dominance condition of these correspondences in 
\cite[Lem.~5.2.3, Prop.~5.2.4, Thm.~5.3.2]{KS97} allow to have
finitely many genericity conditions simultaneously.
Since $\xi\in Z_\gamma$, we can choose a
$\bbZ$-flag $F\in\Fl_\gamma$ such that $(F,\xi)\in\Lambda_\gamma$.

Using this sequence and the $\mathbb Z$-flag $F$, we construct the finite
descending $\widehat H$-vertex-graded flag
$E^\bullet=(E^{(k)})_{k=0}^N$ by the removal procedure of
Lemma~\ref{lem:stable-irreducible}, by replacing the $\mathbb Z$-flag $F_\gamma$ there by our $F$. 
At a step $(c,P)$ its quotient
has dimension
\[
 \dim(E^{(k-1)}/E^{(k)})
 =\dim(F_c^{\leqslant P}/F_c^{<P})=v_c(\gamma,P),
\]
the multiplicity of $(c,P)$ in $\gamma$. At a framing step
$r_{c,L}$ the quotient has dimension $\dim W_c(L)$.
Thus the labelled colour quotients recover $\gamma$.
Note also that,
since all colour positions in the chosen interval occur in the removal sequence, 
and each framing vertex is removed once, we have $$E^{(N)}=0$$

Next, we define another  finite descending
$\widehat H$-vertex-graded flag $U^\bullet=(U^{(k)})_{k=0}^N$ on $E$ by
\[
 U^{(0)}=E,\quad
 U^{(k)}=K_{i_k}\bigl(\xi|_{U^{(k-1)}}\bigr).
\]
Here $K_i$ retains the incoming images at $i$ and all other vertex spaces, see \eqref{K1}. 
An induction gives 
\begin{align}\label{UE}U^{(k)}\subseteq E^{(k)},\quad\forall k.\end{align}
Indeed, at step $k$,
all incoming images from $E^{(k-1)}$ lie in $E^{(k)}_{i_k}$ because
$E^{(k)}$ is a subrepresentation and $E^{(k-1)}/E^{(k)}$ is supported at $i_k$.  At a vertex $r_{c,L}$, the space is unchanged
until its unique occurrence and it is zero afterwards. Its removed quotient is therefore  $W_c(L)$.
In particular,  $$U^{(N)}=0.$$
Since \(U^{(0)}=E\) and \(U^{(N)}=0\), the sum of the dimensions of the quotients at all steps \((c,P)\) is
\[
\sum_P v_c(\eta,P)=\dim\underline V_c
=\sum_P v_c(\gamma,P).
\]
Thus \(\eta\in C_\gamma\).
Since $\eta$ and $\gamma$ are in the same weight block, it makes sense to compare them for the partial order \eqref{PO}.
The inclusions \eqref{UE} imply that
$$\gamma\leqslant\eta$$

By the generic maximal-top correspondences of Kashiwara-Saito
\cite[Lem.~5.2.3, Prop.~5.2.4, Thm.~5.3.2]{KS97}, the successive
representations $\xi|_{U^{(k)}}$ belong generically to the components
labelled by the $b_k$ of \eqref{IB}, and
\[
 a_k=\dim(U^{(k-1)}/U^{(k)}).
\]
Since $U^{(N)}=0$, this also proves $$b_N=u_\infty.$$

Lusztig's induction lemma \cite[Lem.~6.4]{Lus91}, extended to all
orientations by Fourier transform \cite[\S5]{Lus91}, places the simple perverse sheaf associated to the element
$b_{k-1}\in B_{\widehat H}(\infty)$, up to shift, in the induction of that of $b_k$ by an
$i_k$-supported quotient of dimension $a_k$. Iterating from
$b_N=u_\infty$ puts $\mathcal P$ in $\mathscr L_\eta$, up to shift.
Thus $\mathcal P$ is a simple shifted summand of both
\(\mathscr L_\gamma\) and \(\mathscr L_\eta\).

Now, since $\mathcal P$ a shifted simple summand  of both $\mathscr L_\gamma$
and $\mathscr L_\eta$, the inclusions and projections
give homogeneous elements $a,$ $b$ such that the projector to $\mathcal P$ factorizes in the following way
\[
 a\in e_\eta\bR e_\gamma,\quad
 b\in e_\gamma\bR e_\eta,\quad ba=p.
\]
Hence, since $p L(\rho,\gamma)\ne0$, we also have  
$$e_\eta L(\rho,\gamma)\ne0.$$

Further, by Theorem~\ref{thm:V}, we have 
$$e_\gamma L(\rho,\gamma)\ne0
,\quad
e_\nu L(\rho,\gamma)=0
,\quad
\forall\nu\in\Up(\gamma).$$
If $\gamma\ne\eta$, then
$\gamma<\eta$ gives $\eta\in\Up(\gamma)$ by Remark~\ref{rem:Up},
contrary to this nonvanishing. Hence 
$$\gamma=\eta.$$
As $U^{(k)}\subseteq E^{(k)}$ and the labelled colour quotients
of both flags have the same multiplicities ($\eta=\gamma$), we obtain
$U^{(k)}=E^{(k)}$ for every $k$; the framing quotients agree as well.
Thus every removal in $E^\bullet$ is maximal.
After a graded change of basis taking $F$ to $F_\gamma$, the
maximal removal criterion proved in Lemma~\ref{lem:stable-irreducible}
applies to $E^\bullet$ and gives $F_\xi=F$.
Thus $\xi\in Z_\gamma^{\rm st}$. 
Hence $\gamma\in\Gamma^{\rm st}$.

(b) Let $N_\gamma=\mathring M(\rho,\gamma)^{\xi_\gamma}$ or $M(\rho,\gamma)^{\xi_\gamma}$.
We have $e_\eta N_\gamma=0$ unless $\eta\leqslant\gamma$.
For Morse modules this is Lemma~\ref{lem:Morse}(c), for specialized
microlocal modules this follows from
\cite[Cor.~5.4.10]{KS90}, Lemma~\ref{lem:Z}(f) and \eqref{eq:microlocal-K}.
Since $e_\eta L(\rho,\eta)\neq0$,
no constituent $L(\rho,\eta)$ can occur unless $\eta\leqslant\gamma$.
Theorems~\ref{thm:morse-JH}(a), \ref{thm:microlocal-JH}(a)
then give
\[
 [N_\gamma]=[L(\rho,\gamma)]+
 \sum_{\substack{\eta\in\Gamma\\\eta<\gamma}}
 n_{\eta\gamma}[L(\rho,\eta)],\quad n_{\eta\gamma}\in\mathbb N.
\]
Thus both transition matrices are unitriangular and give integral
bases. Corollary~\ref{cor:microlocal-morse-K0} identifies them.
The independence of the generic points follows from the
characteristic cycle formulas in Theorems~\ref{thm:morse-JH}(c),
\ref{thm:microlocal-JH}(c) and the irreducibility of
$Z_\gamma^{\rm st}$ in Lemma~\ref{lem:stable-irreducible}.
\end{proof}

\begin{Remark}\label{rem:genericity-necessary}\hfill
\begin{enumerate}[label=$\mathrm{(\alph*)}$,leftmargin=8mm]
\item
The genericity assumption in Corollary~\ref{cor:microlocal-morse-K0}
cannot be omitted, see the example in \S\ref{ex:A1-K0-special}.
Smoothness and irreducibility of $Z_\gamma^{\rm st}$ do not imply
constancy of the microlocal or Morse multiplicities.

\item
For $\eta\in\Gamma$, let $\lambda(\eta)\in\mathcal I_\Gamma$
be the IC label  in Theorem~\ref{thm:morse-JH}(b). 
Theorem \cite[Thm.~4.4.2]{VV25} yields a
decategorification isomorphism $\DEC$ from the (complexified) ungraded Grothendieck group
of $\scrO_\mu^\rho$ to a weight subspace in a finite dimensional representation
of the simple Lie algebra of type $\dQ$.
We put 
$$b_\eta^\vee=\DEC([L(\rho,\eta)]).$$
This is a dual canonical basis element.
Choose $\xi$ generic in $Z_\gamma^{\rm st}$, and put
$\Lambda_\gamma=\overline{Z_\gamma^{\rm st}}$.
The closure is taken in $T^*\dN$.
Theorems~\ref{thm:morse-JH}(c) and \ref{thm:microlocal-JH}(c) give
\[
 \DEC([\mathring M(\rho,\gamma)^\xi])
 =\DEC([M(\rho,\gamma)^\xi])
 =\sum_{\eta\in\Gamma}m_{\Lambda_\gamma}(\IC_{\lambda(\eta)})b_\eta^\vee.
\]
Let $V_\rho^\vee$ be the span of the $b_\eta^\vee$ and
$V_\rho$ its dual, with basis $(b_\eta)$. Let $F_\rho$ be the vector
space with basis $([\Lambda])$, indexed by the irreducible components
of the microsupports of the sheaves $\IC_{\lambda(\eta)}$ with $\eta\in\Gamma$.
Consider the linear map
\[
 C_\rho:V_\rho\longrightarrow F_\rho,\quad
 b_\eta\longmapsto\operatorname{CC}(\IC_{\lambda(\eta)}).
\]
Let $(\epsilon_\Lambda)$ be the basis dual of $([\Lambda])$.
By Corollary~\ref{cor:microlocal-morse-K0}, the modules
$\mathring M(\rho,\gamma)^\xi$ and 
$M(\rho,\gamma)^\xi$ have equal classes in the ungraded Grothendieck
group, although they need not be isomorphic. The preceding formula
identifies the image of this common class under $\DEC$
with the image 
$$C_\rho^\vee(\epsilon_{\Lambda_\gamma})$$ 
of $\epsilon_{\Lambda_\gamma}$ by the transpose map
$C_\rho^\vee:F_\rho^\vee\to V_\rho^\vee.$
By Proposition~\ref{prop:microlocal-basis}, the matrix
$(m_{\Lambda_\gamma}(\IC_{\lambda(\eta)}))_{\gamma,\eta\in\Gamma}$
is unitriangular for a linear extension of $\leqslant$.
Since its diagonal coefficients are one, the class
$$\DEC([\mathring M(\rho,\gamma)^\xi])=\DEC([M(\rho,\gamma)^\xi])$$ belongs to the
dual canonical basis if and only if
\[
 m_{\Lambda_\gamma}(\IC_{\lambda(\eta)})=\delta_{\eta,\gamma}
 \quad\text{for every }\eta\in\Gamma.
\]

\end{enumerate}
\end{Remark}

\section{Proof of Theorem $\ref{thm:main}$}\label{sec:surj}
We prove Theorem $\ref{thm:main}$ for quivers of type ADE. 
The arguments used in the proof apply more generally to any finite quiver \(\Omega\), including loops and multiple arrows, 
with arbitrary integer gradings. We first establish these general geometric statements, 
then apply them with \(\Omega=\underline Q\) or \(\Omega=Q\).
The extension to more general symmetric quivers is treated separately in \S\ref{sec:general}.

\subsection{Quiver variety and their graded fixed loci}\label{subsec:ungraded}
The graded quiver varieties are fixed components of quiver variety.
Let $\Omega$ be any finite quiver, including loops and multiple arrows.
Let $\bar\Omega$ be the double of $\Omega$, and let
$V=\bigoplus_cV_c$ and $W=\bigoplus_cW_c$ be vector spaces indexed by
the vertices of $\Omega$. Put
\[
 v=(\dim V_c)_c,\quad
 w=(\dim W_c)_c.
\]
The Nakajima quiver variety is the quotient
\[
 \M_\Omega(v,w)
 =\mu^{-1}(0)^{\rm st}/\dG,
 \quad
\mu=\sum_{h\in\bar\Omega}
           \varepsilon(h)B_{\bar h}B_h+AA^*,
\]
where
\[
 B_h:V_c\to V_{c'},\quad
 A:W\to V,\quad A^*:V\to W
 \quad h:c\to c'
\]
Stability means that $V$ is generated by $A(W)$ under the
maps $B_h$. These maps are not required to be homogeneous.  
Fix an arbitrary integer grading $W_c=\bigoplus_LW_c(L)$, and let
$T_\rho=\bbC^\times$ act by
\[
 t\cdot(B,A,A^*)
 =\big(tB,\ A\rho(t)^{-1},\ t^2\rho(t)A^*\big),
 \quad \rho(t)|_{W_c(L)}=t^L.
\]
 It commutes with $\dG$, multiplies
$\mu$ by $t^2$, and preserves stability, so it acts on
$\M_\Omega(v,w)$. 
For a grading $V_c=\bigoplus_PV_c(P)$ given by a type $\gamma$ we put
$$v_\gamma=(\dim V_c(P))_{c,P},\quad w_\rho=(\dim W_c(L))_{c,L}.$$
Recall that $\M_{\dQ}^\bullet(v_\gamma,w_\rho)$ is the
quotient of the stable zero-moment data satisfying
\[
 B_h:V_c(P)\to V_{c'}(P-1),\quad
 A_c:W_c(L)\to V_c(L),\quad
 A^*_c:V_c(P)\to W_c(P-2)
\]
by $\dG_{v_\gamma}=\prod_{c,P}GL(V_c(P))$. The bullet distinguishes the
graded variety from its ungraded ambient variety.  
The following is well-known.

\begin{Lemma}\label{lem:fixed}
We have 
$$\M_\Omega(v,w)^{T_\rho}\cong\bigsqcup_\gamma\M_\Omega^\bullet(v_\gamma,w_\rho).$$
The sum runs over  gradings $v_\gamma$ with $\sum_Pv_c(P)=\dim V_c.$ 
The group $T_\rho$ acts on the fibre of the
tautological bundle $\calV_c(P)$ with weight $P$. 
\qed
\end{Lemma}

\begin{Lemma}\label{lem:kirwan}
For each grading $v_\gamma$ the local Kirwan map $\kappa_\gamma$ in \eqref{LKM1} is surjective. 
\end{Lemma}

\begin{proof}
By \cite{MN18} the global Kirwan map $H^*_{\dG}\to H^*(\M_\Omega(v,w),\bbC)$ is surjective. 
Since $H^*_{\dG}$ is concentrated in even degrees we have $H^{\rm odd}(\M_\Omega(v,w),\bbC)=0$.
By parity $H^*_{T_\rho}(\M_\Omega(v,w),\bbC)$ is free over $H^*_{T_\rho}=\bbC[a]$ and
$$H^*_{T_\rho}(\M_\Omega(v,w),\bbC)/(a)=H^*(\M_\Omega(v,w),\bbC).$$
The image $I$ of the equivariant global Kirwan map 
$$H^*_{\dG\times T_\rho}\to H^*_{T_\rho}(\M_\Omega(v,w),\bbC)$$ is a graded $\bbC[a]$-submodule with 
$$I+aH^*_{T_\rho}(\M_\Omega(v,w),\bbC)=H^*_{T_\rho}(\M_\Omega(v,w),\bbC).$$
Since
$H^*_{T_\rho}(\M_\Omega(v,w),\bbC)$ is bounded below, the graded Nakayama lemma gives 
\begin{align}\label{I=}I=H^*_{T_\rho}(\M_\Omega(v,w),\bbC).\end{align}
The localization theorem \cite{Qui71}, Lemma \ref{lem:fixed} and \eqref{I=} yield a surjection by right exactness
\begin{align}\label{surj1}
\bbC[\dg]^{\dG}=H^*_{\dG\times T_\rho}\otimes_{\bbC[a]}\bbC_{a=1}\twoheadrightarrow
\textstyle\bigoplus_\gamma H^*(\M_\Omega^\bullet(v_\gamma,w_\rho),\bbC).
\end{align}
On the summand indexed by $\gamma$, a symmetric polynomial $f$ in the Chern roots of the $V_c$ 
goes to $f$ evaluated at the $T_\rho$-equivariant Chern
roots of $\calV_c$, specialised at $a=1$.
By Lemma \ref{lem:fixed} these are the Chern roots of the $\calV_c(P)$ shifted by $P$. 
Each $H^*(\M_\Omega^\bullet(v_\gamma,w_\rho),\bbC)$ is a local Artinian algebra. We have 
$$\bbC[\dg]^{\dG}=\bbC[\mathfrak t]^{\mathfrak S},\quad \mathfrak S=\prod_c\mathfrak S_{\dim V_c}$$ 
The map \eqref{surj1} is a product of surjections, each factoring through the completion at some point 
$p_\gamma\in\mathfrak t/\mathfrak S$
\[
\widehat{\bbC[\mathfrak t]^{\mathfrak S}_{p_\gamma}}\twoheadrightarrow H^*(\M_\Omega^\bullet(v_\gamma,w_\rho),\bbC).
\]
The stabilizer of $p_\gamma$ in $\mathfrak S$ is $\mathfrak S_{v_\gamma}=\prod_{c,P}\mathfrak S_{v_c(P)}$,
since the positions are distinct integers.
Thus,
$$\widehat{\bbC[\mathfrak t]^{\mathfrak S}_{p_\gamma}}=\widehat{\bbC[\mathfrak t]^{\mathfrak S_{v_\gamma}}_{p_\gamma}},$$ 
which, after translating by $p_\gamma$, is $\widehat H^*_{\dG_{v_\gamma}}$, the map being the
Kirwan map. 
\end{proof}

For $\Omega=\dQ$ the graded
data are those of \S\ref{subsec:graded} with
$\dG_{v_\gamma}=\dG_\gamma$ 
and the map is $\kappa_\gamma$.

\subsection{Proof of Theorem \ref{thm:main}(a)}\label{subsec:proof-main-map}
For $D=1$, compose \eqref{OPS} with 
Proposition~\ref{prop:coh} to get the map
\begin{align}\label{PHG2}
 \Phi_\gamma:\widehat\B_\gamma\longrightarrow H^*(\M_{\dQ}^\bullet(v_\gamma,w_\rho),,\bbC).
\end{align}
Propositions~\ref{prop:psi}(c) and \ref{prop:cohZ} give
$\Phi_\gamma\circ r_\gamma=\kappa_\gamma$, hence
$\Ker(r_\gamma)\subseteq\Ker(\kappa_\gamma)$. \qed

\subsection{Proof of Theorem \ref{thm:main}(b)}\label{subsec:surj}
By Theorem \ref{thm:main}(a), we have
$$\kappa_\gamma=\Phi_\gamma\circ r_\gamma.$$
Further $\kappa_\gamma$ is surjective by Lemma
\ref{lem:kirwan}. Hence $\Phi_\gamma$ is surjective. \qed

\subsection{Proof of Theorem \ref{thm:main}(c)}\label{subsec:inj}
The ungraded quiver variety is $\M_Q(v,w)$ for the total dimension vectors, as in \S\ref{subsec:global-conjecture}.
In Lemma \ref{lem:fixed} we  take $\Omega=Q$.
We use the non equivariant dimension bound \eqref{C2} below as an
hypothesis. The cited Hikita results supply it in types
$A,D,E_6$, and in types $E_7,E_8$ under the framing restrictions
\begin{align}\label{FR}\begin{array}{ll}
 E_7: w_3=w_4=w_5=0,\quad
 E_8: w_2=w_3=w_4=w_5=w_6=w_7=0,
 \end{array}
\end{align}
with the vertex numbering of \cite[Cor.~2.17]{DK25}, componentwise
for a disconnected diagram. Theorem A.7 of \cite{DK25} has the
same hypotheses.  For other framings, the argument is conditional on \eqref{C2}. 
These results concern the
commutative algebra at $\hbar=0$, whereas $\widehat\B_\gamma$ is
quantized at $\hbar=1$. 
For the deformation argument we use the notation of
\S\ref{subsec:global-conjecture}: 
$$A=A_{F,\hbar}
,\quad
\B=\B^\theta(A)
,\quad
A_s=A\otimes_S\bbC_s
,\quad
s_\rho=(\rho,1).$$
We use the following two facts:
\begin{equation}\label{C1}
\parbox[c]{0.87\linewidth}{\normalfont
The algebra $A$ is free over the central algebra $H^*_{F\times\bbC^\times}$,
the specialization $A_0=A\otimes_S\bbC_0$ is the commutative Coulomb branch, and the Cartan quotient
$q_{F,\hbar}:R_{F,\hbar}\to\B$ is surjective, see \cite{BFN} and \cite[Thm.~B]{Wee19}.\par}
\end{equation}
\begin{equation}\label{C3}
\parbox[c]{0.87\linewidth}{\normalfont
The completion of $\B^\theta(A_{s_\rho})$ at the maximal ideal of
$H^*_{\dG}$ determined by $\gamma$ is $\widehat\B_\gamma$,
compatibly with $\widehat R_\gamma$. See  \eqref{H2}, Corollary~\ref{prop:localB} and Lemma~\ref{VQ}.\par}
\end{equation}
We also use the following hypothesis:
\begin{equation}\label{C2}
\parbox[c]{0.87\linewidth}{\normalfont
Assume $\dim\B^\theta(A_0)\leqslant\dim H^*(\M_Q(v,w),\bbC)$. 
 The bound is known in types
$A,D,E_6$, and in types $E_7,E_8$ under \eqref{FR}, by
\cite{KTWWY19a}, \cite{DK25}. Otherwise it is an hypothesis, see Remark~\ref{rem:gen-C2}.\par}
\end{equation}

\begin{Lemma}\label{lem:inj-Nak}
Assume $N_0:=\dim \B^\theta(A_0)<\infty$. Then:\hfill
\begin{enumerate}[label=$\mathrm{(\alph*)}$,leftmargin=8mm,itemsep=2mm]
\item
$\B$ is generated as an $H^*_{F\times\bbC^\times}$-module by $N_0$ homogeneous elements, 
\item 
$\B\otimes_S\bbC_s=\B^\theta(A_s)$ for all $s\in\operatorname{Spec}H^*_{F\times\bbC^\times}$,
\item
$\dim \B^\theta(A_s)\leqslant N_0$ for all $s\in\operatorname{Spec}H^*_{F\times\bbC^\times}$.
\end{enumerate}
\end{Lemma}

\begin{proof}
The algebra $\B$ is a quotient of $R_{F,\hbar}$ by a homogeneous ideal \eqref{C1}.
Hence $\B$ is nonnegatively graded. The B-algebra is a cokernel of maps between
charge-graded pieces. Thus, by freeness of $A$ and right exactness we have  
\begin{align}\label{BS}\B\otimes_{H^*_{F\times\bbC^\times}}\bbC_s=\B^\theta(A_s).\end{align}
For $s=0$ this is 
$$\B/(H^*_{F\times\bbC^\times})_+\B=\B^\theta(A_0).$$
 By the graded
Nakayama and a degree by-degree induction, $\B$ being bounded below, the lifts of a basis of 
$\B/(H^*_{F\times\bbC^\times})_+\B$ generate $\B$.
\end{proof}

\begin{proof}[Proof of Theorem $\ref{thm:main}(c)$]
The local maps have been constructed
using Theorem~\ref{thm:V} and Corollary \ref{prop:localB}.
By \eqref{BS} we have
$$\B_{s_\rho}=\B\otimes_{H^*_{F\times\bbC^\times}}\bbC_{s_\rho}=\B^\theta(A^\rho)$$
It  is a commutative quotient of $R_{s_\rho}=H^*_{\dG}$. 
It is finite dimensional by Lemma \ref{lem:inj-Nak}. Hence it is the product of its localizations at
its finitely many maximal ideals.
Distinct types $\gamma$ give distinct maximal ideals $\mathfrak m_\gamma$. 
By \eqref{C3}, Lemma \ref{lem:inj-Nak} and \eqref{C2}, we have
\[
\sum_\gamma\dim\widehat\B_\gamma=\sum_\gamma\dim(\B_{s_\rho})_{\mathfrak m_\gamma}\leqslant\dim \B_{s_\rho}\leqslant N_0\leqslant\dim H^*(\M_Q(v,w),\bbC).
\]
By Lemmas \ref{lem:kirwan} and \ref{lem:fixed}, we have
$$\dim H^*(\M_Q(v,w),\bbC)=\dim H^*(\M_Q(v,w)^T,\bbC)=\sum_\gamma\dim H^*(\M_{\dQ}^\bullet(v_\gamma,w_\rho),\bbC)$$  
By Theorem \ref{thm:main}(b), for each $\gamma$ we have
$$\dim\widehat\B_\gamma\geqslant\dim H^*(\M_{\dQ}^\bullet(v_\gamma,w_\rho),\bbC)$$ 
So all inequalities are equalities, i.e., we have
$$\dim\widehat\B_\gamma=\dim H^*(\M_{\dQ}^\bullet(v_\gamma,w_\rho),\bbC)$$ 
and the surjection $\Phi_\gamma$ is bijective.
\end{proof}

\begin{Remark}
The proof of Theorem \ref{thm:main} shows that 
$\B^\theta(A^\rho)\cong H^*(\M_Q(v,w)^{T_\rho},\bbC).$ 
Hence the algebra $\B^\theta(A^\rho)$ is supported at integral points.
\end{Remark}

\section{Symmetric types}\label{sec:general}

In this section $(C,D)$ is of symmetric type.
Fix any pair $(\rho,\gamma)$ of reds and blacks
at integer positions as in Definition \ref{def:blacks}. 
We explain which results of Sections \ref{sec:cotangent}-\ref{sec:surj} extend. 
The shifted quantum affine algebras used in \cite{VVsh} is only defined in finite type.
In general we replace it by the quantized Coulomb branch of the quiver gauge theory.

\subsection{Setting}\label{subsec:gen-setting}
Use the colours, flags, order, incidence sheaves and algebras of
\S\ref{subsec:unfolded}, now with
\[
 \dN=\bigoplus_c\Hom(\dW_c,\dV_c)\oplus
 \bigoplus_{h:i\to j}\ \bigoplus_{r\in\bbZ/2}
 \Hom(\dV_{(i,r)},\dV_{(j,r+1)}).
\]
Every arrow, including a loop, changes parity. The same formulas
define $\Cyc(\rho)$, the local Kirwan maps and the graded varieties.
A loop contributes two maps $V_i(P)\to V_i(P-1)$.
The ungraded variety is $\M_Q(v,w)$, not the variety of the colour cover.

Let $A=A_{F,\hbar}$ as in \eqref{AFH}.
Let $s_\rho$ be the point of $\operatorname{Spec}H_{F\times\bbC^\times}^*$ given by the red positions and 
$\hbar=1$.
Let $\B^\theta(A^\rho)$ be the B-algebra of 
$A^\rho=A\otimes_{H_{F\times\bbC^\times}^*}\bbC_{s_\rho}$. 
We take $\theta$ with $\theta_i>0$ for all $i$, e.g., the determinant cocharacter $\theta_{\det}=(1,\dots,1)$. 
It is generic. 
Recall the Cartan quotient in \eqref{CQ}.
The completed Cartan quotient is the map
\begin{align}\label{CCQ2}q_{\rho,\gamma}:\widehat R_\gamma\twoheadrightarrow
\widehat{\B^\theta(A^\rho)}_\gamma\end{align}
The {charge} relative to $\gamma$ of a type $\eta$ of dimension $v$, is 
$$\mathrm{ch}(\eta,\gamma)=\big((|\eta|_i-|\gamma|_i)/2\big)_{i\in I}\in\bbZ^I$$
see \S\ref{subsec:unfolded}. We put
\[
\Up^\theta(\gamma)=\{\eta:\ \theta\cdot\mathrm{ch}(\eta,\gamma)>0\}.
\]
In this section $\widehat\B_\gamma$ denotes $e_\gamma\Delta(\rho,\gamma)$ as in \eqref{BD}.

\begin{Lemma}\label{VQ}
Assume $D=1$.
Let $\rho$ be integral. 
\begin{enumerate}[label=$\mathrm{(\alph*)}$,leftmargin=8mm]
\item The category of integral weight modules of $A^\rho$
is identified with the  category of $\bR$-modules. The identification
preserves products, the idempotents $e_\gamma$, the Cartan algebras
$\widehat R_\gamma=\widehat H^*_{P_\gamma}$ (completion at 0) and the charge grading. 
\item 
There is a continuous isomorphism satisfying $\ \vartheta_{\rho,\gamma}\circ q_{\rho,\gamma}=r_\gamma$
\[
 \vartheta_{\rho,\gamma}:
 \widehat{\B^\theta(A^\rho)}_\gamma
 \xrightarrow{\cong} \widehat\B_\gamma
\]
\end{enumerate}
\end{Lemma}

\begin{proof}
Part (a) is homological localization
\cite[Thm.~4.4, \S4.3]{Web24}, which does not require finite type.
For (b), apply the finite-level construction of
Corollary~\ref{prop:localB} to
$A^\rho/A^\rho A^\rho_{[>0]}$, whose degree-zero space is
$\B^\theta(A^\rho)$. Part (a) sends positive-charge relations to
factorizations through $\Up^\theta(\gamma)$.
The bounded-weight argument in \cite[Thm.~4.3.1(b)]{VV25}
annihilates cyclotomic idempotents at every finite level, so passing
to the cyclotomic quotient adds no relation. Taking $\gamma$-pieces
and inverse limits gives (b). The right Cartan actions give the
multiplicative and Cartan compatibilities, as in that corollary.
\end{proof}

\begin{Remark}\label{rem:VQ}\hfill
\begin{enumerate}[label=$\mathrm{(\alph*)}$,leftmargin=8mm]
\item 
By Lemma~\ref{lem:Z}, Proposition~\ref{prop:psi}(b),(d) and
continuity, for $\Lambda_{F_\gamma}^{\rm st}\ne\emptyset$ the map $\Psi_\gamma$ in \eqref{PSG} induces
\[
 \widehat\Psi_\gamma:\widehat\B_\gamma\longrightarrow H^*_{P_\gamma}(\Lambda_{F_\gamma}^{\rm st},\bbC).
\]
\item The references \cite{BFN} and \cite{Web24} are stated for arbitrary $(G,N)$.
The geometric localization theorem of \cite{VV25} extends from its
finite type setting to arbitrary finite quivers.
The fixed-point, localization and convolution arguments are unchanged.
The description of category $\scrO$ by quiver Hecke algebras in \cite{KTWWY19b} is specific to finite type and is not needed here.
The comparison of K-theoretic and homological completions in \cite[Thm.~D]{DK25}
holds for any $(G,N)$.
\end{enumerate}
\end{Remark}

\subsection{What holds in arbitrary symmetric types}\label{sec:gen-verbatim}
Assume $D=1$.
The following hold with the same proofs:
\begin{enumerate}[label=$\mathrm{(\alph*)}$,leftmargin=8mm]
\item $\S\ref{sec:cotangent}$: Minimal $\bbZ$-flag, infinitesimal uniqueness, smoothness, Morse module, 
and the second part of Lemma $\ref{lem:noncyc}$.
\item \S\S$\ref{subsec:retraction}$-$\ref{subsec:Zst-coh}$:
Lemmas~$\ref{lem:leading}$-$\ref{lem:contr}$ and Propositions~$\ref{prop:coh}$, $\ref{prop:cohZ}$ give, 
together with Lemma~$\ref{lem:Z}$, an isomorphism
$H^*_{P_\gamma}(\Lambda^{\rm st}_{F_\gamma},\bbC)\cong H^*(\M_{\dQ}^\bullet(v_\gamma,w_\rho),\bbC)$
compatible with the Kirwan map.
\item $\S\ref{sec:micro}$: Proposition $\ref{prop:psi},$
\item $\S\ref{sec:surj}$: Lemmas $\ref{lem:fixed}$ and $\ref{lem:kirwan}$ for $\Gamma=Q$.
In particular $\kappa_\gamma$ is surjective and $\M_Q(v,w)$, $\M_{\dQ}^\bullet(v_\gamma,w_\rho)$ have
no odd cohomology.
\end{enumerate}

Indeed, in the proofs of parts (a) and (c) we never used the type of $Q$. 
For (b), the only input is the table of minimal centred degrees of \S\ref{subsec:graded}.
This is all that Lemmas \ref{lem:leading}-\ref{lem:contr} use. Lemma \ref{lem:free} is general, and
$$\Lambda_{F_\gamma}^{\rm st,\bullet}/\dG_\gamma\cong\M_{\dQ}^\bullet(v_\gamma,w_\rho)$$ 
is the Nakajima graded quiver variety of $Q$ \cite[\S4]{Nak01}. 
For (d) the Kirwan surjectivity
theorem of \cite{MN18} holds for arbitrary quivers, and the proofs of Lemmas \ref{lem:fixed} and \ref{lem:kirwan} 
do not use the type: for $\Gamma=Q$ the positions are arbitrary integers, and $\dG_v=\dG_\gamma$.

\begin{Theorem}\label{thm:gen-surj}
Assume $D=1$. For each
pair $(\rho,\gamma)$ we have:
\begin{enumerate}[label=$\mathrm{(\alph*)}$,leftmargin=8mm]
\item There is a surjective ring homomorphism 
$\Phi_\gamma:\widehat\B_\gamma\twoheadrightarrow H^*(\M_{\dQ}^\bullet(v_\gamma,w_\rho),\bbC)$
with $\Phi_\gamma\circ r_\gamma=\kappa_\gamma.$
\item 
$\M_{\dQ}^\bullet(v_\gamma,w_\rho)\neq\emptyset\Rightarrow\zeta^{2\gamma}\in\Spec \B^\theta(A^\rho)$ 
(i.e., there is a simple module in the category $\scrO$ of $A^\rho$ with highest weight $\gamma$).
\end{enumerate}
\end{Theorem}

\begin{proof}\hfill

(a) The proofs of Theorem \ref{thm:main}(a) and (b) apply, using \S\ref{sec:gen-verbatim}, 
Lemma \ref{VQ} in place of Theorem~\ref{thm:V} and Corollary \ref{prop:localB}, and Remark \ref{rem:VQ}(a). 

(b) follows from
(a): a nonempty graded quiver variety has nonzero cohomology, so
its unital comparison map forces $\widehat\B_\gamma\ne0$.
Lemma~\ref{VQ} identifies this with the existence of the indicated
highest weight simple. This implication uses no finite type assumption.
\end{proof}

\begin{Proposition}\label{prop:gen-points}
Assume $D=1$. For $i,j\in I$ let $a_{ij}$ be the number of arrows of the double quiver $\bar Q$ from $j$ to $i$. 
Suppose that for all $i\in I$ and
$m\in\bbZ$ we have
\begin{equation}\label{eq:gen-ldom}
w_i(m)+\sum_{j\in I}a_{ij}\,v_j(m+1)-v_i(m+2)-v_i(m)\ \geqslant\ 0 .
\end{equation}
Then, we have $\Lambda_{F_\gamma}^{\rm st,\bullet}\neq\emptyset$ and $\M_{\dQ}^\bullet(v_\gamma,w_\rho)\neq\emptyset$.
\end{Proposition}

\begin{proof}
Repeat the descending degree construction of
Proposition~\ref{prop:topdown}, counting each arrow separately.
Explicitly, set $A^*=0$ and
\begin{align*}
E_i^{(m)}&=W_i(m)\oplus\bigoplus_{h\in\bar Q,\,t(h)=i}V_{s(h)}(m+1),\\
 \tau_i^{(m)}&=(A_i^{(m)},(B_h^{(m)})_{t(h)=i}):E_i^{(m)}\to V_i(m),\\
 \sigma_i^{(m)}&=(0,(\varepsilon(\bar h)B_{\bar h}^{(m+1)})_{t(h)=i}):
 V_i(m+2)\to E_i^{(m)}.
\end{align*}
The moment equation is $\tau_i^{(m)}\sigma_i^{(m)}=0$.
For a loop, the two opposite arrows contribute to the commutator.
Once target degrees above $m$ have been treated, the map $\sigma_i^{(m)}$
is fixed, and \eqref{eq:gen-ldom} gives
\[
 \dim\operatorname{coker}\sigma_i^{(m)}
 \geqslant w_i(m)+\sum_j a_{ij}v_j(m+1)-v_i(m+2)
 \geqslant v_i(m).
\]
Choose the map $\tau_i^{(m)}$ to be a surjection from this cokernel.
Each arrow is chosen exactly once, being completely determined by its target degree.
It enters
another $\sigma$ only at the next step. Thus no earlier equation
changes, even for loops. Since $V$ and $W$ are finite dimensional  the induction terminates.
Surjectivity of all $\tau_i^{(m)}$ implies, again by induction on descending degrees,
that the framing generates $V$. Lemma~\ref{lem:stab}, via
\S\ref{sec:gen-verbatim}(b), gives $F_\xi=F_\gamma$.
The dimension estimate is independent of orientation signs, proving
all assertions.
\end{proof}

For $Q$ without loops and without multiple arrows, 
the matrix $(a_{ij})$ is the adjacency matrix and \eqref{eq:gen-ldom} is the inequality in the proof of
Proposition \ref{prop:topdown}. 
The entries $2-a_{ii}$ and $-a_{ij}$ for $i\ne j$ form the
symmetric quiver Cartan matrix. With loops its diagonal can be
nonpositive, so it is of Borcherds type. Condition \eqref{eq:gen-ldom}
says that the exponent of $Z_{i,m+1}$ in $Z^{w}A^{-v}$ is nonnegative, with 
$$A_{i,a}=Z_{i,a-1}Z_{i,a+1}\prod_jZ_{j,a}^{-a_{ij}}$$
This is the $\ell$-dominance condition for the
$q$-Cartan matrix of $Q$. For example, for the Jordan quiver ($I=\{i\}$, one loop, $a_{ii}=2$) it reads 
$$w(m)+2v(m+1)-v(m+2)-v(m)\geqslant0.$$
In finite type it is
the hypothesis of Proposition \ref{prop:topdown}. It does not use Nakajima's criterion for regular strata, which is a finite type statement.

\subsection{The global comparison map}
\begin{Theorem}\label{thm:gen-global-map}
Assume $D=1$.
The kernel of the Cartan quotient \eqref{C1}
$$q_{F,\hbar}:R_{F,\hbar}\to \B^\theta_{F,\hbar}$$ 
is contained in the kernel of the global Kirwan map \eqref{GKM}
$$\kappa_{F,\hbar}:R_{F,\hbar}\to H^*_{F\times\bbC^\times}(\M_Q(v,w),\bbC).$$ 
The map $\kappa_{F,\hbar}$ induces a surjective homomorphism of graded $R_{F,\hbar}$-algebras
\[
\Phi_{F,\hbar}:\ \B^\theta_{F,\hbar}\twoheadrightarrow H^*_{F\times\bbC^\times}(\M_Q(v,w),\bbC).
\]
\end{Theorem}

\begin{proof}
Put $\mathcal A=F\times\bbC^\times$ and
$S=H^*_{\mathcal A}=\bbC[\lambda_1,\ldots,\lambda_m,\hbar]$. 
All degrees in this proof are cohomological, not charge
degrees. The charge relations are homogeneous, so it suffices to
show $\kappa_{F,\hbar}(f)=0$ for homogeneous $f\in\Ker (q_{F,\hbar})$.
By \S\ref{sec:gen-verbatim}(d), equivariant formality
gives a finite homogeneous $S$-basis $e_a$ of $H^*_{\mathcal A}(\M_Q(v,w),\bbC)$.
Write $\deg f=2d$, $\deg e_a=2r_a$ and
\[
 \kappa_{F,\hbar}(f)=\sum_a p_a(\lambda,\hbar)e_a,
 \quad\deg p_a=d-r_a.
\]
For integral $\rho$, specialize at $s=s_\rho=(\rho,1)$.
Right exactness of the charge quotient gives
$\B\otimes_S\bbC_s=\B^\theta(A_s)$.
Thus $f_s$ vanishes in this B-algebra and, by Lemma~\ref{VQ}(b),
its image $f_{s,\gamma}$ vanishes in every $\widehat\B_\gamma$.
The local comparison of Theorem~\ref{thm:gen-surj}(a) yields
\begin{align}\label{form1}
 \kappa_\gamma(f_{s,\gamma})=0.
\end{align}

Lemma~\ref{lem:fixed} gives
$\M_Q(v,w)^{T_\rho}=\bigsqcup_\gamma\M_{\dQ}^\bullet(v_\gamma,w_\rho)$
for the one parameter group $T_\rho:t\mapsto(\rho(t^{-1}),t^{-1})$.
Localization at $s$ and the splitting
$\mathcal A=F\times T_\rho$ identify
\[
 H^*_{\mathcal A}(\M_Q(v,w),\bbC)\otimes_S\bbC_s
 \cong\bigoplus_\gamma H^*(\M_{\dQ}^\bullet(v_\gamma,w_\rho),\bbC).
\]
Here the fixed components are equivariantly formal by
\S\ref{sec:gen-verbatim}(d), and $s$ has coordinates
$(0,1)$ in the splitting. On the $\gamma$-summand the Kirwan map $\kappa_{F,\hbar}$ specializes to 
$$R_s\to\widehat R_\gamma\xrightarrow{\kappa_\gamma}H^*(\M_{\dQ}^\bullet(v_\gamma,w_\rho),\mathbb C).$$
Equation~\eqref{form1} therefore implies $p_a(\rho,1)=0$ for all
integral $\rho$. Density of $\bbZ^m$ gives $p_a(\lambda,1)=0$.
Homogeneity then gives
$p_a(\lambda,\hbar)=\hbar^{d-r_a}p_a(\lambda/\hbar,1)=0$
for $\hbar\ne0$, hence identically. Thus $\Ker (q_{F,\hbar})\subseteq\Ker(\kappa_{F,\hbar})$.
The resulting map $\Phi_{F,\hbar}$ is surjective because $\kappa_{F,\hbar}$ is surjective
by \cite{MN18}, equivariant formality and graded Nakayama.
\end{proof}

\begin{Remark}
The map $\Phi_{F,\hbar}$ 
in Theorem \ref{thm:gen-global-map} is the comparison map due to Botta and Tamagni, for the flavour torus $F=T_W$.
Our proof uses neither quasimaps nor their construction.
\end{Remark}

\subsection{Injectivity}
The proof of Theorem \ref{thm:main}(c) uses the type of $Q$ only through \eqref{C2}, which we now take as a hypothesis.

\begin{Theorem}\label{thm:gen-inj}
Assume $D=1$.
Assume condition \eqref{C2} holds. 
Then for each $\rho$ and each type $\gamma$ of dimension $v$, the maps 
$$\Phi_\gamma:\widehat\B_\gamma\to H^*(\M_{\dQ}^\bullet(v_\gamma,w_\rho),\bbC)
,\quad
\Phi_{F,\hbar}:\B^\theta_{F,\hbar}\to H^*_{F\times\bbC^\times}(\M,\bbC)$$
are isomorphisms. 
Further, we have 
 $$\gamma\in\Spec \B^\theta(A^\rho)\iff\M_{\dQ}^\bullet(v_\gamma,w_\rho)\neq\emptyset.$$
\end{Theorem}

\begin{proof}
 Lemma \ref{lem:inj-Nak} holds for $Q$, since it only uses \eqref{C1}, which holds for each quiver, and so does the proof of Theorem \ref{thm:main}(c) in
\S\ref{subsec:inj}, with Lemma \ref{VQ}, Theorem \ref{thm:gen-surj} and \S\ref{sec:gen-verbatim}(d) in place of \eqref{C3}, Theorem \ref{thm:main}(b) and Lemma \ref{lem:kirwan}.
This yields the local isomorphisms.
By Theorem \ref{thm:gen-global-map}, the map $\kappa_{F,\hbar}$ induces a surjection 
$$\Phi_{F,\hbar}:\B^\theta_{F,\hbar}\to H^*_{F\times\bbC^\times}(\M,\bbC)$$ 
By Lemma
\ref{lem:inj-Nak} there is a surjection
$(H^*_{F\times\mathbb C^\times})^{N_0}\twoheadrightarrow \B^\theta_{F,\hbar}$.
By \eqref{C2} we have
$$H^*_{F\times\bbC^\times}(\M,\bbC)\cong (H^*_{F\times\mathbb C^\times})^{N},\quad
N=\dim H^*(\M,\bbC)\geqslant N_0$$  The composite
$$(H^*_{F\times\mathbb C^\times})^{N_0}\twoheadrightarrow \B^\theta_{F,\hbar}\twoheadrightarrow (H^*_{F\times\mathbb C^\times})^N$$
 is a surjection of free modules, 
thus  $N_0=N$, and a surjective endomorphism of a finitely
generated module is an isomorphism. Therefore, both maps are isomorphisms.

\end{proof}

\begin{Remark} \label{rem:gen-C2} Condition \eqref{C2} is known in the following cases.\hfill
\begin{enumerate}[label=$\mathrm{(\alph*)}$,leftmargin=8mm]
\item Types $A,D,E_6$, and types $E_7,E_8$ with the framing
restrictions of \cite[Cor.~2.17, Rem.~2.18, Thm.~A.7]{DK25},
as specified in \S\ref{subsec:inj}. 
\item The Jordan quiver: the case $w=1$ is Hikita's original theorem \cite{Hik17}, the case
of general $w$ considered in \cite{KS24}. 
\item  In \cite{HikNotes} condition \eqref{C2} is derived from an IA computation for all quivers from a dimension count for
quivers with a loop at each vertex together with a loop-removal argument. 
That argument uses the comparison map for the quiver $Q$ with an added loop at each vertex and an additional scaling flavour 
$\bbC^\times$. 
\end{enumerate}
\end{Remark}

\section{Non symmetric types}\label{sec:ns}

In this section $(C,D)$ is of any finite type.

\subsection{The unfolded geometry}\label{subsec:ns-geom}
The data of \S\ref{subsec:unfolded} are defined for any finite type Cartan pair. 
Colours of vertex $i$ are classes of centred positions modulo $2d_i$.
An arrow $h:i\to j$ gives a component $\Hom(\dV_c,\dV_{c'})$ of $\dN$ for each pair of colours 
$c=(i,r),$ $c'=(j,s)$ satisfying \eqref{parity}. 
The underlying graph of the unfolded quiver $\dQ$ of \S\ref{subsec:unfolded} has the colours as vertices and an
edge for each such pair. Call it the \emph{unfolded graph}.
The minimal symmetrizers are
\[
 D_{B_n}=\operatorname{diag}({2,\ldots,2},1),\quad
 D_{C_n}=\operatorname{diag}({1,\ldots,1},2).
\]
By \eqref{parity}, the unfolded graph is two disjoint copies of
$A_{2n-1}$ for $B_n$, of $D_{n+1}$ for $C_n$,
of $E_6$ for $F_4$, and of $D_4$ for $G_2$.
In the symmetric loopless Dynkin case, it is two copies of $Q$.

The new feature is that the minimal centred degrees are not constant. 
A component $B_h$ from a basis vector at position $P$ and colour of vertex $i$ may
land at any position $Q<P$ with a colour of vertex $j$. The largest one is $P-\delta$ with $\delta$ the least positive element of
$P-s+2d_j\bbZ$, which depends on $P$ modulo $2d_j$ when $d_j>d_i$. 
The framing $A^*$ has centred degrees at least $2d_i$. 
Components of non minimal centred degree correspond to the positive powers of $z_i$ in the loop
space, and to the action of $\varepsilon_i$ in the generalized preprojective algebras of \cite{GLS}.

\subsection{What holds in non symmetric types}
The arguments of
Lemmas~\ref{lem:Z}, \ref{lem:stable-irreducible}, \ref{lem:free}(a)
and Propositions~\ref{prop:cohZ}, \ref{prop:psi}(a)-(c) apply without change.
Only the finiteness required in Proposition~\ref{prop:psi}(d)
remains to be checked.

\begin{Lemma}\label{F}
Let $(C,D)$ be of finite type.
\begin{enumerate}[label=$\mathrm{(\alph*)}$,leftmargin=8mm]
\item 
The group $\dG$ acts freely and properly on $Z^{\rm st}_\gamma$ with quotient $\mathfrak L_{\rho,\gamma}$.
\item
The variety $Z^{\rm st}_\gamma$ and $\mathfrak L_{\rho,\gamma}$ are quasi-projective smooth
and irreducible when nonempty.
\item
$Z^{\rm st}_\gamma\to \mathfrak L_{\rho,\gamma}$ is a principal $\dG$-bundle. 
\item
We have  isomorphic finite dimensional graded algebras with nilpotent positive degree part
\[
H^*_{P_\gamma}(\Lambda^{\rm st}_{F_\gamma},\bbC)\cong H^*_{\dG}(Z^{\rm st}_\gamma,\bbC)\cong H^*\big(\mathfrak L_{\rho,\gamma},\bbC\big)
\]
\end{enumerate}
\end{Lemma}

\begin{proof}
The space $T^*\dN$ is the doubled framed representation space of
$\dQ$, with gauge group $\dG$ and moment map $\mu$. Put
\begin{align}\label{eq:ungraded-generated-stable}
\Lambda^{\rm st}=\{\xi\in T^*\dN:\mu(\xi)=0,\ 
\dV\ \text{is spanned by all}\ w(\xi)A(u)\,,\,u\in\dW\,,\, w\ \text{a word in }x,y\}.
\end{align}
This is the ungraded analogue of \eqref{LB}. Its quotient
$\Lambda^{\rm st}\to\mathfrak M_{\dQ}$ is the principal gauge
bundle of the Nakajima variety for generated stability
\cite{Kin94}, \cite[\S3]{Nak01}.
By Lemmas~\ref{lem:words} and \ref{lem:Z},
\[
 Z_\gamma^{\rm st}=\Lambda^{\rm st}\cap Z_\gamma\cap\Omega_\gamma
\]
is an invariant locally closed subset. It is saturated for the
geometric quotient, so its quotient is locally closed in
$\mathfrak M_{\dQ}$ and its projection is the restriction of the
principal bundle. The smoothness of $Z_\gamma^{\rm st}$ is
Lemma~\ref{lem:Z}(d), and that of the quotient follows by smooth
descent along the principal bundle. Irreducibility follows from
Lemma~\ref{lem:stable-irreducible}, proving (a)-(c).
Equivariant cohomology therefore equals the  cohomology
of this quasi-projective quotient, which is finite dimensional
and has nilpotent positive-degree ideal. The other identification
is Proposition~\ref{prop:cohZ}.
\end{proof}

If $D=1$ then Proposition~\ref{prop:coh} compares the cohomology of 
$\mathfrak L_{\rho,\gamma}$ and $\M_{\dQ}^\bullet(v_\gamma,w_\rho)$. 
If $D\neq 1$ there is no analog of $\M_{\dQ}^\bullet(v_\gamma,w_\rho)$, so we use $\mathfrak L_{\rho,\gamma}$. 
The Kirwan map of $\Lambda^{\rm st}_{F_\gamma}$ is the restriction map
\begin{equation}\label{eq:ns-local-kirwan}
\kappa^{\rm ns}_\gamma:\ \widehat R_\gamma\longrightarrow H^*_{P_\gamma}(\Lambda^{\rm st}_{F_\gamma},\bbC)
\end{equation}
Its image is the tautological subalgebra, i.e., the subalgebra generated by the tautological classes.
Recall the map $r_\gamma:\widehat R_\gamma\to \widehat\B_\gamma$ in \eqref{RG}.

\begin{Theorem}\label{thm:ns-map}
Let $(C,D)$ be of finite type. For every integral pair $(\rho,\gamma)$ the following hold.
\begin{enumerate}[label=$\mathrm{(\alph*)}$,leftmargin=8mm]
\item $\Ker(r_\gamma)\subseteq\Ker(\kappa^{\rm ns}_\gamma)$.
There is a ring homomorphism
\[
\Phi^{\rm ns}_\gamma:\ \widehat{\B(\calA^\rho_\mu)}_\gamma\longrightarrow H^*_{P_\gamma}(\Lambda^{\rm st}_{F_\gamma},\bbC)
\text{ with } \Phi^{\rm ns}_\gamma\circ r_\gamma=\kappa^{\rm ns}_\gamma.
\]
\item 
$\mathrm{Im}(\Phi^{\rm ns}_\gamma)=\mathrm{Im}(\kappa^{\rm ns}_\gamma)$.
\item 
$\Phi^{\rm ns}_\gamma\neq0\iff \Lambda^{\rm st}_{F_\gamma}\neq\emptyset\iff (\rho,\gamma)
\text{ satisfies condition \eqref{bfE} for the orientation } Q.$
\item
For every weight block $\Gamma$ we have $\Gamma^{\rm st}=\Gamma$. More precisely, we have
\[
 Z_\gamma^{\rm st}\ne\emptyset
 \quad\Longleftrightarrow\quad
 \zeta^{2\gamma}\in\Spec\B(\calA^\rho_\mu)
 \quad\Longleftrightarrow\quad
 L(\rho,\gamma)\in\scrO_\mu^\rho.
\]

\end{enumerate}
\end{Theorem}

\begin{proof}\hfill

(a) is  Theorem \ref{thm:main}(a), with Proposition \ref{prop:psi} and Lemma \ref{F} in place of Proposition \ref{prop:coh}.

(b) follows from (a) and the surjectivity of the map $r_\gamma:\widehat R_\gamma\to\widehat\B_\gamma$.

 (c) $\Phi^{\rm ns}_\gamma$ is a unital ring map, so it is nonzero iff its target is nonzero.
 
(d) The B-algebra classification in \S\ref{subsec:algebras} gives
the second equivalence. Theorem~\ref{thm:point}(a) gives the forward
implication in the first, and Proposition~\ref{prop:microlocal-basis}(a)
gives its converse.

\end{proof}

\subsection{A surjection onto a graded quiver variety of $\dQ$}\label{subsec:ns-dQ}
Use the centred gradings of \S\ref{subsec:graded} on the unfolded quiver $\dQ$. A point $\xi\in T^*\dN$ is said
\emph{$\dQ$-homogeneous} when both edge directions have degree $-1$,
the incoming framing degree $0$, and the outgoing framing degree $-2$.
Let $\Lambda_{\dQ}^{\rm st}$ be the locus of such points satisfying
$\mu^\bullet=0$ and generation by the framing.
Its quotient is the graded quiver variety $\M^\bullet_{\dQ}(v_\gamma,w_\rho)$ of
\S\ref{subsec:ungraded}, with 
$$v_c(P)=\text{ number of blacks of colour $c$ at $P$}
,\quad
w_c(L)=\text{ number of reds of colour $c$ at $L$}$$ 
The Kirwan map is
\[
 \kappa_\gamma^{\dQ}:\widehat R_\gamma\longrightarrow
 H^*(\M^\bullet_{\dQ}(v_\gamma,w_\rho),\bbC).
\]

\begin{Theorem}\label{thm:ns-dQ}
Assume that $d_{ij}=1$ for each edge $i-j$ of $Q$. 
\begin{enumerate}[label=$\mathrm{(\alph*)}$,leftmargin=8mm]
\item There is a surjective ring homomorphism
\[
\Phi^{\dQ}_\gamma:\ \widehat{\B(\calA^\rho_\mu)}_\gamma\twoheadrightarrow H^*\big(\M^\bullet_{\dQ}(v_\gamma,w_\rho),\bbC\big),\quad
\Phi^{\dQ}_\gamma\circ r_\gamma=\kappa^{\dQ}_\gamma .
\]
\item
If $\M^\bullet_{\dQ}(v_\gamma,w_\rho)\neq\emptyset$ then $\zeta^{2\gamma}\in\Spec\B(\calA^\rho_\mu)$.
This holds if $(\rho,\gamma)$ is
{$\dQ$-dominant}, i.e., if 
$$w_c(k-1)-v_c(k-1)-v_c(k+1)+\sum_{c'-c}v_{c'}(k)\geqslant0,\quad\forall k,c$$ \end{enumerate}
\end{Theorem}

\begin{proof}
Let $\xi\in\Lambda_{\dQ}^{\rm st}$. Homogeneity satisfies the
four inclusions \eqref{eq:NF}, \eqref{eq:NFperp}, hence
\begin{align}\label{step1}
 \xi\in L_{F_\gamma}.
\end{align}
Every moment map term has degree $-2$, so $\mu^\bullet=0$ if and only if $\mu=0$.
By Lemma~\ref{lem:words}, $F_\xi\subseteq F_\gamma$.
The assumption $d_{ij}=1$ makes adjacent colours have opposite
parities. A nonzero homogeneous arrow from position $P$ ends at
$P-1$, the largest permitted target position for either bound rule.
Thus a nonzero framed word $wA(u)$ with $u\in W_c(L)$ and
$\operatorname{length}(w)=m$ has both position and bound $L-m$.
Generation by these homogeneous words and Lemma~\ref{lem:words}
give the reverse flag inclusion. Therefore
\begin{align}\label{LQLF}
 \Lambda_{\dQ}^{\rm st}\subseteq\Lambda_{F_\gamma}^{\rm st}.
\end{align}

Restriction to the Levi subgroup $\dG_\gamma\subset P_\gamma$
and then to this locus gives
\begin{align}\label{HLM}
\begin{aligned}
 H^*_{P_\gamma}(\Lambda_{F_\gamma}^{\rm st},\bbC)
 \xrightarrow{\cong}H^*_{\dG_\gamma}(\Lambda_{F_\gamma}^{\rm st},\bbC)
 \rightarrow H^*_{\dG_\gamma}(\Lambda_{\dQ}^{\rm st},\bbC)\cong H^*(\M^\bullet_{\dQ}(v_\gamma,w_\rho),\bbC).
\end{aligned}
\end{align}
These maps preserve the graded tautological classes.
Composing with $\Phi_\gamma^{\rm ns}$ gives
$\Phi_\gamma^{\dQ}\circ r_\gamma=\kappa_\gamma^{\dQ}$.
The latter is surjective by Lemma~\ref{lem:kirwan}, applied
to $\dQ$, proving (a) without assuming surjectivity of the
intermediate restriction.

A nonempty target gives $\widehat\B_\gamma\ne0$, hence $\gamma\in\Spec\B$ by Corollary~\ref{prop:localB}.
Finally, substituting $m=k-1$ in the $\dQ$-dominance inequalities
gives exactly the hypotheses of Proposition~\ref{prop:gen-points} for $\dQ$ with $D=1$. 
Hence, under this condition
$\Lambda_{\dQ}^{\rm st}\neq\emptyset$, proving (b).
\end{proof}

\begin{Remark}\hfill
\begin{enumerate}[label=$\mathrm{(\alph*)}$,leftmargin=8mm]
\item The hypothesis $d_{ij}=1$ on every edge in
Theorem \ref{thm:ns-dQ} holds for the minimal symmetrizer in type
$C_n$: the symmetrizing integers are $1,\ldots,1,2$, with the
vertex of weight $2$ at an end of the Dynkin diagram. In type
$B_n$ the hypothesis fails for $n\geqslant3$.


\item In the non symmetric case $\M^\bullet_{\dQ}(v_\gamma,w_\rho)$  can be empty even when condition \eqref{bfE} holds. 
In particular, the $\dQ$-dominance is not implied by $\ell$-dominance. 
Even if  $d_{ij}=1$ holds that theorem does not establish Theorem \ref{Theorem:431} 
for every $\ell$-dominant $(\rho,\gamma)$. 

\item
The map \eqref{HLM} is not an isomorphism in general.
In type $B_2$ with the datum of Proposition~\ref{prop:B2-local-kirwan}, the local
target is $\bbC^\times\times\mathbb A^2$, whereas
$\M^\bullet_{\dQ}(v_\gamma,w_\rho)$ is empty. 

\item
In every finite type, the variety
$\mathfrak L_{\rho,\gamma}$
embeds as a smooth locally closed Lagrangian subvariety of $\mathfrak M_{\dQ}(v,w)$.
Indeed, the embedding follows from Lemma~\ref{F}, the isotropy descends by
symplectic reduction from Lemma~\ref{lem:Z}, and we have
\[
 \dim\mathfrak L_{\rho,\gamma}=\dim\dN-\dim\dG
 =\tfrac12\dim\mathfrak M_{\dQ}(v,w).
\]
\item 
For $D=1$, Lemma~\ref{lem:attracting-stratum} identifies the variety $\mathfrak L_{\rho,\gamma}$
with an attracting stratum in $\mathfrak M_{\dQ}(v,w)$ with for the action of the one parameter subgroup $t\mapsto T_\rho(t)$. 
This action still exists in non symmetric types, but
$\mathfrak L_{\rho,\gamma}$ may not be an attracting stratum anymore. 
In the example in (d), the fixed locus is empty and
$H^1(\mathfrak L_{\rho,\gamma},\bbC)=\bbC$, hence $\mathfrak L_{\rho,\gamma}$
cannot be the attracting stratum to any fixed component of any
one parameter action on $\mathfrak M_{\dQ}(v,w)$ because
such varieties have no odd cohomology by \cite{MN18}.
Indeed, the varying minimal centred degrees obstruct the contraction of
Lemma~\ref{lem:contr}. 
An arrow component $T:V_i(P)\to V_j(Q)$ has centred degree $P-Q$.
In type $B_2$, with $i$ the long vertex ($d_i=2$) and $j$ the short one ($d_j=1$),
a component $B$ along $j\to i$ has minimal centred degree
$u=1$ or $3$, and the reverse component has centred degree $4-u=2d_i-u$, because
a composable pair contributing to \(BB^*\) has the form 
$$V_i(P+4)\xrightarrow{B^*}V_j(P+u)\xrightarrow{B}V_i(P),$$ 
with \(u=1\) or \(3\). Its centred degree is therefore $4$.
The degree 3 component in Proposition~\ref{prop:B2-local-kirwan}
is essential for stability, so retaining only degree 1 components
can destroy stability. 
\end{enumerate}
\end{Remark}

\subsection{The variety $\mathfrak L_{\rho,\gamma}$ as a moduli space of graded modules}\label{subsec:ns-modules}
We follow the conventions in  \cite[\S2]{NW23},
\cite[\S3.1]{VVsh}.
Set
$$z_i=z^{d_i},\quad
\calO_i=\bbC\llb z_i\rrb
,\quad
K_i=\bbC\llp z_i\rrp.$$
For an edge $i-j$ we put  
$$\calO_{ij}=\bbC\llb z_{ij}\rrb,\quad z_{ij}=z_i^{|c_{ij}|}=z_j^{|c_{ji}|}$$ 
The matter of \cite{NW23}, \cite[(3.1.3)]{VVsh} and its loop spaces are
\[
N_\calO=\bigoplus_i\Hom(W_i,V_i)\otimes\calO_i\ \oplus\bigoplus_{h:i\to j}\Hom_{\calO_{ij}}\big(V_i\otimes\calO_i,V_j\otimes\calO_j\big),\quad N_K=N_\calO\otimes_{\calO}K,
\]
Let $N_K^*$ be the dual for the residue pairing and 
$$G_K=\prod_iGL(V_i)(K_i).$$
The $G_K$-action on $N_K$ yields a moment map
$$\mu_K:N_K\times N_K^*\to\mathfrak g_K^*$$ 
A {lattice} is a
family $L=(L_i)$ of $\calO_i$-lattices $L_i\subset V_i\otimes K_i$.
Let $$N_\calO(L)\subset N_K$$ be the regular points with respect to $L$, i.e., mapping $W\otimes\calO$ into $L$ and
$L_i$ into $L_j$. Let  
$$N_\calO(L)^\perp\subset N_K^*$$ be its annihilator.

Let $\widetilde\Pi=\widetilde\Pi(C^{\rm T},D)$ be the generalized preprojective
algebra with $H_i=\calO_i$ and $\varepsilon_i=z_i$.
In this notation $D$ records the loop degrees $2d_i$. It may not
be a GLS symmetrizer of $C^{\rm T}$. At truncation order
$N=kN_0$ the GLS pair is $(C^{\rm T},D_N)$ with
$D_N=\operatorname{diag}(N/d_i)$.
The transpose is the
convention change of \cite[Rem.~2.1]{NW23}. 
A \emph{framed $\widetilde\Pi$-module of rank $v$} is a lattice $L$ together with 
a pair $$(x,y)\in N_\calO(L)\times N_\calO(L)^\perp$$ such
that $\mu_K(x,y)=0$. It is {stable} if $L$ is generated by the image of the framing $A(W\otimes\calO)$, 
as a module over the algebras $\calO_i$ and the maps $B_h$ and $B^*_h$.
Two modules are isomorphic if they differ by an element of $G_K$.
The arrows of the double quiver act by the components $B_h$ of $x$ and the 
components $B^*_h$ of $y$. 
The $\varepsilon_i$ act as multiplication by $z_i$.
The relation (H2) of \cite{GLS} expresses that the $B_h$ are $\calO_{ij}$-linear. 
The preprojective relation is $\mu_K=0$ which is given by
\[
\sum_{h:j\to i}\sum_{k=0}^{f_{ij}-1}z_i^{f_{ij}-1-k}B_hB_h^*z_i^k-\sum_{h:i\to j}\sum_{k=0}^{f_{ij}-1}z_i^{f_{ij}-1-k}B^*_hB_hz_i^k+AA^*=0,
\]
see \cite[\S2]{NW23}. 
Here $f_{ij}=|c_{ij}|/\gcd(|c_{ij}|,|c_{ji}|)$ and the
multiplicities $g_{ij}$ of \cite{GLS} are $1$ in finite type.
The element $s=s_{\rho,\gamma}$ in \S\ref{subsec:data} acts on 
$$V\otimes K,\quad W\otimes K,\quad N_K,\quad N_K^*,$$
through $\zeta^{2\gamma}$, $\zeta^{2\rho}$ and the loop rotation, with the
twists of \cite[(3.1.5), (3.1.6)]{VVsh}. 
A framed module is {$s$-fixed} if $s$ maps it to an isomorphic one.
Its {type} is defined below.

\Needspace{7\baselineskip}
\begin{Proposition}\label{prop:ns-modules}
Let $(C,D)$ be of finite type. 
\begin{enumerate}[label=$\mathrm{(\alph*)}$,leftmargin=8mm]
\item Let $p=(L,x,y)$ be a stable framed $\widetilde\Pi$-module whose
isomorphism class is fixed by $s$. The $g_t$'s such that $t\cdot p=g_t\cdot p$ define a cocycle.
The operators $h_t=g_t^{-1}\circ t$ define a  torus
action preserving $L$. Each $L_i$ has a $\calO_i$-basis homogeneous
for this action, with $\deg z_i=2d_i$. 
The tuple of the degrees of the basis elements is a type $\gamma$. 
The action of $h_t$ is diagonal in this basis.
\item 
Up to isomorphism, stable fixed framed $\widetilde\Pi$-modules of type $\gamma$
are identified with points of $\mathfrak L_{\rho,\gamma}.$
\item Each matrix coefficient of an arrow with exponent larger than the
minimal one is obtained from a coefficient of minimal
exponent by multiplication by a positive power of $z_i$ or $z_j$.
\end{enumerate}
\end{Proposition}

Here a matrix coefficient is taken between specified homogeneous
weight spaces. We use the fixed lattice description of \cite{VV25}, \cite[\S2.4]{VVsh}.
The following lemma refines the lattice/flag identification of
Lemma~$\ref{lem:semilinear-flags}$.

\begin{Lemma}\label{lem:ns-lattice}\hfill
\begin{enumerate}[label=$\mathrm{(\alph*)}$,leftmargin=8mm]
\item For an $s$-stable lattice $L$ with $\bbZ$-flag $F$, we have
$$(N_K)^s=\dN,\quad N_\calO(L)^s=\dN_F,\quad (N_\calO(L)^\perp)^s=\dN_F^\perp.$$
\item We have $(G_K)^s=\dG$. The stabilizer of $L$ in $(G_K)^s$ is
$P_F$, the stabilizer of $F$.
\end{enumerate}
\end{Lemma}

\begin{proof}\hfill

(a) The first claim is \cite[\S2.4]{VVsh}, \cite{VV25}.
Regularity of $x$ with respect to $L$ means that each component maps
$L_{i,P}$ into $L_{j,Q}$ for all degrees, which in centred positions
is \eqref{eq:NF}. See \cite[(2.4.22)]{VVsh} and
\S\ref{subsec:data}. The residue pairing pairs $s$-weight $\lambda$
with $-\lambda$, and on the fixed parts gives the trace pairing
of $\dN$. Thus the fixed part of the annihilator is the
annihilator of $N_\calO(L)^s=\dN_F$.

(b) By Lemma~\ref{lem:semilinear-flags}, the group $G_K^s$ is the group
of graded automorphisms of $V\otimes K$, hence
$\prod_cGL(\dV_c)=\dG$. Such an automorphism preserves $L$ if and only if  it preserves each $F^{\leqslant P}$.
\end{proof}

\begin{proof}[Proof of Proposition $\ref{prop:ns-modules}$]\hfill

(a) An automorphism fixing $p$ fixes the framing and commutes with
all arrow and coefficient maps. Thus, it fixes the generated
lattice $L$ and its localization. Hence it is the identity. Thus $g_t$ is uniquely determined. 
The external torus acts on $G_K$ through loop rotation. Hence
\[
 (tu)\cdot p=t\cdot(g_u\cdot p)
          =(t\cdot g_u)g_t\cdot p,
\]
which gives the cocycle law.
Consequently
\[
 h_t h_u=g_t^{-1}(t\cdot g_u^{-1})\circ(tu)
        =g_{tu}^{-1}\circ(tu)=h_{tu}.
\]
This is a semilinear action on $L$ over coefficient rotation.
Diagonalize it first on $L_i/z_iL_i$. At each finite quotient
$L_i/z_i^nL_i$ the torus is linearly reductive. Thus a weight
basis lifts through the equivariant surjections to a compatible
weight basis. Nakayama's lemma makes these a free $\calO_i$-basis.
In that basis
\[
 h_t\bigl(f(z_i)e_a\bigr)=f(t^{2d_i}z_i)t^{P_a}e_a.
\]
The weights $P_a$ are the initial weights on $L_i/z_iL_i$, and are therefore
independent of the homogeneous lifts. They define the type $\eta$.
This establishes the diagonal normalization. 

(b) By Lemma \ref{lem:ns-lattice}, the $s$-fixed lattice $L$
determines a $\bbZ$-flag $F$ of type $\gamma$. The regularity
conditions on the arrow maps are $x\in\dN_F$ and
$y\in\dN_F^\perp$. Since the graded gauge group $\dG$ acts
transitively on flags of type $\gamma$, we may choose a graded gauge
transformation sending $F$ to the standard flag $F_\gamma$.
By the equivariant Frobenius identifications verified in
\S\ref{subsec:ns-lattice-action}, for $\xi=(x,y)$ we have
\[
 \mu_K(\xi)=0\quad\Longleftrightarrow\quad\mu(\xi)=0.
\]
Thus, after this change of frame, the pair $\xi$ belongs to $\Lambda_{F_\gamma}$. 

The submodule generated by the framing under the action of $z_i$, $B_h$ and
$B_h^*$ is graded. Its degree pieces form the smallest increasing
family containing the framing and satisfying \eqref{eq:Gdef}.
This is $F_\xi$. Therefore
module generation is equivalent to $F_\xi=F_\gamma$.

Conversely, starting from a point $\xi\in\Lambda_{F_\gamma}$, we can reconstruct the corresponding fixed lattice of type $\gamma$. Two 
such points define isomorphic fixed lattices if and only if when they differ by an element of $P_\gamma$. Indeed, an isomorphism of fixed lattices 
must preserve their gradings and their standard flag, so it belongs to $P_\gamma$. 
The graded lattice/flag identification commutes with changes of
homogeneous frame and base change, so the correspondence also holds in families.

(c) By the computation of the weights of $N_\calO$ in \cite[(3.1.3), (3.1.5)]{VVsh}, \cite[\S2]{NW23} we have
$$\Hom_{\calO_{ij}}(V_i\otimes\calO_i,V_j\otimes\calO_j)=\Hom(V_i,V_j)\otimes\bbC\llb z^{d_{ij}}\rrb\,z^{d_i+b_{ij}}$$ 
as a graded space. Thus, maps between two weight spaces have exponents
$$e=e_0+d_{ij}k,\quad k\geqslant0,$$ 
with $e_0=d_i+b_{ij}$ the minimal exponent. The centred degree of an arrow component of exponent $e$ is an increasing affine function of $e$ 
by \cite[(3.1.5), (3.1.6)]{VVsh}. Thus 
positive excess means $k>0$. In finite type one of $d_i,d_j$ divides the other, 
hence $z^{d_{ij}}=z_i$ or $z_j$, and a map of exponent
$e_0+d_{ij}k$ is the matrix coefficient of exponent $e_0$ multiplied by $z_i^k$ 
on the source or $z_j^k$ on the target. The same holds for $y$ and for the framing, where
$d_{ij}$ is replaced by $d_i$.
\end{proof}

\begin{Remark}\label{rem:ns-modules}\hfill
\begin{enumerate}[label=$\mathrm{(\alph*)}$,leftmargin=8mm]
\item The modules in Proposition \ref{prop:ns-modules} are free over $H_i=\calO_i$.
Thus they are modules over the untruncated algebra
$\widetilde\Pi(C^{\rm T},D)$ and not over a finite truncation
$\Pi(C^{\rm T},D_N)$.
The finite truncations are the global models $\mathfrak M_N$ of
\S\ref{app:target}. The space of all stable framed
$\widetilde\Pi$-modules is
infinite dimensional.
Its $s$-fixed components should be the finite dimensional
varieties $\mathfrak L_{\rho,\gamma}$.
\item If $D=1$ the homogeneous modules are Nakajima's graded quiver data and,
by Proposition \(\ref{prop:coh}\), the variety 
$\mathfrak L_{\rho,\gamma}$ is homotopy equivalent to the corresponding graded quiver variety.
\end{enumerate}
\end{Remark}

\subsection{Existence of stable points}\label{subsec:ns-existence}
Theorem~\ref{thm:ns-map}(d) characterizes non emptiness by existence
of a simple highest weight module. We now give a direct sufficient
criterion: $\ell$-dominance, without assuming that the simple exists.
This is the non symmetric analogue of Proposition~\ref{prop:gen-points}.
The argument uses all matrix units, not only their homogeneous part,
and does not impose $A^*=0$.

\begin{Theorem}\label{thm:ns-existence}\label{conj:E}
Let $(C,D)$ be of finite type.
For each pair $(\rho,\gamma)$ such that $Z_{\rho,\gamma}$ is $\ell$-dominant  we have
$Z^{\rm st}_\gamma\neq\emptyset.$
\end{Theorem}

To prove the theorem, we first prove a lemma on  preprojective algebras.
Then, we use the auxiliary  graph $\widehat H$ introduced in the proof of Lemma \ref{lem:stable-irreducible}, and
we apply the  geometric realization of the crystal $B_{\widehat H}(\infty)$
for symmetric generalized Cartan matrices in \cite{KS97} to this graph.
Note that $\widehat H$ is not of finite type in general.

Let $H$ be a finite graph without loops, with vertex set $J$ and symmetric
Cartan matrix
\[
C_{u,v}^H=2\delta_{uv}-\#\{\text{edges of $H$ between $u$ and $v$}\}.
\]
Choose an orientation. Let $\Pi(H)$ be the preprojective algebra. 
For a module $M$ and a vertex $u$ set
\[
\operatorname{top}_u(M)
=M_u\Big/\sum_{h:\,t(h)=u}\operatorname{Im}M_h,
\]
where the sum runs over the doubled graph. Recall from \eqref{K1}
that $K_u(M)$ is the kernel of the maximal quotient supported at $u$.
Thus
\[
M/K_u(M)\cong S_u^{\oplus\dim\operatorname{top}_u(M)},
\]
where $S_u$ is the simple module supported at $u$.

\begin{Lemma}\label{lem:ns-finite-word}
Let $\mathbf i=(i_1,\ldots,i_N)\in J^N$ and
$\mathbf a=(a_1,\ldots,a_N)\in\mathbb N^N$.
Suppose that for each consecutive occurrences $p<q$ of a vertex $u$ in $\mathbf i$ we have
\begin{equation}\label{eq:ns-word-inequality}
a_p+a_q+\sum_{p<k<q}C_{u,i_k}^Ha_k\leqslant0.
\end{equation}
Then, there is a finite dimensional nilpotent $\Pi(H)$-module $M$ such that
\begin{align}\label{eq:ns-word-filtration}
M^{(k-1)}/M^{(k)}\cong S_{i_k}^{\oplus a_k}
,\quad
M^{(0)}=M
,\quad 
M^{(k)}=K_{i_k}(M^{(k-1)})
,\quad 
M^{(N)}=0
,\quad
1\leqslant k\leqslant N\end{align}
\end{Lemma}

\begin{proof}
We use the strict embedding of crystals
\[
\Psi_i:B(\infty)\hookrightarrow B(\infty)\otimes B_i
\]
in \cite[Thm.~3.2.2]{KS97}.
Iterate these embeddings
\[
\Psi_{\mathbf i}:B(\infty)\hookrightarrow
B(\infty)\otimes B_{i_N}\otimes\cdots\otimes B_{i_1}.
\]
This identify elements of $B(\infty)$ with tensors in the right hand side.
 Let $b_i(-x)$ be the element of $B_i$ of weight $-x\alpha_i$ for each $x\in\bbZ$.
For a tensor whose first factor is $u_\infty$ let
$\mathbf x=(x_1,\ldots,x_N)$ be the coordinates of the remaining factors. 
We define the score of $\mathbf x$ at $p$ to be the integer
\[
\sigma_p(\mathbf x)=x_p+\sum_{j>p}C_{i_p,i_j}^H x_j.
\]
The score of $u_\infty$ is $0$.
The tensor product rule of \cite[\S2.2]{NZ97} implies that the Kashiwara operator
$\widetilde f_u$ acts on the positions $p$ with $i_p=u$. More precisely, it acts on
the smallest such position having maximal score $\sigma_p(\mathbf x)$. 
See also \cite[(2.32)-(2.33)]{NZ97}.

Starting with all coordinates equal to zero, apply the operators
\[
\widetilde f_{i_1}^{a_1},\quad
\widetilde f_{i_2}^{a_2},\quad\ldots,\quad
\widetilde f_{i_N}^{a_N}
\]
in this order. We claim that $\widetilde f_{i_k}^{a_k}$ increases
only the $k$-th coordinate, from $0$ to $a_k$.
Indeed, during the $k$-th step the coordinates have the form
\[
\mathbf x=(a_1,\ldots,a_{k-1},t,0,\ldots,0),
\quad 0\leqslant t<a_k.
\]
The score at position $k$ is $\sigma_k(\mathbf x)=t$. We must
check that no earlier occurrence of $i_k$ has a score as large as $t$.
For two consecutive occurrences $p<q$ with $i_p=i_q=i_k=u$, since $C_{u,u}^H=2$, we get
\[
\sigma_p(\mathbf x)-\sigma_q(\mathbf x)
=x_p+x_q+\sum_{p<j<q}C_{u,i_j}^H x_j.
\]
If \(q<k\), all coordinates between \(p\) and \(q\) have reached their final values, hence \eqref{eq:ns-word-inequality} gives 
$$\sigma_p(\mathbf x)\leqslant\sigma_q(\mathbf x).$$ 
If \(q=k\), the \(k\)-th coordinate is \(t<a_k\), hence \eqref{eq:ns-word-inequality} gives
\[
\sigma_p(\mathbf x)-\sigma_k(\mathbf x)
\leqslant t-a_k \]
By \eqref{eq:ns-word-inequality}, this difference is nonpositive
when $q<k$. When $q=k$, the same inequality gives
\[
\sigma_p(\mathbf x)-\sigma_k(\mathbf x)=a_p+t+\sum_{p<j<k}C^H_{u,i_j}a_j\\
\leqslant t-a_k<0.
\]
Applying these inequalities successively shows that every earlier
occurrence $p<k$ of $i_k$ has
\[
\sigma_p(\mathbf x)<\sigma_k(\mathbf x)=t.
\]
Every later occurrence has score $0$. Thus position $k$ is the
first position with maximal score, even when $t=0$. The
tensor product rule therefore increases $x_k$, proving the claim.

Repeating this argument for $k=1,\ldots,N$ gives
\[
u_\infty\otimes b_{i_N}(-a_N)
\otimes\cdots\otimes b_{i_1}(-a_1).
\]
Since the embeddings are strict, this tensor is the image of an
element $b\in B(\infty)$ under $\Psi_{\mathbf i}$.

Next, we apply \cite[Thm.~3.2.2]{KS97}.
 Let $b^*$ be the image of $b$ under Kashiwara's involution. Setting $b_0=b^*$ and
$b_k=\widetilde e_{i_k}^{a_k}b_{k-1}$, we obtain
\[
\varepsilon_{i_k}(b_{k-1})=a_k,\quad b_N=u_\infty.
\]
Further, the geometric realization \cite[Thm.~5.3.2]{KS97} identifies $b_k$ with
an irreducible component $Z_k$ of a nilpotent preprojective representation
variety. For a generic representation $M\in Z_k$, we have
$$\varepsilon_u(b_k)=\dim\operatorname{top}_u(M)$$ and
$\widetilde e_u^{\max}b_k$ is  the irreducible
component obtained generically by replacing $M$ with $K_u(M)$,
by \cite[\S5.1, Prop.~5.2.4]{KS97}.
The correspondences of \cite[Lem.~5.2.3]{KS97} which define the crystal are dominant on the
corresponding components. Hence, a generic representative in the
component of $b_0$ can be chosen so that its successive maximal top
kernels lie in the required generic loci for $b_1,\ldots,b_N$.
The final component has dimension vector zero, so the resulting module
satisfies \eqref{eq:ns-word-filtration}.
\end{proof}

\begin{proof}[Proof of Theorem $\ref{thm:ns-existence}$]
\hfill

Let $\mathscr C$ be the set of colours and $c\in\mathscr C$. 
For $c=(i,r)$ put $p_c=2d_i$
and let $v_c(P)$ and $w_c(P)$ be the black and red multiplicities at
position $P\in r+p_c\bbZ$. Let $H$ be the graph of $\dQ$.
Adjacent colours have no position in common,
by the condition \eqref{parity}.
Recall that $\ell$-dominance means that every exponent in the monomial 
$Z_{\rho,\gamma}$ of \eqref{ZRG} is nonnegative. 
In terms of the multiplicities $v_c(P)$ and $w_c(P)$, this condition is equivalent to
\begin{equation}\label{eq:ns-interval-dominance}
v_c(P)+v_c(P+p_c)
\leqslant w_c(P)+\sum_{c'\sim c}\ \sum_{P<Q<P+p_c}v_{c'}(Q)
\quad\text{for all $c,P$.}
\end{equation}
where we write $c'\sim c$ when $c'$ and $c$ are joined by an edge in the graph $H$.

The auxiliary graph $\widehat H$ is constructed as in the proof of Lemma \ref{lem:stable-irreducible}:
For each colour vertex $c$ and each position $P$ with $w_c(P)>0$,
we add to $H$ a new vertex $r_{c,P}$, joined to $c$ by exactly one edge and
joined to no other vertex. The vertex set of $\widehat H$ is
\[
\widehat{\mathscr C}=\mathscr C\sqcup
\{r_{c,P}\,:\,w_c(P)>0\}.
\]
Define a dimension vector $\widehat v$ on $\widehat H$ by
\[
\widehat v_c=\sum_P v_c(P)=\dim V_c,\quad
\widehat v_{r_{c,P}}=w_c(P)=\dim W_c(P),
\]
and define the decreasing height sequence
$$\mathbf i=(i_1,\ldots,i_N)$$ as in Lemma~\ref{lem:stable-irreducible},
with the same enlarged interval and zero multiplicity colour steps.
We attach the multiplicity $v_c(P)$ to each colour occurrence $(c,P)$,
and $w_c(P)$ to each framing occurrence $r_{c,P}$.
In this order they form a sequence 
$$\mathbf a=(a_1,\ldots,a_N)\in\mathbb N^N.$$
Between consecutive occurrences of a colour $c$, at positions $P+p_c$
and $P$, the entries whose vertices are adjacent to $c$ are the entry $r_{c,P}$,
when present, and the entries $(c',Q)$ with $c'$ adjacent to $c$ in $H$ and
$P<Q<P+p_c$. Thus
\eqref{eq:ns-word-inequality} for this pair is
\eqref{eq:ns-interval-dominance}. Each framing vertex occurs only once, with multiplicity $w_c(P)$,
so imposes no further inequality.
Apply Lemma \ref{lem:ns-finite-word}, orienting the colour graph by the
chosen $Q$ and each new edge from the framing towards its colour vertex. 
We obtain a filtred $\Pi(\widehat H)$-module $M$ as in
\eqref{eq:ns-word-filtration}.

At the colour vertices, the filtration defines a $\bbZ$-flag $F=(F_c^{\leqslant P})_{c,P}$
with
\[
\dim F_c^{\leqslant P}/F_c^{<P}=v_c(P).
\]
Identify its colour spaces with $\dV_c$ so that $F=F_\gamma$,
and the framing vertex spaces with $W_c(P)$.
The maps of $M$ along the doubled colour edges define $B=(B_h)$ and $B^*=(B_h^*)$.
Those along the framing edges define
$A_{c,P}:W_c(P)\to\dV_c$ and $A^*_{c,P}:\dV_c\to W_c(P)$.
Using the orientation signs fixed above, combine the framing maps
into $A=(A_c)$ and $A^*=(A_c^*)$, where
$W_c=\bigoplus_P W_c(P)$. 
We have thus defined a point
\[
\xi=(x,y),\quad x=(A,B)\in\dN,\quad
 y=(A^*,B^*)\in\dN^*.
\]
We claim that $\xi\in\Lambda^{\rm st}_{F_\gamma}.$
The theorem follows.

Now, we prove the claim.
Each term of \eqref{eq:ns-word-filtration} is a submodule, and each
successive quotient is supported at a single vertex. Since there are
no loops, each arrow maps the space at its source vertex into the
subspaces corresponding to later entries of the word.
Consequently, the maps on both directions of a colour edge strictly
lower position. The framing maps satisfy
\[
A_{c,P}:W_c(P)\longrightarrow F_c^{\leqslant P},\quad
A^*_{c,P}\bigl(F_c^{\leqslant P}\bigr)=0.
\]
These are precisely \eqref{eq:NF} and \eqref{eq:NFperp}.
On colour edges strict and weak position inequalities agree,
since adjacent colours have no common allowed position.
The preprojective relation at each colour vertex is, with the
orientation signs fixed above,
\[
\sum_h[B_h,B_h^*]+\sum_P A_{c,P}A^*_{c,P}=0.
\]
Therefore $\xi\in L_{F_\gamma}$ and $\mu(\xi)=0$, so
$\xi\in\Lambda_{F_\gamma}$.
The extra relations at the framing vertices cause no difficulty: we have
constructed a module satisfying all relations of $\Pi(\widehat H)$,
and only the colour-vertex relations are needed for $\mu=0$.

For adjacent vertices $c'$ and $c$ in $H$, let
$$T_{c'c}:M_{c'}\to M_c$$ be the linear map assigned by $M$ to
the arrow $c'\to c$ in the doubled quiver of $\widehat H$.
Via the chosen identifications, we
regard $T_{c'c}$ as a map $V_{c'}\to V_c$.
Fix an allowed position $P$ of $c$ with $v_c(P)>0$, and put
$R=P+p_c$. Consider the two consecutive occurrences of $c$ in the word,
at positions $R$ and $P$. The occurrence at $R$ is included even when
$v_c(R)=0$. We apply the operator $K_c$, defined in \eqref{K1}, at the occurrence
corresponding to $R$.

At this step $K_c$ replaces $F_c^{\leqslant R}$ by
$F_c^{\leqslant P}$. The neighbouring spaces are $F_{c'}^{<R}$,
and the remaining framing vertices have positions $Q\leqslant P$.
Hence, the definition of $K_c$ gives
\begin{equation}\label{eq:ns-top-flag}
F_c^{\leqslant P}
=A_c(W_c^{\leqslant P})+
\sum_{c'\sim c}T_{c'c}\bigl(F_{c'}^{<P+p_c}\bigr).
\end{equation}
This identity holds also at the highest nonzero black positions for which 
the multiplicity is zero.

By Lemma \ref{lem:words}, we have $F_\xi\subseteq F$ termwise.
We first show
\begin{equation}\label{eq:ns-graded-generation}
F_{\xi,c}^{\leqslant P}+F_c^{<P}=F_c^{\leqslant P}\quad\forall\,c,P.
\end{equation}
If this fails, choose a largest position $P$ and a colour $c$ at which it fails. 
Such a position exists since the
nonzero graded pieces have finite support. Necessarily $v_c(P)>0$.
For each neighbour $c'$, equality holds at each allowed position
strictly above $P$. Eliminating the corresponding graded pieces, from
the largest downwards, gives
\[
F_{c'}^{<P+p_c}=F_{\xi,c'}^{<P+p_c}+F_{c'}^{\leqslant P}.
\]
On the other hand, we have
\[
T_{c'c}\bigl(F_{\xi,c'}^{<P+p_c}\bigr)\subseteq F_{\xi,c}^{\leqslant P}
,\quad
T_{c'c}\bigl(F_{c'}^{\leqslant P}\bigr)\subseteq F_c^{<P},
\]
Finally, we have 
$$A_c(W_c^{\leqslant P})\subseteq F_{\xi,c}^{\leqslant P}.$$
Substituting these inclusions into \eqref{eq:ns-top-flag}, we obtain
\[
F_c^{\leqslant P}\subseteq F_{\xi,c}^{\leqslant P}+F_c^{<P}.
\]
This contradicts the chosen failure of \eqref{eq:ns-graded-generation}.
Increasing induction through the finitely many nonzero graded pieces
now gives $F_{\xi,c}^{\leqslant P}=F_c^{\leqslant P}$ for each $c,P$.
Thus 
$$F_\xi=F_\gamma.$$
We deduce that the point $\xi$ lies in
$\Lambda^{\rm st}_{F_\gamma}$, i.e., the claim is proved.
\end{proof}

\section{Quasimap modules}\label{sec:QMM}
\subsection{Based quasimap modules}
For a point $\bar\xi\in \M_Q(v,w)$, let 
$$\operatorname{QM}^{\bar\xi}_\beta(\mathbb P^1, \M_Q(v,w))$$ be the stack of quasimaps
$\mathbb P^1\to \M_Q(v,w)$ of degree $-\beta$ nonsingular at $\infty$
with value $\bar\xi$ there. Set
$$\operatorname{QM}^{\bar\xi}=\bigsqcup_\beta\operatorname{QM}^{\bar\xi}_\beta
,\quad
\operatorname{QM}^{\bar\xi}_\beta=\operatorname{QM}^{\bar\xi}_\beta(\mathbb P^1, \M_Q(v,w))$$
Botta-Tamagni \cite[Thm.~E]{BT26} define a Donaldson-Thomas perverse sheaf
$\varphi_{\operatorname{QM}^{\bar\xi}}$ on $\operatorname{QM}^{\bar\xi}$ and
construct a Coulomb branch action on the critical cohomology
\[M^{\bar\xi}=H_T^*\bigl(\operatorname{QM}^{\bar\xi},\varphi_{\operatorname{QM}^{\bar\xi}}\bigr).\]
Since this Coulomb action is only constructed for quiver varieties, the following results
and conjectures are restricted to symmetric types.
Let $T_\rho\cong\mathbb C^\times$ be the one parameter subgroup \eqref{TR}.
Let $u=\hbar|_{T_\rho}$ be the equivariant parameter. 
For each type $\gamma$ put
\begin{align}\label{JG}
j_\gamma:\M_{\dQ}^\bullet(v_\gamma,w_\rho)\hookrightarrow \M_Q(v,w)^{T_\rho}\hookrightarrow \M_Q(v,w).
\end{align}
For each point $\xi\in\Lambda_{F_\gamma}^{\rm st,\bullet}$ let
$\bar\xi$ denote its image in $\M_{\dQ}^\bullet(v_\gamma,w_\rho)$.
We consider the obvious inclusions
$$i_\xi:\{\xi\}\hookrightarrow Z_\gamma^{\rm st}
,\quad
i_{\bar\xi}:\{\bar\xi\}\hookrightarrow \M_{\dQ}^\bullet(v_\gamma,w_\rho)$$
Let $M^{\bar\xi}_\rho$ be the specialization of $M^{\bar\xi}$ at $s_\rho=1$ and $\hbar=1$.
Recall the notation \eqref{eq:microlocal-K}. We set 
$$n=\dim\dN=\dim Z_\gamma^{\rm st}
,\quad
\mathcal K_\eta=\muhom(\mathscr L_\gamma,\mathscr L_\eta)|_{Z_\gamma^{\rm st}}$$
The co-specialized microlocal $\bR$-module is given by
\[M^!(\rho,\gamma)^\xi=\bigoplus_\eta H^*\bigl((i_\xi)^!\mathcal K_\eta[2n]\bigr).\]
The point $\xi\in\Lambda_{F_\gamma}^{\rm st,\bullet}$ is contained in the
$T_\rho$-invariant linear Lagrangian  $L_{F_\gamma}$.
It is $P_\gamma$-invariant, but may not be $\dG$-invariant.
Hence, the point $\bar\xi$ may not be polarized in the sense of \cite[\S1.3.5]{BT26}.

\begin{Conjecture}\label{conj:BT-specialized}
Assume $D=1$. 
Choose $\xi\in\Lambda_{F_\gamma}^{\mathrm{st},\bullet}$. 
There is an  $\bR$-module isomorphism
\[
 M^{\bar\xi}_\rho\cong M^!(\rho,\gamma)^\xi.
\]
The constant quasimap class corresponds to the complex orientation class $1_\gamma$ in
\[
 e_\gamma M^!(\rho,\gamma)^\xi
 =H^0\bigl((i_\xi)^!\bbC_{Z_\gamma^{\rm st}}[2n]\bigr)\cong\bbC.
\]
The quasimap summand of degree $\beta$ corresponds to
\[
 \bigoplus_{\eta:\,\mathrm{ch}(\eta,\gamma)=\beta}
 e_\eta M^!(\rho,\gamma)^\xi.
\]
The comparison preserves the type pieces and quasimap degree.
\end{Conjecture}

\begin{Remark}
By Lemma~\(\ref{lem:free}\)(c), every point $\bar\xi$
of the \(T_\rho\)-fixed locus \(\M_Q(v,w)^{T_\rho}\) belongs to a unique graded stratum \(\M_{\dQ}^\bullet(v_\gamma,w_\rho)\) and admits a representative
 \(\xi\in\Lambda_{F_\gamma}^{\mathrm{st},\bullet}\). Thus the conjecture concerns all \(T_\rho\)-fixed evaluation points $\bar\xi$ in
 the construction of Botta-Tamagni. 
\end{Remark}

\subsection{The polarized case}
In this section we give a proof of Conjecture \ref{conj:BT-specialized} if the point $\bar\xi$ is polarized.
More precisely, we say that $\bar\xi$ is BT-polarized if there is
an orientation polarization $N^{\rm BT}\subset T^*\dN$ as in \cite[\S1.3.5]{BT26} such that $\bar\xi\in L^{\rm BT}$, where
$$L^{\rm BT}=\M_Q(v,w)\cap[N^{\rm BT}/\dG].$$
Then, we define the \emph{polarized based quasimap space} of degree $-\beta$ to be
\[\mathfrak Q_\beta^{\bar\xi}=\operatorname{QM}_{-\beta}^{\bar\xi}(\mathbb P^1,L^{\rm BT})\]
and set
$$\mathfrak Q^{\bar\xi}=\bigsqcup_\beta\mathfrak Q_\beta^{\bar\xi}.$$
The stack $\mathfrak Q^{\bar\xi}$
 is a Deligne-Mumford stack which parametrizes the maps $\mathbb P^1\to [N^{\rm BT}/\dG]$ nonsingular at $\infty$
with value $\bar\xi$ there.
By \cite[Props.~1.3.8, 2.4.8, 3.2.2]{BT26}, its
Borel-Moore homology realizes $M^{\bar\xi}$ as a module over the quantized Coulomb branch.

We define as in \eqref{X}-\eqref{pi} the subspace $N_F^{\rm BT}\subset N^{\rm BT}$ for any $\mathbb Z$-flag $F$, and the incidence map
$$\pi_\eta^{\rm BT}:\mathscr X_\eta^{\rm BT}\to N^{\rm BT}.$$
We abbreviate $N^{\rm BT,\, st}=(N^{\rm BT})^{\rm st}$ and
$$
\mathcal X_\eta=\mathscr X_\eta^{\rm BT,\, st}/\dG
,\quad
\mathcal N=N^{\rm BT,\, st}/\dG
,\quad
\mathscr X_\eta^{\rm BT,\, st}=(\pi_\eta^{\rm BT})^{-1}(N^{\rm BT,\, st}).$$
Thus, the map $\pi_\eta^{\rm BT}$ yields a proper map 
$$f_\eta:\mathcal X_\eta\to\mathcal N.$$
We identify the cotangent bundle of $N^{\rm BT}$ with $T^*\dN$. 
In this presentation we have 
$$\xi=(x,0),\quad x\in N^{\rm BT,\, st}.$$
Let $b=[x]$ be the image of $x$ in $\mathcal N$.
We consider the polarized inclusion
\[
 \jmath:\mathcal N\longrightarrow \M_Q(v,w),\quad
 b\longmapsto\bar\xi.
\]
The point $b$ is $T_\rho$-fixed.
Set
\[
 \calN_b^-=\{s\in\mathcal N:\lim_{t\to\infty}T_\rho(t)s=b\}.
\]
Put
\[
d_\eta=\dim\mathscr X_\eta
,\quad
 d_\eta^{\rm BT}=\dim\mathscr X_\eta^{\rm BT}
 ,\quad
 \delta_\eta^{\rm BT}
 =(d_\eta-d_\gamma)-(d_\eta^{\rm BT}-d_\gamma^{\rm BT}).
\]
This shift vanish when the BT polarization is the original one. 

\begin{Lemma}\label{BTP}
For every type $\eta$
there is a natural isomorphism
\begin{equation}\label{C4}
 \Phi_\eta:\ (i_\xi)^!\mathcal K_\eta[2n]\xrightarrow{\cong}
 R\Gamma\bigl(f_\eta^{-1}(\calN_b^-),
                 \omega_{f_\eta^{-1}(\calN_b^-)}\bigr)
                 [\delta_\eta^{\rm BT}].
\end{equation}
These isomorphisms intertwine the Yoneda composition with the 
Borel-Moore convolution action. The
latter action is transported to $\bR$ by partial Fourier transform.
\end{Lemma}

\begin{proof}
The incidence spaces in this proof are relative to the unfolded framed quiver $\dQ$, with gauge group $\dG$.
If the polarization is first chosen on $Q$, since $D=1$, it lifts to a polarization of $\dQ$ in the obvious way. 
This gives the polarization $N^{\rm BT}$ of $T^*\dN$. 
This polarization is gauge invariant, so any representative $\xi$ of a polarized $\bar\xi$ belongs to it.

The stability  convention in \cite[\S\S1.1.4, 1.3.5]{BT26} is by cogeneration. 
Under the duality
\[
 (B_h,B_h^*,A,A^*)\longmapsto
 ((B_h^*)^t,-B_h^t,(A^*)^t,-A^t)
\]
this cogeneration condition is identified with our generation condition. 
We transport the setting of \cite{BT26} into ours by this duality.

We first work in the BT orientation. Let
$\mathscr L_\eta^{\rm BT}$ and $\mathcal K_\eta^{\rm BT}$ denote the
incidence and microlocal complexes for that orientation.
We return to the original orientation at the end. 
The change preserves the set $F_\xi$.
Indeed, since $D=1$, adjacent colours have opposite parities, so at
a source position $P$ the target flag has
$F_d^{<P}=F_d^{\leqslant P}$. Thus the conditions for $B_h$ and $B_h^*$ in
\eqref{eq:Gdef} coincide. Reversing an edge exchanges these maps up
to sign and leaves the preserved flag condition unchanged. 
The homogeneity is also preserved, since both arrow directions lower degree by one.
In BT coordinates, we have
$$\xi=(x,0),\quad F_\xi=F_\gamma.$$
Consider the $\dG$-torsor
$q:N^{\rm BT,\, st}\to\mathcal N.$
We have a Cartesian square
\[
\begin{array}{ccc}
\mathscr X_\eta^{\rm BT,\, st}&\xrightarrow{\ q_\eta\ }&\mathcal X_\eta\\[2mm]
{\scriptstyle\pi_\eta^{\rm BT}}\downarrow&&\downarrow{\scriptstyle f_\eta}\\[2mm]
N^{\rm BT,\, st}&\xrightarrow{\ q\ }&\mathcal N .
\end{array}
\]
Set 
$\mathcal G_\eta=(f_\eta)_*\omega_{\mathcal X_\eta}$.
Since $f_\eta$ is proper, smooth base change gives
\begin{align}\label{GL}
q^!\mathcal G_\eta=q^!(f_\eta)_*\omega_{\mathcal X_\eta}
\cong
(\pi_\eta^{\rm BT})_*\,(q_\eta)^!\omega_{\mathcal X_\eta}
\cong(\pi_\eta^{\rm BT})_*\,\omega_{\mathscr X_\eta^{\rm BT,\, st}}
=\mathscr L_\eta^{\rm BT}\big|_{N^{\rm BT,\, st}}.
\end{align}
Set $g=\dim\dG$. Since $q$ is smooth of relative complex
dimension $g$, we have
\[
q^!\mathcal G_\eta=q^*\mathcal G_\eta[2g].
\]
Let $\mathcal N^0\subset\mathcal N^{T_\rho}$ be the component containing $b$,
and $\mathcal N^+$ and $\mathcal N^-$ be the attracting and repelling strata of $\mathcal N^0$.
The inclusions and limit maps are
$$j^\pm:\mathcal N^\pm\hookrightarrow\mathcal N
,\quad
\ell^\pm:\mathcal N^0\hookrightarrow \mathcal N^\pm
,\quad
p^\pm:\mathcal N^\pm\to \mathcal N^0$$ 
Thus, we have
$$(p^-)^{-1}(b)=\mathcal N_b^-\subset\calN^-.$$
Since $F_\xi=F_\gamma$, every flag in $(\pi_\gamma^{\rm BT})^{-1}(x)$
contains $F_\gamma$ by Lemma~\ref{lem:words}.
They are equal since both have type $\gamma$. 
Thus 
$$(\pi_\gamma^{\rm BT})^{-1}(x)=\{(F_\gamma,x)\}.$$
The embedding $(F,u)\mapsto(F,(u,0))$ identifies $d\pi_\gamma^{\rm BT}$
with the restriction of $dp_\gamma^{\rm BT}$, which is injective by
Lemma~\ref{lem:smooth}(c) in the BT orientation.
Properness therefore makes the map
$$\pi_\gamma^{\rm BT}:\mathscr \M_{\dQ}^\bullet(v_\gamma,w_\rho)^{\rm BT}\to N^{\rm BT}$$ 
a closed embedding near $(F_\gamma,x)$.

On the zero cotangent section, the compensated action of
Lemma~\ref{lem:attracting-stratum} has weights
\[
 P-Q-1\quad\text{on }B:V_c(P)\to V_d(Q),\quad
 L-P\quad\text{on }A:W_c(L)\to V_c(P).
\]
Its limit as $t\to0$ exists if and only if the negative weight terms
vanish. By parity, this is equivalent to $Q<P$ for arrows and
$P\leqslant L$ for framings, exactly the flag conditions \eqref{eq:NF}.
Near $x$, the weight-zero framed words still generate, so the limit
is stable and lies in $\mathcal N^0$. Conversely,
Lemma~\ref{lem:attracting-stratum} recovers the flag of each nearby
point attracted to $\mathcal N^0$. 
Hence, in a neighbourhood of \(x\), the image of \(\pi_\gamma^{\rm BT}\) 
coincides with \(q^{-1}(\mathcal N^+)\). 
Since \(\pi_\gamma^{\rm BT}\) is a closed embedding near its fibre over \(x\), 
it identifies its source locally with this attracting locus.

Choose a neighbourhood $U\subset N^{\rm BT,\, st}$ of $x$ on which the image of
$\pi_\gamma^{\rm BT}$
is the smooth closed subvariety
$$S=U\cap\widetilde{\mathcal N}^+.$$
Write $j_S:S\hookrightarrow U$.
Definition~\eqref{L} and Verdier duality give
$$\mathscr L_\gamma^{\rm BT}|_{U}\cong(j_S)_*\omega_S
,\quad
\mathbb D(\mathscr L_\gamma^{\rm BT}|_{U})\cong(j_S)_*\bbC_S
,\quad
\dim S=d_\gamma^{\rm BT}
,\quad
\omega_S=\bbC_S[2d_\gamma^{\rm BT}].$$
Both sheaf identifications are local at $x$, with extension
by zero along $j_S$.
Consider the zero section and the projection
$$z_0:N^{\rm BT,\, st}\hookrightarrow T^*N^{\rm BT,\, st}
,\quad
\pi:T^*N^{\rm BT,\, st}\to N^{\rm BT,\, st}.$$
The conicity and \cite[(2.6)]{KS14} imply that
\[
 (z_0)^!\muhom(A,B)[2n]
 \cong \pi_!\muhom(A,B)[2n]
 \cong\mathbb D A\otimes B.
\]
Near the point $\xi$, the sheaf $\muhom(\mathscr L_\gamma^{\rm BT},\mathscr L_\eta^{\rm BT})$
is supported on the stable conormal locus.
Taking the costalk  gives
\begin{align}\label{KBTG}
 (i_\xi)^!\mathcal K_\eta^{\rm BT}[2n]
 \cong (i_{x,\,\widetilde{\mathcal N}^+})^!(\mathscr L_\eta^{\rm BT}|_{\widetilde{\mathcal N}^+})
 \cong (i_{b,\,\mathcal N^+})^!(j^+)^*\mathcal G_\eta.
\end{align}
Here $i_{x,\,\widetilde{\mathcal N}^+}$ and $i_{b,\,\mathcal N^+}$ are the point inclusions.
The last isomorphism is \eqref{GL}.
By construction the sheaf $\mathcal G_\eta$ is $T_\rho$-equivariant.
Braden's theorem \cite{Bra03} gives
\begin{align}\label{BG}(\ell^+)^!(j^+)^*\mathcal G_\eta\cong (p^-)_*(j^-)^!\mathcal G_\eta.\end{align}
Put 
$$\mathcal X_\eta^-=f_\eta^{-1}(\mathcal N^-)
,\quad
r_\eta=p^-\circ(f_\eta|_{\mathcal X_\eta^-}):\mathcal X_\eta^-\to\mathcal N^0
,\quad
\iota_b:\{b\}\hookrightarrow\mathcal N^0.$$
Hence, we have
$$r_\eta^{-1}(b)=f_\eta^{-1}(\mathcal N_b^-)=:Y_\eta.$$ 
Thus, using \eqref{KBTG} and \eqref{BG},
the base change for $f_\eta$ yields the first isomorphism in
\begin{align}\label{KRG}
(i_\xi)^!\mathcal K_\eta^{\rm BT}[2n]
 \cong(\iota_b)^!(r_\eta)_*\omega_{\mathcal X_\eta^-}
 \cong R\Gamma\bigl(Y_\eta,
                    \omega_{Y_\eta}\bigr).
\end{align}
The second one is
\[
 (\iota_b)^!(r_\eta)_*\omega_{\mathcal X_\eta^-}
 \cong\mathbb D\bigl((\iota_b)^*(r_\eta)_!\bbC_{\mathcal X_\eta^-}\bigr)
 \cong\mathbb D R\Gamma_c(r_\eta^{-1}(b),\bbC)
 \cong R\Gamma(r_\eta^{-1}(b),\omega_{r_\eta^{-1}(b)}).
\]

We now check the  compatibility of the isomorphism \eqref{KRG} with the Coulomb actions. Put
\[
 Z_{\nu\eta}=\mathcal X_\nu\times_{\mathcal N}\mathcal X_\eta,
 \quad s_\eta:Z_{\nu\eta}\to\mathcal X_\eta,
 \quad s_\nu:Z_{\nu\eta}\to\mathcal X_\nu.
\]
Fix any equivariant Yoneda class
$$a:\mathscr L_\eta^{\rm BT}\to\mathscr L_\nu^{\rm BT}[k].$$
Under the isomorphism \eqref{GL}, a restriction to $N^{\rm BT,\,st}$ and a
descent along $q$ send $a$ to a morphism
$$\mathcal G_\eta\to\mathcal G_\nu[k].$$
Under base change and adjunction, it is identified with a morphism over $Z_{\nu\eta}$
\[
 c_a:(s_\eta)^*\omega_{\mathcal X_\eta}
             \longrightarrow (s_\nu)^!\omega_{\mathcal X_\nu}[k]
\]
Set 
$$Z_{\nu\eta,b}=Y_\nu\times_{\mathcal N_b^-}Y_\eta
,\quad
k_\eta:Y_\eta\hookrightarrow\mathcal X_\eta
,\quad
k_Z:Z_{\nu\eta,b}\hookrightarrow Z_{\nu\eta}.$$
Let $s_{\eta,b},$ $s_{\nu,b}$ be the base changes of $s_\eta,$ $s_\nu$.
The exceptional base change of $c_a$ is the composed map
\[
\begin{aligned}
 c_{a,b}:(s_{\eta,b})^*\omega_{Y_\eta}
 =(s_{\eta,b})^*(k_\eta)^!\omega_{\mathcal X_\eta}
 \longrightarrow (k_Z)^!(s_\eta)^*\omega_{\mathcal X_\eta}
 \xrightarrow{(k_Z)^!c_a} (k_Z)^!(s_\nu)^!\omega_{\mathcal X_\nu}[k]
 \cong (s_{\nu,b})^!\omega_{Y_\nu}[k].
\end{aligned}
\]
Applying the pullback unit for $s_{\eta,b}$, then $c_{a,b}$, and the proper
trace for $s_{\nu,b}$ gives the maps
\[
\begin{aligned}
 R\Gamma(Y_\eta,\omega_{Y_\eta})
 \longrightarrow R\Gamma(Z_{\nu\eta,b},(s_{\eta,b})^*\omega_{Y_\eta})
 \xrightarrow{c_{a,b}} R\Gamma(Z_{\nu\eta,b},(s_{\nu,b})^!\omega_{Y_\nu})[k]
 \longrightarrow R\Gamma(Y_\nu,\omega_{Y_\nu})[k].
\end{aligned}
\]
This is the Borel-Moore convolution operator. 
Applied to
$$Y_\lambda\times_{\mathcal N_b^-}Y_\nu\times_{\mathcal N_b^-}Y_\eta,$$
the same construction identify the composition of Yoneda classes with composition of convolution operators.

Now, we return to the original orientation.
Changing the orientation means dualizing some arrow spaces.
Let $r$ be the dimension of the
reversed arrow space and $a_\eta$ the dimension of its
flag compatible subspace. We have 
$$d_\eta^{\rm BT}-d_\eta=r-2a_\eta.$$
A Fourier transform, shifted by $r$, sends the
constant sheaf of this subspace to that of its annihilator with shift $r-2a_\eta$. 
It commutes with the proper pushforward.
Hence
\[
 \mathsf F(\mathscr L_\eta)
 \cong\mathscr L_\eta^{\rm BT}[d_\eta-d_\eta^{\rm BT}],
 \quad
 \mathscr L_\eta^{\rm BT}=(\pi_\eta^{\rm BT})_*
                      \omega_{\mathscr X_\eta^{\rm BT}}.
\]
Compatibility with $\muhom$ under the symplectic identification
\cite[\S3.7, Ex.~VII.2]{KS90} contributes the difference of the
shifts for $\eta$ and $\gamma$, namely $\delta_\eta^{\rm BT}$.
This proves \eqref{C4}.

To make the grading explicit, put
$\epsilon_\alpha=d_\alpha-d_\alpha^{\rm BT}$.
A degree-$k$ class from type $\eta$ to type $\nu$ becomes,
under the Fourier transform $\mathsf F$, a class of degree
$$k+\epsilon_\nu-\epsilon_\eta=k+\delta_\nu^{\rm BT}-\delta_\eta^{\rm BT}.$$
Hence, it yields a degree $k$ class between the shifted complexes in
\eqref{C4}. The Fourier transform preserves composition, so the
compatibility just proved gives
\[
 \Phi_\nu[k]\circ a_{\rm mic}
   =a_{\rm BM}\circ\Phi_\eta
\]
for every $a\in e_\nu\bR^k e_\eta$, with the transported action on
the right. 
\end{proof}

Choose $\xi\in\Lambda_{F_\gamma}^{\mathrm{st},\bullet}$ such that
$\bar\xi$ is BT-polarized. Let $x$, $b$ be as above. 
Let 
$$(\mathfrak Q_\beta^{\bar\xi})^{T_\rho}\subseteq\mathfrak Q_\beta^{\bar\xi}$$ be the fixed locus
for the diagonal action of $T_\rho(t)$ on the target
and $z\mapsto t^2z$ on the source $\mathbb P^1$.
We call it the $T_\rho$-fixed locus in $\mathfrak Q_\beta^{\bar\xi}$.
We need the following.

\begin{Lemma}\label{lem:fixed locus}
$(\mathfrak Q_\beta^{\bar\xi})^{T_\rho}
 \cong\bigsqcup_{\eta:\,|\gamma|-|\eta|=2\beta} f_\eta^{-1}(\calN_b^-).$
\end{Lemma}

\begin{proof}
For a quasimap $q$ in $(\mathfrak Q_\beta^{\bar\xi})^{T_\rho}$, let $\mathcal E_i$ be its associated
gauge bundles over $\mathbb P^1$.
Put 
$$D=\operatorname{Spec}\mathbb C\llb z\rrb
,\quad
D^\times=\operatorname{Spec}\mathbb C\llp z\rrp.$$
The set of sections
$\mathcal E_{i,0}=\Gamma(D,\mathcal E_i)$
is a $\mathbb C\llb z\rrb$-lattice in
$\Gamma(D^\times,\mathcal E_i)$.
Choose a homogeneous trivialization $\tau_i$ over
$\mathbb P^1\smallsetminus\{0\}$, compatible with the framing at $\infty$.
It identifies the completed generic fibre $\Gamma(D^\times,\mathcal E_i)$ with
$V_i\otimes\mathbb C\llp z\rrp$ and gives the embedded lattice
\[
 L_i(q,\tau_i)=\tau_i\bigl({\mathcal E}_{i,0}\bigr)
 \subset V_i\otimes\mathbb C\llp z\rrp.
\]
Different choices of trivialization give gauge equivalent embedded lattices. 
The equivariant structure makes
${\mathcal E}_{i,0}$ homogeneous: it admits a homogeneous
$\mathbb C\llb z\rrb$-basis $(e_{i,a})_{a=1}^{v_i}$ with
\[
 {\mathcal E}_{i,0}
 =\bigoplus_{a=1}^{v_i}\mathbb C\llb z\rrb e_{i,a},
 \quad \deg(z^ke_{i,a})=P_{i,a}+2k.
\]
The multiset of degrees $P_{i,a}$ is intrinsic to the equivariant
gauge bundle and defines a type $\eta$.
After applying the trivialization $\tau_i$, identify the degree $P$
subspace of $V_i\llp z\rrp$ with the colour space
$\dV_{(i,P\bmod 2)}$. The degree $P$ subspaces of
$L_i(q,\tau_i)$ then form an increasing $\bbZ$-flag
$F(q,\tau)\in\Fl_\eta$, as in Lemma~\ref{lem:semilinear-flags}.
This $\mathbb Z$-flag is only defined up to a change of frame.
We consider the component of the fixed locus with type $\eta$.
In the chosen homogeneous trivialization over
$\mathbb C^\times=\mathbb P^1\smallsetminus\{0,\infty\}$, let
$$B_{PQ}(1):V_i(P)\to V_j(Q)
,\quad
A_{LP}(1):W_i(L)\to V_i(P)$$ 
be the arrow and framing of the representation at $z=1$.
The fixedness equations for the  $T_\rho$-action,
see the proof of Lemma~\ref{lem:attracting-stratum},
assign to these maps the $z$-exponents $k_B,$ $k_A$ given by
\[
 2k_B=P-Q-1,\quad 2k_A=L-P.
\]
Since $D=1$, an arrow joins positions of opposite parity, whereas
$L$ and $P$ have the same parity for a framing.  Thus both exponents
are integers.  As $z\mapsto t^2z$, equivariance forces the unique
Laurent monomials
\[
 B_{PQ}(z)=z^{(P-Q-1)/2}B_{PQ}(1),\quad
 A_{LP}(z)=z^{(L-P)/2}A_{LP}(1).
\]
Consequently, the restriction of a fixed quasimap to
$\mathbb C^\times$ is uniquely determined, up to gauge, by its representation at $z=1$.
A fixed quasimap has no defect away from $0$ because the defect is invariant
under rotation of the source, and the quasimap is regular at $\infty$.
Hence it gives an equivariant family of stable representations on
$\mathbb P^1\smallsetminus\{0\}$, determined by its value at $z=1$. 
Modulo changes of frame, this representation at $z=1$ determins a point $s\in\calN$.
The prescribed value $\bar\xi$ at $\infty$ is the point $b\in\calN$, i.e., we have 
$$\lim_{t\to\infty}T_\rho(t)s=b.$$ 
Thus a quasimap $q$ in $(\mathfrak Q_\beta^{\bar\xi})^{T_\rho}$ yields a point $s\in \calN_b^-$ under the evaluation $z\mapsto 1$.

Conversely, a point \(s\in\mathcal N_b^-\) determines, by equivariance, a stable quasimap on \(\mathbb P^1\smallsetminus\{0\}\). 
To extend it across \(0\), choose a homogeneous \(\dG\)-lattice \(L\) in the completed generic fibre \(V\otimes\mathbb C\llp z\rrp\). 
Such a lattice specifies the extension of the gauge bundle to the formal disc \(D\). By Lemma \(\ref{lem:semilinear-flags}\), 
homogeneous lattices of type \(\eta\) correspond to flags \(F\in\Fl_\eta\). 
The quasimap section extends regularly across $0$ precisely when
its representative $x_s\in N^{\rm BT,\, st}$ belongs to $N_F^{\rm BT}$.

Hence, up to changes of frame, the set of fixed quasimaps is identified with the set of pairs
$$\{(F,x_s)\,:\,[x_s]\in \calN_b^-\,,\,x_s\in N^{\rm BT}_F\}.$$
Comparing determinant degrees with the standard lattice gives
$|\gamma|_i-|\eta|_i=2\beta_i$ for all $i$, hence
\[
 (\mathfrak Q_\beta^{\bar\xi})^{T_\rho}
 \cong\bigsqcup_{\eta:\,|\gamma|-|\eta|=2\beta} f_\eta^{-1}(\calN_b^-).
\]
\end{proof}

We can now prove the following comparison theorem.

\begin{Theorem}\label{thm:BT-polarized}
Assume $D=1$ and choose
$\xi\in\Lambda_{F_\gamma}^{\mathrm{st},\bullet}$ such that
$\bar\xi$ is BT-polarized. Then Conjecture~$\ref{conj:BT-specialized}$
holds at $\xi$, with the action transported to the original $\bR$
by partial Fourier transform. The isomorphism is ungraded and
preserves the type pieces, quasimap degree and distinguished class.
\end{Theorem}

\begin{proof}
The comparison of the Coulomb and Hecke actions in the polarized
Borel-Moore model is supplied by \cite[Props.~2.4.8, 3.2.2]{BT26}.
The first proposition identifies the Coulomb action with the BFN
correspondence action after dimensional reduction. The second one identifies
this action with the Hecke action. Under the fixed locus identification in Lemma \ref{lem:fixed locus}, 
lattice modifications become flag changing
correspondences, and the tautological classes become diagonal classes. 
Thus, after inverting the equivariant parameter $u$, the BFN action agrees with flag
convolution in the chosen polarization. Transport this action to the
original $\bR$ by partial Fourier transform.
After specialization at $u=1$ we forget cohomological degrees;
the virtual shifts and $\delta_\eta^{\rm BT}$ only regrade the type
pieces. Hence
\[
 M^{\bar\xi}_\rho
 \cong\bigoplus_\eta H^*\!\left(
 R\Gamma\bigl(f_\eta^{-1}(\calN_b^-),
       \omega_{f_\eta^{-1}(\calN_b^-)}\bigr)
       [\delta_\eta^{\rm BT}]\right)
\]
as ungraded $\bR$-modules, with the flag action on the right hand side.
By the comparison \eqref{C4} proved above, the maps
$H^*(\Phi_\eta^{-1})$ identify this module with
\[
 \bigoplus_\eta H^*((i_\xi)^!\mathcal K_\eta[2n])
 =M^!(\rho,\gamma)^\xi.
\]
Both comparisons preserve the type pieces and the quasimap degree
decomposition. The relative Fourier shift $\delta_\eta^{\rm BT}$
changes the grading of the $\eta$-piece but satisfies
$\delta_\gamma^{\rm BT}=0$, so the normalized $\gamma$-piece
remains in degree zero. The standard lattice is
identified with the constant quasimap. Since the $\gamma$-piece
is one dimensional, we may rescale
the family $\Phi$ by a common nonzero scalar so that
$\Phi_\gamma$ sends $1_\gamma$ to the constant quasimap class.
This preserves all intertwining relations, and the resulting module
isomorphism sends the constant quasimap class to $1_\gamma$.
\end{proof}

\begin{Example}\label{ex:non-BT-point}
Let $Q$ be $1\to2$, with reds at positions $0$ at vertex $1$,
and $-1$ at vertex $2$, and blacks at positions $0,-2$ at vertex $1$,
and $-1,-1$ at vertex $2$.
Write 
$$V_1=\langle e_0,e_{-2}\rangle
,\quad
V_2=\langle f_1,f_2\rangle.$$
Choose
\[
 A_1(1)=e_0,\quad A_2(1)=f_2,\quad
 B(e_0)=f_1,\quad B^*(f_2)=e_{-2},
\]
with all other components zero. Since $B^*B=BB^*=0$, the
moment map equation holds. The four vectors
$e_0,f_2,Be_0,B^*f_2$  span $V_1\oplus V_2$. Hence
$\xi\in\Lambda_{F_\gamma}^{\rm st,\bullet}$ and its class $\bar\xi\in\mathfrak L_{\rho,\gamma}$.
Nevertheless $B\ne0$ and $B^*\ne0$ on the same edge.
Thus $\bar\xi$ is not BT-polarized.
\end{Example}

\begin{Remark}
On the finite corner indexed by $\Gamma$, at any sufficiently generic point we have
$$(i_\xi)^!\mathcal K_\eta[2n]\cong (i_\xi)^*\mathcal K_\eta$$
At BT-polarized points
satisfying this genericity condition.
Theorem~\ref{thm:BT-polarized}
therefore gives the characteristic cycle multiplicities of
Theorem~\ref{thm:microlocal-JH} and the ungraded Morse classes of
Corollary~\ref{cor:microlocal-morse-K0}. At special BT-polarized
points it instead gives costalk multiplicities. 
\end{Remark}

\subsection{Non based quasimap modules}
Now, we consider a relative version. 
By \cite[\S0.1.6.4]{BT26} we have the evaluation morphism
\[
 \operatorname{ev}_{\beta,\infty}:\operatorname{QM}_{-\beta}^{\mathrm{ns},\infty}(\mathbb P^1,\M_Q(v,w))
 \longrightarrow \M_Q(v,w)
\]
We set
$$d=\dim \mathfrak L_{\rho,\gamma}=n-\dim\dG
,\quad
m=\dim \M_{\dQ}^\bullet(v_\gamma,w_\rho).$$
Let $s_\gamma:\M_{\dQ}^\bullet(v_\gamma,w_\rho)\hookrightarrow \mathfrak L_{\rho,\gamma}$ be the fixed locus embedding.
The complex $\mathcal K_\eta$ descends to $\mathcal F_\eta$ on $\mathfrak L_{\rho,\gamma}$.
Set
\[
 \mathcal D_\eta=(s_\gamma)^!\mathcal F_\eta[2(d-m)].
\]
The transitivity of exceptional pullbacks gives
\[
 H^*\bigl((i_{\bar\xi})^!\mathcal D_\eta[2m]\bigr)
 \cong e_\eta M^!(\rho,\gamma)^\xi,
 \quad \mathcal D_\gamma\cong\bbC_{\M_{\dQ}^\bullet(v_\gamma,w_\rho)}.
\]
Recall the charge in \eqref{CEG} and the map $j_\gamma$ in \eqref{JG}.

\begin{Conjecture}\label{conj:BT-Microlocal}
\label{conj:BT-relative}Assume $D=1$ and \eqref{bfE}. \hfill
\begin{enumerate}[label=$\mathrm{(\alph*)}$,leftmargin=8mm]
\item
The relative quasimap construction gives $T_\rho$-equivariant complexes $\mathscr K_\beta\in D^b(\M_Q(v,w))$ 
carrying a Coulomb branch action such that $\mathscr K_0\cong\bbC_{\M_Q(v,w)}$.
At ${\rm BT}$-polarized evaluation points $\bar\xi\in \M_Q(v,w)$, the exceptional base change
identifies their costalk cohomology with the based quasimap module.
\item
After equivariant localization and specialization at $u=1$ there
is an isomorphism of $\bR$-modules
\begin{equation}\label{eq:BT-relative-microlocal}
 \bigoplus_{\beta\in\bbZ_{\geqslant0}^I}
 H^*_{T_\rho}\!\left(\M_{\dQ}^\bullet(v_\gamma,w_\rho),(j_\gamma)^*\mathscr K_\beta\right)
 [u^{-1}]\Big|_{u=1}
 \ \cong\ \bigoplus_\eta H^*(\M_{\dQ}^\bullet(v_\gamma,w_\rho),\mathcal D_\eta).
\end{equation}
The degree $\beta$ summand corresponds to the types $\eta$ with
$\mathrm{ch}(\eta,\gamma)=\beta$.
The degree zero identification is $H^*(\M_{\dQ}^\bullet(v_\gamma,w_\rho),\bbC)$ on both sides.
\item
For $\xi$ as in Conjecture~$\ref{conj:BT-specialized}$,
there is an isomorphism of $\bR$-modules
\[
 \bigoplus_{\beta\in\bbZ_{\geqslant0}^I}
 H^*_{T_\rho}\!\bigl(\{\bar\xi\},
 (i_{\bar\xi})^!(j_\gamma)^*\mathscr K_\beta[2m]\bigr)
 [u^{-1}]\Big|_{u=1}
 \ \cong\ M^!(\rho,\gamma)^\xi.
\]
It takes the type pieces to the quasimap degrees, and the constant quasimap class to $1_\gamma$.
\end{enumerate}
\end{Conjecture}

At BT-polarized points, part (c) follows from (a) and Theorem~\ref{thm:BT-polarized}.

\appendix
\section{The non symmetric global comparison map}
\label{app:ns-global-technical}
We  give detailed constructions and proofs underlying
the non symmetric global equivariant comparison map summarized in \S\ref{subsec:ns-modules}.
Let $(C,D)$ be of finite type. 

\subsection{The target}\label{app:target}
Put
\[
 N_0=\operatorname{lcm}(d_i)_{i\in I},\quad N=kN_0,\quad
 m_i=N/d_i,\quad m_{ij}=N/(d_i|c_{ij}|).
\]
The last integer is used only for edges $i-j$. 
Since $C$ is of finite type, we have
$$d_i|c_{ij}|=\max(d_i,d_j).$$ 
Hence $N_0$ is divisible by each 
integer $d_i|c_{ij}|$. The GLS symmetrizer for $C^{\rm T}$ is
$D_N=\operatorname{diag}(m_i)$. 
Set
\begin{align*}
H_i&=\bbC[\varepsilon_i]/(\varepsilon_i^{m_i}),\quad
 H_{ij}=\bbC[\varepsilon_{ij}]/(\varepsilon_{ij}^{m_{ij}}),
 \quad \varepsilon_{ij}=\varepsilon_i^{|c_{ij}|}
                         =\varepsilon_j^{|c_{ji}|},\\
 R_N&=\bigoplus_i\Hom_{H_i}(H_i\otimes W_i,H_i^{v_i})
 \ \oplus\!\bigoplus_{h:i\to j}
       \Hom_{H_{ij}}(H_i^{v_i},H_j^{v_j}),\\
 G_N&=\prod_iGL_{v_i}(H_i).
\end{align*}
Use the Frobenius functionals extracting the top coefficients of the
$H_i$ and $H_{ij}$ to identify the complex dual of $R_N$ with the
reverse arrow representation space as in \cite{VV23}. The Hamiltonian moment map 
$\mu_N:T^*R_N\to\mathfrak g_N^*$
 is then the framed GLS relation, with the conventions
\cite[\S8]{GLS}, \cite[Prop.~2.4]{HJV26}.
A point is stable if the images of the framing generate
$\bigoplus_iH_i^{\oplus v_i}$ under the arrow maps and the $\varepsilon_i$.
Define
\begin{align}\label{MN}\mathfrak M_N=\mu_N^{-1}(0)^{\rm st}/G_N.
\end{align}
This definition does not presume that a truncated module lifts to a
lattice. The lifting assertion is proved below.

\begin{Proposition}\label{prop:MN}\hfill
\begin{enumerate}[label=$\mathrm{(\alph*)}$,leftmargin=8mm]
\item $\mathfrak M_N$ is a quasi-projective fine moduli space of stable
framed locally free $\Pi(C^{\rm T},D_N)$-modules of rank $v$.
The quotient map is a Zariski locally trivial principal $G_N$-bundle.
\item $\mathfrak M_N$ is smooth of dimension
$2\dim R_N-2\dim G_N$ when it is nonempty.
\end{enumerate}
\end{Proposition}

\begin{proof}
Put $d=\sum_i m_i v_i$, and let $E$ be the space of framed words
of length at most $d$, including the letters $\varepsilon_i$.
At a stable point these words span $M=\bigoplus_i H_i^{v_i}$:
otherwise their spans would stabilize before dimension $d$ and
form a proper submodule containing the framing. Hence $E\twoheadrightarrow M$.
Its kernel and the outgoing framing determine the module.
The relations define a closed condition on such quotients, while
$H_i$-freeness is open. This gives a quasi-projective fine moduli
space inside a Grassmannian times an affine space. Its tautological
$H_i$-free modules have a Zariski locally trivial bundle of
$H_i$-frames, namely $\mu_N^{-1}(0)^{\rm st}$. This proves (a).

A stabilizer fixes the framing and commutes with every letter,
so it fixes all framed words. It is trivial, also infinitesimally.
The Hamiltonian identity
\[
 \langle d\mu_{N,p}(v),\xi\rangle=\omega(\xi\cdot p,v)
\]
then makes $d\mu_N$ surjective at stable points. Thus the stable
zero fibre is smooth of dimension $2\dim R_N-\dim G_N$.
Its quotient by the free $G_N$-action is smooth of dimension
$2\dim R_N-2\dim G_N$, proving (b).
\end{proof}

\subsection{The tower of truncations}
Reduction modulo $z_i^{kN_0/d_i}$ gives the following tower.
Part (a) does not require a torus action.

\begin{Proposition}\label{prop:ns-truncation-tower}
\hfill
\begin{enumerate}[label=$\mathrm{(\alph*)}$,leftmargin=8mm]
\item
For each $k\geqslant1$ the reduction defines a map
\[
 r_k:\mathfrak M_{(k+1)N_0}\longrightarrow\mathfrak M_{kN_0}
\]
It is a torsor under a vector bundle. 
It gives an isomorphism on cohomology.
\item
Suppose that an algebraic torus $A$ acts compatibly on the quiver
data defining the varieties $\mathfrak M_{kN_0}$ for all $k\geqslant1$. 
Then the reduction maps $r_k$ are
$A$-equivariant and induce isomorphisms on
$A$-equivariant cohomology. For each subtorus $T\subseteq A$, the map
\[
 \mathfrak M_{(k+1)N_0}^{T}\longrightarrow\mathfrak M_{kN_0}^{T}
\]
is an affine bundle over each connected component of its
target and gives an isomorphism of the
cohomology. The flavour and loop-rotation actions to which this
applies are constructed in \S\ref{subsec:ns-lattice-action}.
\end{enumerate}
\end{Proposition}

\begin{proof}
Put $\tau=z^{N_0}$, $a_i=N_0/d_i$ and
$a_{ij}=N_0/(d_i|c_{ij}|)$. Consider the finite free
$\bbC[\tau]$-algebras
\[
 \mathcal H_i=\bbC[\tau,u_i]/(u_i^{a_i}-\tau),\quad
 \mathcal H_{ij}=\bbC[\tau,u_{ij}]/(u_{ij}^{a_{ij}}-\tau),
\]
where $u_{ij}=u_i^{|c_{ij}|}=u_j^{|c_{ji}|}$.
Using their top-coefficient Frobenius pairings, Proposition~\ref{prop:MN}
gives a family $\mathcal X\to\mathbb A^1_\tau$ with central fibre
$\mathfrak M_{N_0}$. Framed generation kills infinitesimal
stabilizers, so the relative moment map is submersive and
$\mathcal X$ is smooth over $\mathbb A^1_\tau$.

Set $R_k=\bbC[\tau]/(\tau^k)$. Since
$\mathcal H_i/(\tau^k)=\bbC[u_i]/(u_i^{ka_i})$,
$\mathfrak M_{kN_0}$ parametrizes $R_k$-sections of $\mathcal X$:
stability is detected modulo $\tau$ by Nakayama, and the moment
equations use the relative pairing followed by $[\tau^{k-1}]$.
By smoothness, lifts of a section to $R_{k+1}$ form a torsor under
\[
 q_k^*T_{\mathfrak M_{N_0}}\otimes I_k,
 \quad I_k=(\tau^k)/(\tau^{k+1}),
\]
where $q_k$ is reduction to $\mathfrak M_{N_0}$
\cite[Tag~02GZ]{Stacks}. Such torsors induce cohomology
isomorphisms, proving (a).

For (b), these torsors are $A$-equivariant. Over a $T$-fixed
base point, linear reductivity makes the fixed locus of the
affine $T$-action a nonempty affine space directed by the
invariant translation space. These spaces form a vector bundle
on each connected fixed component. Thus $r_k^T$ is an affine
bundle, and the same argument on Borel constructions gives
the equivariant cohomology isomorphisms.
\end{proof}

\subsection{Fixed points and tautological classes}
\label{subsec:ns-lattice-action}
The graded truncation model has the following equivariant
compatibilities which follow from the construction. Put
\[
 T_\rho(t)=(\rho(t^{-1}),t^{-1}).
\]
\begin{enumerate}[label=$\mathrm{(\alph*)}$,leftmargin=8mm]
\item The $F\times\bbC^\times$-action on the arrow, framing, dual-arrow and
gauge spaces extends to the family over $\mathbb A^1_\tau$ used in
Proposition~\ref{prop:ns-truncation-tower}. It preserves stability
and the zero locus of the moment map, and rotation acts by
\begin{equation}\label{eq:ns-action-normalization}
 \lambda\cdot z_i=\lambda^{2d_i}z_i,\quad
 \lambda\cdot\tau=\lambda^{2N_0}\tau.
\end{equation}
The dual-arrow action is transported through the Frobenius pairing.

\item For a homogeneous lattice $L$ with associated flag $F_L$, the
$s$-fixed parts of the loop spaces identify with the finite-dimensional
spaces of \S\ref{subsec:ns-modules}:
\begin{equation}\label{eq:ns-fixed-linear-data}
 \begin{gathered}
 G_K^s\cong\dG,\quad N_K^s\cong\dN,\quad (N_K^*)^s\cong\dN^*,\quad
 N_\calO(L)^s\cong\dN_{F_L},\quad
 (N_\calO(L)^\perp)^s\cong\dN_{F_L}^\perp.
 \end{gathered}
\end{equation}
These identifications are compatible with the gauge action, pairings,
moment equations, families and changes of homogeneous frame.

\item Let $\mathcal V_i$ be the tautological $H_i$-module on
$\mathfrak M_{N_0}$ and put
\begin{align}\label{EI}
 \mathcal E_i=\mathcal V_i/\varepsilon_i\mathcal V_i.
\end{align}
It has rank $v_i$. On a fixed family of type $\gamma$, its
centred degree $P$ summand is $\calV_i(P)$. With the compensated
$T_\rho$-action this summand has weight $P$, as in
Lemma~\ref{lem:fixed}. Write
\begin{align}\label{a}H^*_{T_\rho}=\bbC[a]\end{align}
and denote the Chern roots of $\calV_i(P)$ by $u_{i,P,r}$.
The equivariant splitting principle gives
\begin{equation}\label{eq:ns-tautological-normalization}
 c_{T_\rho}(\mathcal E_i)
   =\prod_P\prod_{r=1}^{\dim V_i(P)}(1+Pa+u_{i,P,r}).
\end{equation}
This is understood after restriction to the fixed family, or a
sufficiently large truncation. At $a=1$ its roots are $P+u_{i,P,r}$.
\end{enumerate}

Indeed, for (a), use the relative algebras $\mathcal H_i$ and
$\mathcal H_{ij}$ from the proof of
Proposition~\ref{prop:ns-truncation-tower}, and give their generators
degrees
\[
 \deg u_i=2d_i,\quad
 \deg u_{ij}=2d_i|c_{ij}|,\quad \deg\tau=2N_0.
\]
The relations 
$$u_i^{N_0/d_i}=u_{ij}^{N_0/(d_i|c_{ij}|)}=\tau$$ are homogeneous.
The Frobenius functionals extracting the coefficients of
$u_i^{N_0/d_i-1}$ and $u_{ij}^{N_0/(d_i|c_{ij}|)-1}$ have
homogeneous degrees $-2(N_0-d_i)$ and $-2(N_0-d_i|c_{ij}|)$. 
The gauge action and moment equations
are also homogeneous, and the framed generation is invariant.
Part (b) follows from Lemma~\ref{lem:ns-lattice}.
For (c), Proposition~\ref{prop:ns-modules}(a) gives, on a
homogeneous basis, $h_t(e_{i,P})=t^Pe_{i,P}$ and
$h_t(z_ie)=t^{2d_i}z_ih_t(e)$. Therefore
$\mathcal E_i=\mathcal V_i/\varepsilon_i\mathcal V_i$ has
centred-degree-$P$ summand $\calV_i(P)$, on which $T_\rho$
acts with weight $P$. Its equivariant Chern roots are
$Pa+u_{i,P,r}$, proving \eqref{eq:ns-tautological-normalization}.

By these identifications, Proposition~\ref{prop:ns-modules}(b)
identifies the $T_\rho$-fixed locus  of type $\gamma$ for  $\mathfrak M_{N_0}$ with $\mathfrak L_{\rho,\gamma}$. 
In the symmetric case, this should not be confused with
Lemma~\ref{lem:attracting-stratum} which identifies
$\mathfrak L_{\rho,\gamma}$ with a $T_\rho$-{attracting stratum}
in the unfolded quiver variety. The type records the centred degrees of the generators
$L_i/z_iL_i$. For a fixed type, sufficiently large truncations
determine the homogeneous coefficients and changes of frame.
Therefore the  affine bundle invariance of Proposition~\ref{prop:ns-truncation-tower}
and Lemma~\ref{F} yield
\begin{equation}\label{eq:ns-fixed-cohomology}
 \begin{aligned}
 H^*((\mathfrak M_{N_0})^{T_\rho},\bbC)
 \cong\bigoplus_\gamma H^*(\mathfrak L_{\rho,\gamma},\bbC).
 \end{aligned}
\end{equation}
Here $\rho$ is fixed and $\gamma$ runs over the occurring types.

The $F\times\bbC^\times$-equivariant bundles $\mathcal E_i$ define the
Cartan Kirwan map
\begin{equation}\label{eq:ns-global-kirwan}
 \kappa^{\rm ns}:R_{F,\hbar}\longrightarrow
 H^*_{F\times\bbC^\times}(\mathfrak M_{N_0},\bbC),
\end{equation}
by sending symmetric polynomials to the corresponding
Chern classes of $\mathcal E_i$, and flavour and rotation parameters
to their equivariant classes. By
\eqref{eq:ns-tautological-normalization}, the restriction to the
$T_\rho$-fixed component gives the local Kirwan map in \eqref{eq:ns-local-kirwan}
\[
 \kappa^{\rm ns}_\gamma:\widehat R_\gamma\longrightarrow
 H^*(\mathfrak L_{\rho,\gamma},\bbC).
\]

\subsection{From local to global comparison maps}
We descend the global Kirwan map $\kappa^{\rm ns}$ to the
homological B-algebra \eqref{BTFH}. The Cartan inclusion gives the
surjective graded quotient
\begin{align}\label{CCQH}
 q:R_{F,\hbar}\twoheadrightarrow\B^\theta_{F,\hbar}.
\end{align}
For integral $\rho$ let $s_\rho=(\rho,1)$ in
$\operatorname{Spec}H^*_{F\times\bbC^\times}$.
Specialization at $s_\rho$, followed by completion at the gauge
character $\gamma$, gives
\[
 \B^\theta_{\rho,\gamma}
 =\bigl(\B^\theta_{F,\hbar}\otimes_{H^*_{F\times\bbC^\times}}\bbC_{s_\rho}\bigr)^\wedge_\gamma,
 \quad q_{\rho,\gamma}:\widehat R_\gamma\to\B^\theta_{\rho,\gamma}.
\]
The completed-localization comparison with the K-theoretic
Cartan quotient $q_\gamma:\widehat R_\gamma\to\widehat\B_\gamma$ is
\begin{equation}\label{eq:ns-global-completion}
 \begin{gathered}
 \nu_{\rho,\gamma}:\B^\theta_{\rho,\gamma}
    \xrightarrow{\cong}\widehat\B_\gamma,\quad
 \nu_{\rho,\gamma}\circ q_{\rho,\gamma}=q_\gamma.
 \end{gathered}
\end{equation}
Set $J=\Ker q$, and let ${\rm tor}\subset H^*_{F\times\bbC^\times}(\mathfrak M_{N_0},\bbC)$ 
be its ${H^*_{F\times\bbC^\times}}$-torsion
ideal, with 
$$H^*_{F\times\bbC^\times}(\mathfrak M_{N_0},\bbC)_{\rm tf}=H^*_{F\times\bbC^\times}(\mathfrak M_{N_0},\bbC)/{\rm tor}.$$
The following is proved as in the main text.

\begin{Proposition}\label{prop:ns-global-torsion}
For every finite type $(C,D)$, the following hold.
\begin{enumerate}[label=$\mathrm{(\alph*)}$,leftmargin=8mm,itemsep=2mm]
\item There is a unique map
\[
 \Phi^{\rm ns}_{\rm tf}:\B^\theta_{F,\hbar}\longrightarrow H^*_{F\times\bbC^\times}(\mathfrak M_{N_0},\bbC)_{\rm tf}
\]
whose image is the image of the equivariant tautological subalgebra.
\item If $H^*_{F\times\bbC^\times}(\mathfrak M_{N_0},\bbC)$ is ${H^*_{F\times\bbC^\times}}$-torsion free, there is a unique map
\[
 \Phi^{\rm ns}:\B^\theta_{F,\hbar}\longrightarrow H^*_{F\times\bbC^\times}(\mathfrak M_{N_0},\bbC)
 \text{ with }
 \Phi^{\rm ns}\circ q=\kappa^{\rm ns}.
\]
It is surjective exactly when $\kappa^{\rm ns}$ is.
\item At generic integral parameters $\mathrm{(a)}$ recovers the local comparison
maps after localization. 

\item
If $\mathfrak M_{N_0}$ is equivariantly formal,
the completion of the specialization of $\mathrm{(b)}$ at $\gamma$
agrees with $\Phi^{\rm ns}_\gamma$.
\end{enumerate}
\qed
\end{Proposition}

\begin{Remark}\hfill
\begin{enumerate}[label=$\mathrm{(\alph*)}$,leftmargin=8mm,itemsep=2mm]
\item 
By Theorem~\ref{thm:ns-map}(b), we have
\begin{equation*}
 \Phi^{\rm ns}_\gamma\text{ is onto}
 \quad\Longleftrightarrow\quad
 \kappa^{\rm ns}_\gamma\text{ is onto}.
\end{equation*}
A sufficient global criterion is surjectivity of
$H^*_{G_{N_0}}\to H^*(\mathfrak M_{N_0},\bbC)$.
The even degree source forces $H^{\rm odd}(\mathfrak M_{N_0},\mathbb C)=0$.
This criterion is not valid for all data:
for the $B_2$ datum $v=(1,2)$, $w=(0,3)$ of
Proposition~\ref{prop:B2-local-kirwan},
\eqref{eq:ns-fixed-cohomology} gives odd cohomology of
$(\mathfrak M_2)^{T_\rho}$. If $H^{\rm odd}(\mathfrak M_2,\mathbb C)$ vanished,
parity in the Borel spectral sequence and localization
would force the fixed point cohomology to be also even. 

\item
Purity of $H^j(\mathfrak M_{N_0},\mathbb C)$ of weight $j$ would imply equivariant
formality, by the Borel spectral sequence, which is important for surjectivity. A smooth
semiprojective action could establish purity, but one must prove
that scaling limits remain locally free. The purity theorem of
\cite[Thm.~B, \S5.4]{HJV25} uses a different stability.
Thus the purity results there do not settle the formality question.
\end{enumerate}
\end{Remark}

\section{Critical locus descriptions}\label{app:ns-critical}

\subsection{The critical locus description}\label{subsec:ns-crit}
Recall the notation in \eqref{LF}. We abbreviate
$$L_\gamma=L_{F_\gamma}.$$
It is a $P_\gamma$-stable Lagrangian subspace of $T^*\dN$. 
By Lemma \(\ref{lem:words}\), every \(\xi\in L_\gamma\) satisfies \(F_\xi\subseteq F_\gamma\). 
Therefore, by the definition of \(\Omega_\gamma\) in \eqref{OLZ}, we have
\[
L_\gamma\cap\Omega_\gamma=\{\xi\in L_\gamma:\ F_\xi=F_\gamma\}
\]
Recall that the stable locus defined is
\[\Lambda_{F_\gamma}^{\mathrm{st}}=L_\gamma\cap\Omega_\gamma\cap\mu^{-1}(0).\]
By Lemma \ref{lem:smooth}, on $L_\gamma$ the moment map takes values in 
$\frakp_\gamma^\perp=(\dg/\frakp_\gamma)^*$.
Let
\begin{align*}
f_\gamma&:(L_\gamma\cap\Omega_\gamma)\times\dg/\frakp_\gamma\to\bbC,\quad
f_\gamma(\xi,z)=\langle\mu(\xi),z\rangle=\operatorname{tr}\Big(z\big(\textstyle\sum_h[B_h,B_h^*]+AA^*\big)\Big).
\end{align*}
The set $L_\gamma\cap\Omega_\gamma$ is open in $L_\gamma$.
The function $f_\gamma$ is $P_\gamma$-invariant.
We define
$$S_\gamma=(L_\gamma\cap\Omega_\gamma)\times\dg/\frakp_\gamma,\quad T_\gamma=S_\gamma/P_\gamma$$
In the following proposition, let \(K^{\rm crit}_{P_\gamma}(S_\gamma,f_\gamma)\) denote the Grothendieck group of 
\(P_\gamma\)-equivariant coherent graded factorizations. The auxiliary grading torus acts by 
$$t\cdot(\xi,z)=(\xi,t^2z).$$ 
Thus the coordinates of \(\xi\) have degree \(0\), the multiplier \(z\) has degree \(2\), and the potential \(f_\gamma=\langle\mu(\xi),z\rangle\) 
is homogeneous of degree \(2\). A graded factorization has an odd differential of degree \(1\), whose square is multiplication by 
\(f_\gamma\).

\begin{Proposition}\label{prop:ns-crit}
Let $(C,D)$ be of finite type. Set $r=\dim\dg/\frakp_\gamma$. 
\begin{enumerate}[label=$\mathrm{(\alph*)}$,leftmargin=8mm]
\item The group $P_\gamma$ acts freely on $L_\gamma\cap\Omega_\gamma$.
The variety $T_\gamma$ is smooth quasi-projective.
It is a vector bundle over $(L_\gamma\cap\Omega_\gamma)/P_\gamma$.
\item 
 $\operatorname{Crit}(f_\gamma)/P_\gamma=\mathfrak L_{\rho,\gamma}$. It is smooth.
\item 
We have $\varphi_{f_\gamma}\bbC_{T_\gamma}\cong\bbC_{\mathfrak L_{\rho,\gamma}}[1-2r]$, hence
\begin{align*}
 H_c^k(T_\gamma,\varphi_{f_\gamma}\bbC)\cong H_c^{k+1-2r}({\mathfrak L_{\rho,\gamma}},\bbC),\quad
 H^k(T_\gamma,\varphi_{f_\gamma}\bbC)\cong H^{k+1-2r}({\mathfrak L_{\rho,\gamma}},\bbC).
\end{align*}
\item
The dimensional reduction gives
\[
 K^{\rm crit}_{P_\gamma}(S_\gamma,f_\gamma)
\cong K^{P_\gamma}(\Lambda^{\rm st}_{F_\gamma})\cong K(\mathfrak L_{\rho,\gamma}).
\]
\end{enumerate}
\end{Proposition}

\begin{proof}
(a) By Lemma~\ref{lem:words}, the equality $F_\xi=F_\gamma$ implies generation
by framed words, without imposing $\mu=0$. The stabilizer and
infinitesimal stabilizer therefore vanish as in Lemma~\ref{lem:free}(a).
Apply the quotient argument in Lemma~\ref{F} to the generated locus
$Y\subset T^*\dN$, omitting the moment equation.
By the uniqueness of the minimal flag as in Lemma~\ref{lem:Z},
the map $$\dG\times_{P_\gamma}(L_\gamma\cap\Omega_\gamma)\to Y$$ identifies
$(L_\gamma\cap\Omega_\gamma)/P_\gamma$ with a locally closed subvariety of $Y/\dG$.
Since $L_\gamma\cap\Omega_\gamma$ is open in $L_\gamma$, its quotient is
smooth, and $T_\gamma$ is its associated vector bundle.

(b) Differentiating gives
\[
 df_\gamma(\delta x,\delta y,\delta z)
 =\langle\mu(\xi),\delta z\rangle
 -\langle\delta x,z\cdot y\rangle
 +\langle z\cdot x,\delta y\rangle.
\]
The critical equations say $\mu(\xi)=0$,
$z\cdot x\in\dN_{F_\gamma}$ and
$z\cdot y\in\dN_{F_\gamma}^{\perp}$.
By Lemma~\ref{lem:IU}, the bundle map
\[
 \ell_\xi:\dg/\frakp_\gamma\longrightarrow T^*\dN/L_\gamma,
 \quad z\longmapsto z\cdot\xi
\]
is injective on $L_\gamma\cap\Omega_\gamma$. Locally it has a left inverse, so these
linear equations generate the ideal $z=0$. The remaining equation
cuts out $\Lambda_{F_\gamma}^{\rm st}$, smooth by
Lemma~\ref{lem:smooth}, proving the scheme theoretic assertion.

(c) By Lemma~\ref{lem:smooth}, the $r$ components of $\mu$
are normal coordinates $x_a$ near its zero set. With dual
multiplier coordinates $z_a$, we have
$f_\gamma=\sum_{a=1}^r x_a z_a$ locally.
Its Milnor fibre has reduced cohomology $\bbC$ in degree
$2r-1$, and the complex orientation makes the corresponding local
system constant. Thus, for unshifted vanishing cycles, we have
$$\varphi_{f_\gamma}\bbC_{T_\gamma}\cong
\bbC_{\mathfrak L_{\rho,\gamma}}[1-2r]$$
giving both cohomology formulas. 

(d) Dimensional reduction
gives the K-theory isomorphism \cite[\S2.1, Prop.~4.6]{VV23}.
\end{proof}

\subsection{Comparison with the critical cohomology of \cite{VV23}}\label{subsec:ns-VV23}
\subsubsection{Recollection from \cite{VV23}}
Let $\widetilde Q$ be the triple quiver, i.e.,  the double quiver $\bar Q$ with a loop $\varepsilon_i$ at each
vertex. It is framed by arrows $a_i:i\to i'$ and $a_i^*:i'\to i$. Its graded version $\widetilde Q^\bullet$ has gradings 
$$\deg\alpha_{ij}=b_{ij}=d_ic_{ij}
,\quad
\deg\varepsilon_i=2d_i
,\quad
\deg a_i=\deg a_i^*=-d_i$$
Let $\widetilde{\mathfrak M}^\bullet(v,w)$ be the moduli space of stable graded representations of $\widetilde Q^\bullet_f$ of graded 
dimension $v$, framed by $W$. The potentials are
\[
w_2=\sum_{i,j}o_{ij}\,\varepsilon_i^{-c_{ij}}\alpha_{ij}\alpha_{ji},\quad w_1=w_2+\sum_i\varepsilon_ia_i^*a_i,
\]
with $o_{ij}=-o_{ji}=\pm1$ given by an orientation.
Every term of $w_1$ has total degree zero for the grading above.
Write $w_1^\bullet$ for the same cyclic potential on the graded
framed quiver $\widetilde Q_f^\bullet$, and set
$\tilde f^\bullet_1=\operatorname{tr}(w_1^\bullet)$, see \cite{VV23} for details.

\begin{Lemma}[\cite{VV23}]\label{lem:VV23-crit}
A point of $\widetilde{\mathfrak M}^\bullet(v,w)$ is critical for $\tilde f^\bullet_1$ if and only if
\begin{enumerate}[label=$\mathrm{(\alph*)}$,leftmargin=8mm,itemsep=2mm]
\item $\varepsilon_i^{|c_{ij}|}\alpha_{ij}=\alpha_{ij}\varepsilon_j^{|c_{ji}|}$ for all arrows $\alpha_{ij}$ of $\bar Q$,
\item $\sum_jo_{ij}\sum_{k=0}^{|c_{ij}|-1}\varepsilon_i^{k}\alpha_{ij}\alpha_{ji}\varepsilon_i^{|c_{ij}|-1-k}+a_i^*a_i=0$ for all $i$,
\item $a_i\varepsilon_i=0$ and $\varepsilon_ia_i^*=0$.
\qed
\end{enumerate}
\end{Lemma}

By Lemma \ref{lem:VV23-crit}, 
the critical points of $\tilde f^\bullet_1$
are the stable graded framed modules over the Jacobian algebra of $(\widetilde Q,w_1)$. 
Up to the conditions (c), this is
the untruncated generalized preprojective algebra $\widetilde\Pi(C^{\rm T},D)$ when $\gcd(c_{ij},c_{ji})=1$ for all edges, e.g., in finite type. 
Then claim (a) is the relation $(H2)$ of \cite{GLS},
and claim (b) is the preprojective relation as in  \S\ref{subsec:ns-modules}. 
The lemma does not compare the critical loci
as moduli spaces. More precisely:

\begin{enumerate}[label=$\mathrm{(\alph*)}$,leftmargin=8mm]
\item
The set $\operatorname{Crit}(\tilde f_1^\bullet)$ parametrizes points in a {graded} triple quiver variety. 
The operators $\varepsilon_i$ are nilpotent. 
The set $\operatorname{Crit}(f_\gamma)/P_\gamma=\mathfrak L_{\rho,\gamma}$
parametrizes stable $s$-fixed lattices, with injective multiplication by $z_i$.

\item
The derivatives of $w_1^\bullet$ with respect to $\varepsilon_i$
give the preprojective relations. The arrow derivatives of $w_1^\bullet$
impose $\calO_{ij}$-linearity. In our lattice model this linearity is built in.
The costability of \cite{VV23} becomes our generation condition by
transposition. 
\item
An {ungraded} comparison is available after
truncating lattices modulo $z_i^{m_i}$,
using the deformed potential below. On the stratum with $v_i$
Jordan blocks of length $m_i$ one obtains
\eqref{eq:VV23-Jordan-stratum}. 
\end{enumerate}

\subsubsection{Nilpotent framing: a comparison of moduli spaces}
Following \cite{VV23},
for a fixed nilpotent endomorphism $\eta=(\eta_i)$ of the framing space $W'$, we consider the function
\begin{equation}\label{eq:VV23-deformed-potential}
\tilde f_\eta=\operatorname{tr}(w_2)
       +\sum_i\operatorname{tr}_{V'_i}(\varepsilon_i a_i^*a_i)
       -\sum_i\operatorname{tr}_{W'_i}(\eta_i a_i a_i^*).
\end{equation}
Let $\mathcal P_i(\alpha,\varepsilon)$ be the
preprojective expression in Lemma \ref{lem:VV23-crit}(b).
The critical equations of $\tilde f_\eta$ are
\begin{align}
 \varepsilon_i^{|c_{ij}|}\alpha_{ij}
     =\alpha_{ij}\varepsilon_j^{|c_{ji}|},\quad
 \mathcal P_i(\alpha,\varepsilon)+a_i^*a_i=0,\quad
 \varepsilon_i a_i^*=a_i^*\eta_i,
 \quad a_i\varepsilon_i=\eta_i a_i.
 \label{eq:VV23-deformed-critical}
\end{align}
In particular $\eta=0$ imposes annihilation
by the loops on the framing maps, whereas $\eta\neq 0$ imposes intertwining. 
Put
$$N=kN_0
,\quad
m_i=N/d_i
,\quad
V'_i=\bbC^{m_i v_i}
,\quad 
W'_i=H_i\otimes W_i
,\quad
\eta_i=\varepsilon_i|_{H_i\otimes W_i}.$$
Let $\mathcal O_{(m_i^{v_i})}$ be the nilpotent orbit with $v_i$
Jordan blocks of length $m_i$. Inside the stable triple quiver
variety with dimensions $\dim V'_i=m_i v_i$ and framing $W'$, take
\[
 Z_N=\operatorname{Crit}(\tilde f_\eta)
           \cap\{\varepsilon_i\in\mathcal O_{(m_i^{v_i})}
                        \text{ for each }i\}.
\]
The critical locus is computed in the ambient triple quiver variety, not in the orbit stratum.
Recall the variety $\mathfrak M_N$ introduced in \eqref{MN}.

\begin{Lemma}
There is an isomorphism of moduli schemes
\begin{equation}\label{eq:VV23-Jordan-stratum}
 \mathfrak M_N\ \cong\ Z_N,
\end{equation}
after transposing the stability convention and making the usual
opposite arrow sign choices. 
\end{Lemma}

\begin{proof}
Fix each
$\varepsilon_i$ in its prescribed Jordan normal form. Its centralizer
is $GL_{v_i}(H_i)$. The first equation of
\eqref{eq:VV23-deformed-critical} says that the arrows are
$H_{ij}$-linear. The last two say that the framing maps are
$H_i$-linear, and the middle equation is the framed GLS
moment map relation. Quotienting the orbit stratum by
$\prod_iGL(V'_i)$ is therefore the same as fixing these Jordan forms
and quotienting by $G_N$. Finally, the stability in \cite{VV23}
is that no nonzero $(\alpha,\varepsilon)$-submodule lies in
$\Ker a$. Under duality it becomes generation by the incoming
framing, which is our stability. The Frobenius self-duality of $H_i$
identifies the dual coefficient modules with $H_i$-free modules.
This identifies the equations and the moduli functors, proving
\eqref{eq:VV23-Jordan-stratum}.
\end{proof}

Thus local freeness has a precise meaning in the triple quiver model: it is a prescribed {rectangular Jordan type}, together
with a nilpotent framing of the corresponding lengths. 

\subsubsection{A counterexample in type $B_2$}
For symmetric type, recall the following result from \cite{VV23}.
\begin{Lemma}[\cite{VV23}]\label{lem:VV23-sym}
Let $D=1$. 
Up to shifts, we have
\[
H^*\big(\widetilde{\mathfrak M}^\bullet(v_\gamma,w_\rho),\varphi_{\tilde f^\bullet_1}\big)\ \cong\ H^*\big(\M_{\dQ}^\bullet(v_\gamma,w_\rho),\bbC\big)\ \cong\ H^*_{P_\gamma}(\Lambda^{\rm st}_{F_\gamma},\bbC)
\]
Here the stability conditions are matched by transposition.
\end{Lemma}

\begin{proof}
For $D=1$, the critical equations make $\varepsilon$ commute
with every doubled arrow and give $a\varepsilon=0$.
Thus $\operatorname{Im}\varepsilon$ is an unframed subrepresentation,
so costability forces $\varepsilon=0$.
The critical support consequently lies in the open set where the
doubled data, after forgetting loops, are stable. Only on this open
set are the loop variables a vector bundle; restriction to it
does not change vanishing-cycle cohomology.
There the potential
\[
 w_1=\sum_i\varepsilon_i\Bigl(\sum_j o_{ij}\alpha_{ij}\alpha_{ji}
                  +a_i^*a_i\Bigr)
\]
is fibre-linear with moment map section. The degrees agree with
\S\ref{subsec:graded} after shifting $W$ by $d_i$.
Dimensional reduction \cite[Prop.~4.6]{VV23}, followed by
transposition of stability, gives the first isomorphism.
Propositions~\ref{prop:coh} and \ref{prop:ns-crit}(c) give the next one.
\end{proof}

The following example shows that the analogous comparison is false in non symmetric type.
Let $(C,D)$ be of type $B_2$. Let $\rho$ consist of one red of vertex $2$ at centred position $L$. 
For a dimension vector $v=(v_1,v_2)$ let $n(v)$ be the number of types $\gamma$ of 
dimension $v$ with $$\Lambda^{\rm st}_{F_\gamma}\neq\emptyset$$ and let
$\tilde n(v)$ be the number of graded dimension vectors $v_\gamma$ of total
dimension $v$ with $$H^*(\widetilde{\mathfrak M}^\bullet(v_\gamma,w_\rho),\varphi_{\tilde f^\bullet_1})\neq0.$$

\begin{Proposition}\label{prop:VV23-fails} \hfill
\begin{enumerate}[label=$\mathrm{(\alph*)}$,leftmargin=8mm]
\item $n(1,1)=2$ with the  types $\{(2,L),(1,L-1)\}$, $\{(2,L),(1,L-3)\}$, and $\mathfrak L_{\rho,\gamma}=\bullet$ in both cases.
\item $n(2,1)=1$ with the type $\{(2,L),(1,L-1),(1,L-3)\}$, and $n(1,2)=0$.
\item $\tilde n(1,1)=1,$ $\tilde n(2,1)=0,$ and $\tilde n(1,2)=1$.
The nonzero critical cohomology spaces are one dimensional. 
\end{enumerate}
For any injective
assignment $\gamma\mapsto v_\gamma$ of graded dimension vectors of the same total dimension
$$H^*(\widetilde{\mathfrak M}^\bullet(v_\gamma,w_\rho),\varphi_{\tilde f^\bullet_1})
\not\cong
H^*_{P_\gamma}(\Lambda^{\rm st}_{F_\gamma},\bbC)$$
\end{Proposition}

\begin{proof}
Let $Q:1\to2$, with $d_1=2$, $d_2=1$, and one short red
at position $L$.

(a), (b) By Lemma~\ref{lem:words}, a short black is at $L$.
An arrow $B^*$ from it reaches a long black at $L-1$ or $L-3$.
A subsequent $B$ reaches a short black one position below the long one. 

For $v=(1,1)$, these give exactly the two
types in (a). Only $A$ and $B^*$ can be nonzero, and both
must be nonzero. Thus $\mu=0$ and the quotient by the
rescaling torus is a point. 

For $v=(2,1)$, both long colours
must occur, at $L-1$ and $L-3$. Again $\mu=0$, and the
three nonzero coefficients form one gauge orbit.

For $v=(1,2)$, generation requires a second short black,
reached by $BB^*$ from the first. But $A^*=0$ and this
nonzero coefficient of $BB^*$ is the sole contribution to
the corresponding moment map component, a contradiction.

(c) Give the short framing degree $k-1$. In \cite{VV23}
both arrows have degree $-2$, while
$\deg\varepsilon_1=4$, $\deg\varepsilon_2=2$ and
$\deg a_2=-1$. We use Lemma~\ref{lem:VV23-crit} and costability.

For $(1,1)$ both loops vanish, and the unique costable grading is
\[
 V_1(k+2)\xrightarrow{\alpha_{21}}V_2(k)
 \xrightarrow{a_2}W_2(k-1),
\]
with both maps nonzero. 

For $(2,1)$ we have $\varepsilon_2=0$,
the relation $\alpha_{21}\varepsilon_1=0$ makes
$\operatorname{Im}\varepsilon_1$ an unframed submodule.
Costability gives $\varepsilon_1=0$, after which the
nonzero kernel of $V_1\to V_2$ is an unframed submodule.
Thus there are no critical points.

For $(1,2)$ we have $\varepsilon_1=0$. If
$\varepsilon_2 e=f\ne0$, with $e,f$ in degrees $r,r+2$,
then $a_2(f)=0$. Costability requires
$\alpha_{12}(f)\ne0$, placing the long line in degree
$r$. Hence $\alpha_{21}=0$ by degree, leaving that line
unframed, a contradiction. Therefore $\varepsilon_2=0$.
The unique costable grading is
\[
 V_2(k+4)\xrightarrow{\alpha_{12}}V_1(k+2)
 \xrightarrow{\alpha_{21}}V_2(k)
 \xrightarrow{a_2}W_2(k-1),
\]
with all maps nonzero. Each of the two surviving graded
quotients is a point and its potential is identically zero.
The unshifted vanishing cycles give a one dimensional
cohomology up to shift. 

Exchanging the framed vertex and
transposing the Cartan convention gives the same $(1,1)$ count.

To prove the last claim, note that in (a) there are two nonempty types of dimension $(1,1)$,
while, in  (c), there is only one nonzero component of that dimension. 
\end{proof}


\section{Microlocal examples}
\label{app:examples}

\subsection{Non surjectivity of the restriction map \eqref{MM}}
\label{ex:restriction-nonsurjective}
Let $Q$ be the Kronecker quiver $1\rightrightarrows2$ with $D=1$.
Let 
$$\dV_{(1,0)}=\dV_{(2,1)}=\bbC^2
,\quad
\dW_{(1,0)}=\bbC^2,$$
 with all other colour spaces zero.
Thus $\dN$ has coordinates
\[
 A:\bbC^2\to\bbC^2,\quad
 B_1,B_2:\bbC^2\to\bbC^2.
\]
Let $\gamma$ have two blacks at centred positions $(1,0)$ and
$(2,-1)$, and $\eta$ have one black at each of $(1,0),(1,-2),(2,-1),(2,-3)$.
The flag of type $\gamma$ is unique and $\scrX_\gamma=\dN$.
Its stable locus is the open set $U\subset\dN$ where $A$ is
invertible and $B_1(V_1)+B_2(V_1)=V_2$. Hence
$$Z_\gamma^{\rm st}=U\times\{0\}.$$

A flag of type $\eta$ amounts to lines $L_i\subset V_i$ with
$B_1L_1,B_2L_1\subset L_2$. Therefore $\scrX_\eta$ is a
vector bundle over $\mathbb P^1\times\mathbb P^1$ with fibres
of dimension $10$. In particular, it is  irreducible and 
$$d_\eta=\dim\dN=n=12.$$
Since $\mathscr L_\gamma=\bbC_{\dN}[2n]$, the zero section
formula for microlocal Hom identifies the $\eta$-piece of
\eqref{MM} with the restriction
\[
 H^*_{\dG}\bigl(\pi_\eta^{-1}(U),\bbC\bigr)
 \longrightarrow H^*\bigl(\pi_\eta^{-1}(x),\bbC\bigr),
 \quad \xi=(x,0).
\]
Choose $A=B_1=I_2$ and $B_2=\operatorname{diag}(0,1)$.
The corresponding point  $x$ lies in $U$.
Thus, we have $$\xi\in\Lambda_{F_\gamma}^{\rm st,\bullet}.$$
We have $L_2=L_1$, while $L_1$ must be an eigenline of $B_2$. 
Thus $\pi_\eta^{-1}(x)$ consists of two points.
On the other hand, $\pi_\eta^{-1}(U)$ is a nonempty open subset
of the irreducible variety $\scrX_\eta$. 
Hence it is connected.
In degree zero the restriction map \eqref{MM} is therefore
\[
 e_\eta\operatorname{res}_\xi:\bbC\longrightarrow\bbC^2,
 \quad 1\longmapsto(1,1).
\]
It is not surjective.

\subsection{The costalk-to-stalk morphism for Morse modules has not a simple image}
\label{ex:coMorse-A1-nonsimple}
Let $Q=A_1$, $V=\bbC$ and
$W=\bbC w_0\oplus\bbC w_2$, with reds at positions $0,2$.
Let $\gamma$ have one black at $2$ and $\eta$ one at $0$.
Write 
$$(a_0,a_2)\in\dN=\Hom(W,V)=\bbC^2
,\quad
(b_0,b_2)\in\dN^*.$$ 
Then, we have
\[
 \scrX_\gamma=\{a_0=0\},\quad \scrX_\eta=\dN,\quad
 Z_\gamma^{\rm st}=\{a_0=b_2=0,\ a_2\ne0\}.
\]
The off diagonal pairing is
$$f_{\eta,\gamma}=b_0u_0+b_2(u_2-v_2).$$
We have
\begin{equation}\label{eq:offdiag-A1}
 \left.\muhom(\mathscr L_\gamma,\mathscr L_\eta)\right|_{\Omega_\gamma}
 \cong\bbC_{\{a_0=0,\ a_2\ne0,\ b_0=b_2=0\}},
\end{equation}
extended by zero in $\Omega_\gamma$.
Write 
$$Y=\scrX_\gamma=\{a_0=0\}\subset\dN.$$
Choose the stable pair $\xi=(x,y)$ with
\[
 x=(0,1),\quad y=(0,0).
\]
Choose the Morse  function
\[
 f(a_0,a_2)=\tfrac12\bigl(a_0^2+(a_2-1)^2\bigr).
\]
Its differential at $x$ is $y$.
Its graph is transverse at
$\xi$ to the conormal Lagrangian $T_Y^*\dN$.
Put $\mathcal P_\delta=\mathscr L_\delta[-d_\delta]$, as in Remark~\ref{rem:morse-evaluation-dual}(c). For the two indicated
types,
\[
 \mathcal P_\gamma=\bbC_Y[1],\quad
 \mathcal P_\eta=\bbC_{\dN}[2],
\]
where $\bbC_Y$ is extended by zero to $\dN$.
These are nonisomorphic simple perverse sheaves. The restrictions
of $f$ to $Y$ and $\dN$ both have a unique nondegenerate critical
point at $x$. Hence, in a neighbourhood of $x$, we have
\[
 \phi_f(\mathcal P_\gamma)\cong\bbC_{\{x\}},\quad
 \phi_f(\mathcal P_\eta)\cong\bbC_{\{x\}}.
\]
All nonzero Morse components have these skyscraper vanishing cycles.
The remaining incidence space $\{0\}$ misses $x$.
Hence, the costalk-to-stalk morphism is an isomorphism
\[
 \mathring{\mathcal I}_\xi:\mathring M^!(\rho,\gamma)^\xi
 \xrightarrow{\cong}\mathring M(\rho,\gamma)^\xi.
\]
Nevertheless, its image is not simple. Indeed, the normalized
Morse fibres are
\[
 MF^\xi_\gamma\cong\bbC,\quad MF^\xi_\eta\cong\bbC[1].
\]
So, we have
$$H^0(MF^\xi_\gamma)=\bbC
,\quad
H^{-1}(MF^\xi_\eta)=\bbC.$$
We claim that these two contributions cannot be connected by the $\bR$-action.
Indeed, consider a class
$a\in e_\beta\bR^k e_\alpha$ between types whose Morse fibres
are nonzero. This class yields a morphism
\[
 \mathcal P_\alpha\longrightarrow
 \mathcal P_\beta[k+d_\beta-d_\alpha].
\]
Since the vanishing cycles of both perverse sheaves are
concentrated in degree zero, the induced map can be nonzero only
when $k+d_\beta-d_\alpha=0$. In that case it comes from an element in
$\operatorname{Hom}(\mathcal P_\alpha,\mathcal P_\beta)$.
In particular, it vanishes if one of these perverse sheaves is supported on
$Y$ and the other on $\dN$. All nonzero Morse components arise
from one of these two perverse sheaves. Their corresponding
isotypic summands are therefore nonzero proper $\bR$-submodules.
Thus
\[
 \operatorname{Im}\mathring{\mathcal I}_\xi
 =\mathring M(\rho,\gamma)^\xi
\]
is not simple.

\subsection{The costalk-to-stalk morphism for Microlocal modules has a simple image}
\label{subsubsec:A1-slice-image}
We keep the notation of \S\ref{ex:coMorse-A1-nonsimple},
with the homogeneous element $\xi_0=\xi=((0,1),(0,0))$.
Put  
$$S=G\xi_0=\{b_0=0\}\subset Z_\gamma^{\rm st}.$$
The free $G=\bbC^\times$-action sends $(a_2,b_0)$ to
$(ga_2,g^{-1}b_0)$. Thus $t=a_2b_0$ gives
\[
 X_\gamma=Z_\gamma^{\rm st}/G\cong\mathbb A^1,\quad
 \bar\xi_0=t(\xi_0)=0,\quad S/G=\{0\}.
\]
A local section of $Z_\gamma^{\rm st}\to X_\gamma$ identifies a neighbourhood
of $0$ in the quotient with an analytic transverse slice to
$S$ in $Z_\gamma^{\rm st}$.
Let $\mathcal F_s$ denote the descent of the sheaf
$\mathcal K_s=\muhom(\mathscr L_\gamma,\mathscr L_s)|_Z$ to $X_\gamma$.
Here $s$ is an even black position.  By Proposition~\ref{prop:psi}(a) and
\eqref{eq:offdiag-A1}, since $\scrX_s=\dN$ for $s\leqslant 0$ and $\scrX_s=\{0\}$ for
$s\geqslant 4$, these complexes are
\[
 \mathcal F_2\cong\bbC_{\mathbb A^1},\quad
 \mathcal F_s\cong\bbC_{\{0\}},\quad
 \mathcal F_t=0,\quad s\leqslant 0,\quad 
 4\leqslant t.
\]
Here $\bbC_{\{0\}}$ is extended by zero.
Given a small analytic
disc $T$ in $\mathbb A^1$ around $0$, let
$i_{\bar\xi}:\{0\}\hookrightarrow T$ be the inclusion. Since the slice dimension
is $d=1$, the normalized source and the shifted target of
\eqref{eq:co-specialized-slice-map} are given by
\[
 \begin{array}{c|cc}
 &(i_{\bar\xi})^!\mathcal F_s[2]&(i_{\bar\xi})^*\mathcal F_s[2]\\ \hline
 s=2 &\bbC&\bbC[2]\\
 s\leqslant0 &\bbC[2]&\bbC[2]
 \end{array}
\]
while both terms vanish for $s\geqslant4$. The canonical map on
the $s=2$ component is zero for degree reasons, whereas it is
the identity for each $s\leqslant0$. Thus the normalized
co-specialized $2$-piece is in degree $0$, the $s\leqslant0$
co-specialized pieces are in degree $-2$, and the corresponding
unshifted specialized pieces are in degree $0$. Therefore, the image of the map \eqref{eq:co-specialized-slice-map} is
\[
 \operatorname{Im}\mathcal I_{\xi_0}
 =\bigoplus_{s\leqslant0,\ s\in2\bbZ}
 e_s M(\rho,\gamma)^{\xi_0},
 \quad e_s M(\rho,\gamma)^{\xi_0}\cong\bbC.
\]
This image is simple: for even $s,t\leqslant0$, we have
$$\mathscr L_s=\mathscr L_t=\bbC_{\dN}[4],$$and the identity gives a
degree zero class in $e_s\bR e_t$ identifying their one dimensional
pieces. Any nonzero submodule therefore contains them all.
Thus the full costalk-to-stalk map is zero, its sliced
version has simple image, and the co-Morse version is an
isomorphism onto a nonsimple module. 

\subsection{The map \eqref{eq:microlocal-standard-map} is surjective}
\label{subsec:A1-gysin}
We keep the data of \S\ref{ex:coMorse-A1-nonsimple}.
Put
$$X=\dN=\bbC^2
,\quad
Y=\{a_0=0\}
,\quad
i:Y\hookrightarrow X.$$
Then, we have 
$$Z_\gamma^{\rm st}/{\dG}\cong\mathbb A^1
,\quad
H_{\dG}^*(Z_\gamma^{\rm st},\bbC)=e_\gamma M(\rho,\gamma)=\bbC.$$
The isomorphism \eqref{eq:offdiag-A1} also gives $e_\eta M(\rho,\gamma)=\bbC$.
Since 
$$\mathscr L_\gamma=i_*\bbC_Y[2],\quad
\mathscr L_\eta=\bbC_X[4],$$
the purity gives
\[
 e_\eta\bR^0e_\gamma
 =\Hom_{D_{\dG}^b(X)}(i_*\bbC_Y[2],\bbC_X[4])
 \cong H_{\dG}^0(Y,\bbC)=\bbC.
\]
Let $a$ be its Gysin generator. 
The microlocal composition with $a$ yields a map
$\bbC_{Z_\gamma^{\rm st}}\to\bbC_S$, where
$$S=\{a_0=b_0=b_2=0,\ a_2\ne0\}.$$
At $\xi_0=((0,1),(0,0))$ this morphism is nonzero. After rescaling $a$, it is
$1\mapsto1$. Since $S/\dG$ is the origin of $Z_\gamma^{\rm st}/\dG$,
its equivariant cohomology map is also $1\mapsto1$. Thus $a1_\gamma\ne0$.
For every even black position $t\leqslant0$ the same argument applies,
since $\mathscr L_t=\bbC_X[4]$ and $e_tM(\rho,\gamma)=\bbC$.
For $t\geqslant4$, we have $\scrX_t=\{0\}$ and the corresponding piece vanishes.
Hence $M(\rho,\gamma)=\bR1_\gamma$.
We deduce that, in this example, the map
\eqref{eq:microlocal-standard-map} is surjective.

\subsection{The map \eqref{eq:micro-to-Morse} is not an isomorphism}
\label{sec:micro-to-Morse}
We keep the data of \S\ref{ex:coMorse-A1-nonsimple}.
At the point $\xi_0=\xi$, the Morse function of \S\ref{ex:coMorse-A1-nonsimple} gives
$$MF^\xi_\gamma\cong\bbC,\quad
MF^\xi_\eta\cong\bbC[1],$$
whereas
\eqref{eq:offdiag-A1} gives
$$\muhom(\mathscr L_\gamma,\mathscr L_\eta)_{\xi_0}\cong\bbC.$$
Thus the comparison map $\theta_{\xi_0,f}$ in
\eqref{eq:micro-to-Morse} is an isomorphism on the $\gamma$-piece.
It kills the nonzero degree zero $\eta$-class and misses the degree $-1$ Morse class. 
Thus it is neither injective nor surjective.

\subsection{Non generic Morse and specialized microlocal modules have different classes}
\label{ex:A1-K0-special}
In type $A_1$, take 
$$V=W(0)=W(2)=\bbC^2$$ and the centred types
$\gamma=(2,2)$, $\eta=(0,2)$. Write 
$$A=(A_0,A_2)
,\quad
A^*=(J_0,J_2).$$ 
The point $\xi_0=((0,I_2),(0,0))$ is homogeneous
and stable: $A_2$ generates $V$ at position $2$. Put
$$D=\{\det A_0=0\}\subset\operatorname{Mat}_2(\bbC).$$
The incidence spaces are
\[
 \scrX_\gamma=\{0\}\times\operatorname{Mat}_2(\bbC),\quad
 \scrX_\eta=\widetilde D\times\operatorname{Mat}_2(\bbC),
 \quad d_\gamma=4,\quad d_\eta=7,
\]
where 
$$\sigma:\widetilde D=\operatorname{Tot}
(\calO_{\mathbb P^1}(-1)^{\oplus2})\to D$$ is the small resolution
recording a line containing $\operatorname{Im}A_0$.
At the zero covector, microlocal Hom is  sheaf Hom.
Purity for the zero section of $\widetilde D$ and proper base change give
\[
 \muhom(\mathscr L_\gamma,\mathscr L_\eta)_{\xi_0}
 \cong R\Gamma(\mathbb P^1,\bbC)[2].
\]
Choose the Morse function
$$f=q(A_0)+r(A_2-I_2),$$
with $q$ generic quadratic and $r$
nondegenerate quadratic. The Milnor fibre of $q|_D$ has the homotopy type of a double cover
of $\mathbb P(D)\cong\mathbb P^1\times\mathbb P^1$ minus the smooth
$(2,2)$ divisor $\{q=0\}$. It is of genus one. Its Euler characteristic is $8$.
For $\mathcal P=\IC_D=R\sigma_*\bbC_{\widetilde D}[3]$, we have
$$\chi(\mathcal P_0)=-2
,\quad
\chi((\psi_q\mathcal P)_0)=-8,$$
where $\psi_q$ denotes nearby cycles.
The normalized vanishing cycles $\phi_q\mathcal P=\varphi_q\mathcal P[-1]$ are perverse
and supported at $0$. Their rank is $6$.
Thom-Sebastiani therefore gives
$$\varphi_f(\mathscr L_\eta)_{(0,I_2)}[-d_\gamma-1]\cong\bbC^6[3].$$
Consequently, we have
\[
 \dim e_\eta M(\rho,\gamma)^{\xi_0}=2,\quad
 \dim e_\eta\mathring M(\rho,\gamma)^{\xi_0}=6.
\]
Thus, these modules have distinct classes in $K_0(\scrO_\mu^\rho)$.

\subsection{The local Kirwan map is not surjective}
\label{subsec:B2-local-kirwan}
The local Kirwan map need not be onto.
We give an example in type $B_2$.

\begin{Proposition}\label{prop:B2-local-kirwan}
Let $Q:1\to 2$, $D=\operatorname{diag}(2,1)$ and $v=(1,2),$ $w=(0,3).$
Choose the centred positions
\[
 V_1(-2)=\bbC,\quad V_2(1)=V_2(-3)=\bbC,
 \quad W_2(1)=W_2(-1)=W_2(-5)=\bbC.
\]
All other graded pieces are zero. The datum is $\ell$-dominant with
\[
 Z_{\rho,\gamma}=Z_{1,-2}Z_{1,2},\quad
 \mathfrak L_{\rho,\gamma}
       \cong\bbC^\times\times\mathbb A^2.
\]
We have
$\operatorname{Im}\Phi^{\rm ns}_\gamma
  =\operatorname{Im}\kappa^{\rm ns}_\gamma=\bbC
  \subsetneq H^*_{P_\gamma}(\Lambda^{\rm st}_{F_\gamma},\bbC).$
\end{Proposition}

\begin{proof}
Let $L=V_1(-2)$ and write $V_2=\bbC e_1\oplus\bbC e_{-3}$.
The two short blacks and the three reds have the same colour.
The long black has colour $2\mod4$.
The allowed matrices in
$L_{F_\gamma}$ are
\begin{equation}\label{eq:B2-local-matrices}
 A=\begin{pmatrix}p&0&0\\q&a&0\end{pmatrix},\quad
 A^*=\begin{pmatrix}0&0\\b&0\\c&d\end{pmatrix},\quad
 B=\binom{0}{y},\quad B^*=(x,0),
\end{equation}
where the framing basis is ordered by the positions $1,-1,-5$.
Indeed, $A$ weakly lowers positions, $A^*$ strictly lowers them,
and both edge maps strictly lower positions between these colours.
The two moment map components are
\[
 BB^*+AA^*=\begin{pmatrix}0&0\\xy+ab&0\end{pmatrix},
 \quad -B^*B=0.
\]
Consequently the only equation is $$xy+ab=0.$$
A short computation shows  that the minimal flag condition is
\begin{equation}\label{eq:B2-local-stability}
 pxy\ne 0.
\end{equation}
The group $P_\gamma$ consists of the lower triangular invertible
matrices on $V_2$, together with the scalars on $L$.
At each stable point the vectors
\[
 v=A(r_1),\quad \ell=B^*v,\quad w=B\ell
\]
give a basis $\ell$ of $L$ and a basis $(v,w)$ of $V_2$ adapted
to $F_\gamma$. They therefore determine a unique gauge
transformation to the normal form
\begin{equation}\label{eq:B2-local-normal-form}
 A=\begin{pmatrix}1&0&0\\0&-z^{-1}&0\end{pmatrix},\quad
 A^*=\begin{pmatrix}0&0\\z&0\\u&v'\end{pmatrix},\quad
 B=\binom{0}{1},\quad B^*=(1,0).
\end{equation}
Here $(z,u,v')\in\bbC^\times\times\mathbb A^2$.
The moment equation gives $ab=-1$ in this frame. Conversely,
every point of \eqref{eq:B2-local-normal-form} is stable and has zero moment map. 
The frame and its inverse depend regularly on the stable point, since $pxy$ is invertible. 
Thus, this is a global algebraic slice, proving that
$\mathfrak L_{\rho,\gamma}\cong\bbC^\times\times\mathbb A^2$.
Finally, we have
\[
 H^*\big(\bbC^\times\times\mathbb A^2,\bbC\big)
     =\bigwedge\nolimits_{\bbC}(\eta),\quad \deg\eta=1.
\]
Since $H^*_{P_\gamma}$ is concentrated in even degrees,
its image is $\bbC$. Equivalently, the global slice
trivializes all the graded tautological bundles. Theorem
\ref{thm:ns-map}(b), together with the surjectivity of the map $r_\gamma$, 
identifies this image with that of $\Phi^{\rm ns}_\gamma$, which therefore misses $\eta$.
\end{proof}

Thus the non symmetric local comparison map is onto its tautological image.
For the global comparison and its relation to Kirwan
surjectivity, see Proposition~\ref{prop:ns-global-torsion}.

\subsection{The global Kirwan map is not surjective}
\label{subsec:G2-global-kirwan}
The global Kirwan map need not be onto.
We give an example in type $G_2$. 

\begin{Proposition}\label{prop:G2-global-kirwan}
Let $Q:1\to 2$, $D=\operatorname{diag}(3,1)$ and
$v=(1,1),$ $w=(0,1).$
Then $N_0=3$ and $Y=\mathfrak M_3(v,w)$ is a smooth complex
fourfold. Its Poincar\'e polynomial is
\[
 \sum_j\dim H^j(Y,\bbC)t^j=1+2t^2+t^3+t^4.
\]
The Kirwan map
$H^*_{G_3}\longrightarrow H^*(Y,\bbC)$
is not surjective. \end{Proposition}

\begin{proof}
Put
$h=c_1(\mathcal E_2)$ and $k=c_1(\mathcal E_1)$, where
$\mathcal E_i$ is the coefficient bundle of \eqref{EI}. 
We will prove that the image of the Kirwan map is
\[
 H^{\rm even}(Y,\bbC)
 \cong\bbC[h,k]/(h^2,hk,k^3),\quad \deg h=\deg k=2.
\]
The transpose Cartan convention gives
\[
 H_1=\bbC,\quad H_2=H=\bbC[\varepsilon]/(\varepsilon^3),
 \quad H_{12}=\bbC.
\]
The gauge modules at vertices $1,2$ are $L=\bbC$ and $S=H$,
and the free framing module $H\otimes W_2$ is $H$.
Thus the doubled data are
\[
 a:H\to S,\quad c:S\to H\quad\text{($H$-linear)},\quad
 u:L\to S,\quad f:S\to L\quad\text{($\bbC$-linear)}.
\]
We identify $a,c,u$ with elements of $H$ and write $f\in H^*$.
The gauge group $G_3=\bbC^\times\times H^\times$ acts by
\[
 (t,g)\cdot(a,c,u,f)=(ga,cg^{-1},gu\,t^{-1},tfg^{-1}).
\]
Use $\operatorname{tr}_H(h)=[\varepsilon^2]h$ for the framing
duality. Define $\beta(u,f)\in H$ by
\[
 \operatorname{tr}_H(b\beta(u,f))=f(bu)\quad(b\in H).
\]
For $u=u_0+u_1\varepsilon+u_2\varepsilon^2$ and
$q_j=f(\varepsilon^j)$, this is
\[
 \beta(u,f)=u_0q_2+(u_0q_1+u_1q_2)\varepsilon
       +(u_0q_0+u_1q_1+u_2q_2)\varepsilon^2.
\]
Indeed, the moment map pairing with $(\ell,b)\in\bbC\oplus H$ is
\[
 \langle\mu_3,(\ell,b)\rangle
   =-\ell f(u)+\operatorname{tr}_H(bac)+f(bu).
\]
Consequently, the zero-moment equations are
\begin{equation}\label{eq:G2-kirwan-moment}
 f(u)=0,\quad ac+\beta(u,f)=0.
\end{equation}
The generation stability is equivalent to
\begin{equation}\label{eq:G2-kirwan-stability}
 f(Ha)\ne0,\quad Ha+Hu=H.
\end{equation}
In fact, the framing first generates $Ha$ at vertex $2$.
Vertex $1$ is reached precisely when $f(Ha)\ne0$, after which
$u$ and the coefficient loop generate $Hu$ as well. The generated
submodule at vertex $2$ is therefore  $Ha+Hu$.
The second condition in \eqref{eq:G2-kirwan-stability} is
$a_0\ne0$ or $u_0\ne0$. The stable gauge action is free, and
Proposition~\ref{prop:MN} gives
$\dim Y=4$ and smoothness.

Let $U\subset Y$ be the locus $a_0\ne0$. Normalize $a=1$
using $H^\times$. Then $c=-\beta(u,f)$, stability says $f\ne0$,
and the residual action is $(u,f)\mapsto(t^{-1}u,tf)$.
Hence
\begin{equation}\label{eq:G2-kirwan-open}
 U\cong\{(u,f)\in H\times H^*:f\ne0,\ f(u)=0\}/\bbC^\times
   \cong T^*\mathbb P^2.
\end{equation}
The map $f:H\twoheadrightarrow L$ identifies
$\mathcal E_1|_U$ with the universal quotient line on
$\mathbb P(H^*)$, pulled back to its cotangent bundle. Thus
\begin{equation}\label{eq:G2-kirwan-open-cohomology}
 H^*(U,\bbC)=\bbC[k|_U]/((k|_U)^3).
\end{equation}
In particular, the generator of $H^2(U,\bbC)$ extends over $Y$.

Put $D_0=Y\smallsetminus U$. Write
$a=a_0+a_1\varepsilon+a_2\varepsilon^2$.
If $a_0=0$, stability forces $u_0\ne0$, and the constant
coefficient of \eqref{eq:G2-kirwan-moment} gives $q_2=0$.
If $a_1=0$, then $Ha\subseteq\bbC\varepsilon^2$ and $f(Ha)=0$,
a contradiction. Thus $a_1\ne0$ and $q_1\ne0$.
Normalize $u=1$ and the coefficient of $\varepsilon$ in $a$ to $1$.
This fixes the gauge completely. The equations then give the
unique normal form
\begin{equation}\label{eq:G2-kirwan-boundary}
 \begin{gathered}
 u=1,\quad a=\varepsilon+z\varepsilon^2,\quad
 c=r-zr\varepsilon+s\varepsilon^2,\\
 q_0=q_2=0,\quad q_1=-r,\quad
 (z,s,r)\in\mathbb A^2\times\bbC^\times.
 \end{gathered}
\end{equation}
Conversely, these data are stable, since $f(a)=-r\ne0$ and
$u=1$. Thus, we have
\[
 D_0\cong\mathbb A^2\times\bbC^\times.
\]
Concerning the scheme structure, one can check that $D_0$ is a smooth closed Cartier divisor.
It is the zero divisor of the section of $\mathcal E_2$ induced
by the incoming framing modulo $\varepsilon$, so $[D_0]=h$.

For $j:U\hookrightarrow Y$ and $i:D_0\hookrightarrow Y$,
the Gysin sequence contains
\[
 H^2(Y,\bbC)\xrightarrow{j^*}H^2(U,\bbC)
 \longrightarrow H^1(D_0,\bbC)
 \xrightarrow{i_*}H^3(Y,\bbC)\longrightarrow H^3(U,\bbC).
\]
The first arrow is surjective by
\eqref{eq:G2-kirwan-open-cohomology}, and $H^3(U,\bbC)=0$.
Since $H^1(D_0,\bbC)=\bbC$, exactness gives
\begin{equation}\label{eq:G2-kirwan-odd}
 H^3(Y,\bbC)\cong H^1(D_0,\bbC)\cong\bbC.
\end{equation}
The other terms give
\[
 0\longrightarrow\bbC[D_0]\longrightarrow H^2(Y,\bbC)
 \longrightarrow\bbC\,k|_U\longrightarrow0,
 \quad H^4(Y,\bbC)\xrightarrow{\cong}H^4(U,\bbC),
\]
and $H^j(Y,\bbC)=0$ for $j=1$ and $j>4$.
Thus $h,k$ form a basis in degree $2$, and $k^2$ generates
$H^4(Y,\bbC)$. Since $h|_U=0$ and restriction is injective in
degree $4$, we have $h^2=hk=0$. This proves the asserted
Poincar\'e polynomial and even cohomology ring.

Finally, recall that $H^\times$ has unipotent radical $1+\varepsilon H$
and reductive quotient $\bbC^\times$. Hence
\[
 H^*_{G_3}=\bbC[c_1,c_2],\quad \deg c_1=\deg c_2=2.
\]
The Kirwan map sends $c_1$ to $k$ and $c_2$ to $h$.
Its image is therefore $H^{\rm even}(Y,\bbC)$ and misses
\eqref{eq:G2-kirwan-odd}.
\end{proof}

\medskip

\centerline{\textbf{AI disclosure}}

\smallskip

The paper was written by us and all ideas are ours. AI (ChatGPT) 
was used for proofreading purposes, to improve the presentation, and for the examples in Appendix C.

\end{document}